\documentclass[letterpaper,11pt]{article}

\usepackage{latexsym, graphicx, epsfig, amsmath, amsfonts,amssymb,bm}
\usepackage{caption, subcaption}
\usepackage{epstopdf}
\usepackage[margin=1in,vmargin=2cm,includefoot]{geometry}

\usepackage{array}
\usepackage{color}
\usepackage{url}

\def\be{\begin{equation}}
\def\ee{\end{equation}}

\def\y{\mathbf{y}}

\def\PPhi{\boldsymbol{\Phi}}

\def\PPh{\mathbf{\Phi}}

\def\R{{\mathbb R}}

\def\u{\mathbf{u}}

\newcommand{\figref}[1]{Fig.~\ref{#1}}

\usepackage[normalem]{ulem}

\newtheorem{remark}{Remark}[section]
\newtheorem{example}{Example}[section]

\title{Designing Neural Networks for Hyperbolic Conservation Laws}

\author{Zhen Chen\footnotemark[1]\and Anne Gelb\thanks{Department of Mathematics,
		Dartmouth College, Hanover, NH 03755, USA. Emails:
		{\tt zhen.chen@dartmouth.edu, annegelb@math.dartmouth.edu, Yoonsang.lee@dartmouth.edu} Funding: All authors are supported by the DoD MURI grant ONR \#N00014-20-1-2595. AG is also supported by the AFOSR grant \#FA9550-22-1-0411 and the NSF grant DMS \#1912685. YL is also supported by the NSF grant DMS \#1912999. }\and Yoonsang Lee\footnotemark[1]
}

\date{}
\begin{document}
	
	\maketitle
	
	\begin{abstract}
		We propose a new data-driven method to learn the dynamics of an unknown hyperbolic system of conservation laws using deep neural networks. Inspired by classical methods in numerical conservation laws, we develop a new conservative form network (CFN) in which the network learns the flux function of the unknown system. Our numerical examples demonstrate that the CFN yields significantly better prediction accuracy than what is obtained using a standard non-conservative form network, even when it is enhanced with constraints to promote conservation. In particular, solutions obtained using the CFN consistently capture the correct shock propagation speed without introducing non-physical oscillations into the solution. They are furthermore robust to noisy and sparse observation environments.
	\end{abstract}
\section{Introduction}

Hyperbolic systems of conservation laws arise in many applications, in particular where wave motion or advective transport is important, and include problems in gas dynamics, acoustics, optics and elastodynamics. These typically non-linear partial differential equations (PDEs) are well known to cause challenges for numerical simulation. Methods that work well for linear problems break down near discontinuities in the non-linear case unless special care is taken to mitigate the oscillatory behavior that will otherwise cause instabilities. Moreover it is well understood that numerical solvers must be written in {\em flux conserving form}, since otherwise the resulting numerical solution may yield the incorrect shock speed, \cite{hesthaven2017numerical, leveque1992numerical, leveque2002finite}.

Since enormous quantities of data can now be collected, and since there is increased capacity in computational storage and efficiency, data-driven algorithms that recover unknown dynamical systems  from observation data, mainly for ordinary differential equations (ODEs) but more recently for PDES as well,  are becoming more widespread. Many data-driven methods build neural network (NN) models that are trained to fit observed data.   Such an approach is called ``blind'' since it does not take into account any information regarding the underlying system.   As will be demonstrated in our numerical results, while a standard NN technique may yield accurate results within the training time domain, numerical difficulties arise once shocks are formed beyond that time period.  Indeed, the NN may predict non-physical solutions or the wrong shock propagation speed.

Some ODE dynamics-learning methods have been developed to impose physical constraints that may be known about the underlying system. For instance regularization terms can be introduced into the loss function to penalize the NN that would otherwise not satisfy physical constraints \cite{cranmer2020lagrangian, greydanus2019hamiltonian}. Another typical approach seeks to integrate general physical principles into the design of neural networks. For example, the method in \cite{zhang2022gfinns} embeds GENERIC formalism into neural networks to learn dynamical systems and involves two separate generators for reversible and irreversible dynamics, respectively.  Systems of ODEs and PDEs are considered in \cite{baddoo2021physics}, where general physical principles are integrated into what is called {\em dynamic mode decomposition} (DMD), which learns low-rank dynamics from high-dimensional measurements.  This case-by-case study does not have a general framework that can be extended to conservation laws, however.  
There is a also class of data-driven methods that aim to identify the governing PDE equations, see e.g.~\cite{kutz2016SINDy,schaeffer2017learning}.  The time derivative must be approximated from the observation data, thereby limiting observation environments to those where high quality data are obtainable in very short time intervals.
Yet another class of data-driven methods described in \cite{raissi2019physics,rudy2017data} assume that the governing PDE is known, and have the goal to either learn the long term solution or the PDE parameters.  Finally, it is important to point out that all of these methods assume that the solution remains smooth and differentiable in the considered time domain.

In spite of their ubiquitousness in many applications, to the best of our knowledge there are no data-driven methods that are specifically designed to learn unknown systems of hyperbolic conservation laws. In this investigation we therefore seek to develop a method to do so, and in particular for the case in which  the observation data are comprised of a set of perturbed numerical solutions of the true PDE at a finite number of different time instances within a short time (training) domain.   Our new method employs tools from machine learning (ML) to construct a model from the training data to make long-term predictions for a system of hyperbolic conservation laws that extend beyond the short-time domain for which data are observed.

Based on similar behavior observed in traditional solvers for numerical conservation laws, here we propose a new {\em conservative form network} (CFN) for which the network learns the flux function of the unknown system. By incorporating the conservative flux form directly into the network architecture, we are able to mimic the structure of the conservative form scheme found in classical numerical hyperbolic conservation laws. Our new method is distinguishable from these other mentioned approaches as it designs the network to be specifically in conservative form.  Our numerical experiments demonstrate that the data-driven method resulting from our new CFN is conservative. It furthermore correctly predicts the shock propagation speed in extended time domains without introducing non-physical oscillations into the solution.  This is in contrast to data-driven methods that use either the standard ``blind'' NN or those that incorporate a penaly term to promote conservation.

The rest of the paper is organized as follows. Section \ref{sec:setup} gives a brief review of conservation laws along with some standard ideas related to dynamics-learning methods. We introduce our new conservative form network in Section \ref{sec:method}. Section \ref{sec:numexamples} discusses how our experiments are designed. This is followed by three classical examples of one-dimensional conservation laws in Section \ref{sec:examples} used to validate our approach.  Section \ref{sec:conclusion} provides some concluding remarks.

\section{Preliminaries}
\label{sec:setup}
\subsection{Conservation laws}
\label{sec:conserv}
We are interested in learning the dynamics of an unknown hyperbolic system of conservation laws on a spatial domain $x \in (a,b)$ and temporal domain $t \in (0,T)$.   The scalar form is given by 
\begin{equation}\label{setup:eqn}
	 u(x,t)_t + f(u(x,t))_x = 0
\end{equation}
{with appropriate initial and boundary conditions.}\footnote{{We will consider systems of conservation laws in  our numerical experiments.}} 
The main difficulty in solving \eqref{setup:eqn} is due to the formation of shock discontinuities, which will occur even when the initial conditions are smooth.  To retain the proper shock speed, numerical solvers for \eqref{setup:eqn} must be written in flux conserving form \cite{leveque2002finite, leveque1992numerical,hesthaven2017numerical}. 
Specifically, if the spatial domain is discretized as $x_j$, $j = 0,\dots,N$, and time is incremented at $t = t_l$, numerical solvers for the interior of the interval should be of the form
\begin{equation}\label{setup:law}
	\int_{x_j}^{x_{j+1}} u(x,t_{l+1}) dx = \int_{x_j}^{x_{j+1}} u(x,t_l) dx + \int_{t_l}^{t_{l+1}}f(u(x_j,t))dt - \int_{t_l}^{t_{l+1}} f(u(x_{j+1},t))dt.
\end{equation}
In this investigation we seek to approximate the solution 
${\mathbf  u}(t) = \{\bar{u}_j(t)\}_{j = 0}^N$ for $t \in (0,T)$,
where $\bar{u}_j(t)$ is the cell average over a uniform grid $\{x_j\}_{j=1}^{N}$  given by
\begin{equation}
	\bar{u}_j(t) = \int_{x_j-\frac{\Delta x}{2}}^{x_j+\frac{\Delta x}{2}} u(x,t)dx,\quad j=1,\dots, N,\quad \Delta x = \frac{b-a}{N}.
\end{equation}

Studies regarding numerical conservation laws typically assume the flux term is known, with the goal to construct accurate, robust, and efficient solvers for ${\mathbf u}(t)$ by appropriately discretizing \eqref{setup:law}.  Here, by contrast, we are interested in the case where we know apriori that the governing equation is a conservation law, but the flux function itself is unknown.  
The goal is then to determine how the solution will evolve given some early observations regarding the governing PDE.   We will exploit our understanding of conservation laws by designing our numerical method to be in conservative form, as is given by \eqref{setup:law}. Once we are able to construct the numerical flux,  it can be used to predict the evolution of the unknown PDE.

We now introduce some notation for the observable data. We will assume that the solution to the PDE is available at a set of discrete time instances $\{t_l\}_{l=1}^{L}$ for $N_{traj}$ initial conditions resulting in so-called snapshots of the solution,
\begin{equation}\label{setup:data}
	\mathbf{u}^{(k)}(t_l), \quad l=1,\dots,L, \quad k = 1,\dots, N_{traj}.
\end{equation}
The superscript $k$ in \eqref{setup:data} denotes the $k$-th ``trajectory'', which implies all $L$ snapshots are evolved from the same initial state, with $N_{traj}$ denoting the total number of trajectories. The $N_{traj}$ initial conditions in our experiments are chosen by perturbing the true initial conditions of the PDE.
The time step between two consecutive time instances is given by
\begin{equation*}
	\Delta t = t_{l+1} - t_l, \quad l=1,\dots L -1,
\end{equation*}
and for simplicity we assume $\Delta t$  is constant so that $t_l = l\Delta t$.  Our (temporal) training domain is therefore given by
\begin{equation}
\label{eq:training_domain}
\mathcal{D}_{train}  = [0, t_L] = [0,L\Delta t].
\end{equation}
We note that in this investigation we are interested in model prediction after time ${t_L}$.  Our framework may also be used for learning some previous behavior, for instance at time $(n+ \eta)\Delta t$, $n < {L}$ and $\eta \in (0,1)$.  This would be especially useful for cases when the time difference $\Delta t$ is very large.

\subsection{Flow map-based dynamics learning}
\label{sec:flow}

We now briefly review flow map-based deep learning of system dynamics for ordinary differential equations (ODEs) first proposed in \cite{qin2019data}, which will serve as a starting point for our flux learning technique.
To this end, we consider the dynamical system given by
\begin{equation}
\label{eq:ODEmodel}
\frac{d\mathbf{u}}{dt} = g({\bf u}), \quad\quad {\bf u}\in\R^N,
\end{equation}
where ${\bf u} = u(x_j,{t})$ and $g({\bf u}) = -\frac{\partial f}{\partial x}(u(x,{t}))\vert_{x = x_j}$, $j = 0,\dots,N$. The fundamental distinction between the problem formulation in \cite{qin2019data} and the problem in this investigation is that here we are considering a conservation law PDE model instead of a nonlinear system of ODEs.  In either case, 
although we can observe snapshots of the solution ${\mathbf u}$, the flux function $f$ in \eqref{setup:eqn} and correspondingly  $g$ in  \eqref{eq:ODEmodel}  are  unknown.  

The flow map of \eqref{eq:ODEmodel}, $\PPh:\R^{N}\times\R\to\R^{N}$, characterizes its dynamics by mapping the state variable at current time $t=0$ into a future state after some time $\Delta t$ so that  $\mathbf{u}(\Delta t) = \PPh(\mathbf{u}(0),\Delta t)$. The main idea in \cite{qin2019data} is to use a deep neural network to approximate the unknown flow map $\PPhi$ from the observation data.

The flow map-based dynamics deep learning approach begins by regrouping the  observed data in \eqref{setup:data} into pairs of adjacent time instances,
\begin{equation} \label{review:rearange_data}
	\{ \u^{(m)}(0),\u^{(m)}(\Delta t)  \}, \quad m = 1,\dots,M,
\end{equation}
where $M$ is the total number of such data pairs. Then a standard fully connected feed-forward deep neural network whose input and output layers both have  $N$ neurons is constructed. Specifically, we let $\mathcal{N}: \R^N \to\R^N$ be the associated mapping operator and {define the residual network (ResNet) mapping  as (see \cite{he2016deep})}
\begin{equation}\label{eq:resnet}
	\y^{out} =  \left[\mathcal{I} + \mathcal{N}\right]\left(\y^{in}\right),
\end{equation}
where $\mathcal{I}: \R^N \to\R^N$ is the identity operator. Using the data set in \eqref{setup:data} and {setting} $\y^{in}\leftarrow \u^{(m)}(0)$ and $\y^{out}\leftarrow \u^{(m)}(\Delta t)$, we obtain a network model
\begin{equation}
\label{eq:uNN1}
	\u^{(m)}_{NN}(\Delta t;{\Theta}) = \u^{(m)}(0) + \mathcal{N}(\u^{(m)}(0);{\Theta}).
\end{equation}
The network operator $\mathcal{N}$ can be trained by minimizing the  mean square loss function,
\begin{equation} \label{review:loss}
	\mathcal{L}(\Theta) = \frac{1}{M} \sum_{m=1}^{M} \|\mathbf{u}^{(m)}_{NN}(\Delta t;{\Theta}) - \mathbf{u}^{(m)}(\Delta t)   \|^2,
\end{equation}
where $\Theta$ denotes the network parameter set.   Once the network is satisfactorily trained we can obtain a predictive model for any arbitrarily given initial condition $\u(t_0)$.
Network structure variations and network theoretical properties can be found in \cite{qin2019data}. In \cite{qin2021deep}, the flow-map based learning was extended to include variable time stepping as well as other system parameters. Systems with missing variables were discussed in \cite{fu2022modeling}.

%
%
%

\subsection{Flow-maps for PDE models}
\label{sec:zhenmethod}
To provide more general conext for our new method, we briefly describe how the flow-map idea may be extended to learn PDE models using both modal and nodal frameworks, \cite{chen2022deep,wu2020data}, although this is not the approach used in our current  investigation. In particular the focus in \cite{chen2022deep}  was on general deep neural networks (DNN), resulting in the development of a network structure for modeling non-specific types of unknown PDEs. The network structure is based on a user specified numerical PDE solver and consists of a set of multiple disassembly layers and one assembly layer used to model the {hidden spatial} differential operators  in the unknown PDE.  The DNN model {used in \cite{chen2022deep}} defines the  mapping
\begin{equation}
\label{eq:mapping}
\u_{NN}(\Delta t) = \u(0) + \mathcal{A} (\mathcal{F}_1(\u(0)),\dots, \mathcal{F}_{J}(\u(0)), 
\end{equation}
where $\mathcal{F}_{1},\dots,\mathcal{F}_J$ are the NN operators for the disassembly layers and $\mathcal{A}$ is the NN operator for the assembly layer. The mapping in \eqref{eq:mapping} can be viewed as an application of ResNet \eqref{eq:resnet} with $\mathcal{N} = \mathcal{A} \circ (\mathcal{F}_1,\dots,\mathcal{F}_J)$. 
The same training process and loss function \eqref{review:loss} for ResNet \eqref{eq:resnet} can then also be applied to \eqref{eq:mapping}. After satisfactory training, given any new initial condition, predictions can be made by iteratively applying \eqref{eq:mapping}. 

\begin{remark}
\label{rem:zhenmethod}
It is important to note that while the method in \cite{chen2022deep} also seeks to learn PDE dynamics from data, it is not guaranteed to capture the correct shock propagation speed.  A primary motivation for the method given there is to be able to consider environments where data are collected on structure-free grids.   The assembly and disassmebly layers in \eqref{eq:mapping} require a considerable amount of hand-tuning.  This investigation, by contrast, is interested in learning hyperbolic conservation laws by incorporating the form \eqref{setup:law} directly into the neural network.  We also assume a structured grid of data, allowing the use of standard  neural network structures (fully connected feed-forward networks).
\end{remark}

\section{Constructing the network}
\label{sec:method}

Given trajectory data in \eqref{setup:data}, we now seek to construct a neural network $\mathcal{N}$ that learns the evolution of the underlying system. More precisely, we want $\mathcal{N}$ to learn to predict the state value $\mathbf{u}(t_{l+1})$ from the current state value $\mathbf{u}(t_l)$.  {In Section \ref{sec:conservation} we describe our approach for designing the network $\mathcal{N}$ for a system that is known to be conservative but for which the flux is unknown.  The more traditional approach for constructing a network without considering its conservation properties is first reviewed in Section \ref{sec:nonconservativenetwork}.}

\subsection{Standard {(}non-conservative{)} form network {(nCFN)}}
\label{sec:nonconservativenetwork}
A standard approach is to use a deep neural network $\mathcal{G}$ to approximate $f(u(x,t))_x$ in \eqref{setup:eqn} directly {for each cell average ${\bar u}_j(t^n)$} as
\begin{equation}
\label{eq:nonconservative}
	 \mathcal{G}({u}^{n}_{j-p},\dots,{u}^n_{j},\dots,{u}^n_{j+q}) \approx f({\bar u}_j(t^n))_x
\end{equation}
where {$u_j^n = \bar u_j(t_n)$}, $p\geq0, q\geq0$, and then solve the ODE given by
\begin{equation}\label{method:non-conserv}
	\frac{d}{dt}u_j +  \mathcal{G}({u}^{n}_{j-p},\dots,{u}_{j},\dots,{u}_{j+q}) = 0,
\end{equation}
with some pre-determined time integration technique.  {Importantly, \eqref{eq:nonconservative}, which we will refer to as the non-conservative form network (nCFN),  {\em does not} account for conservation  in its design.}   As will be demonstrated in Section \ref{sec:examples}, the nCFN is not able to capture the dynamics of the solution $u(x,t)$ for $t$ that extends beyond the training domain, given by $\mathcal{D}_{train}$ in \eqref{eq:training_domain}.  In particular using the nCFN may result in a numerical solution that does not satisfy the entropy condition for the weak formulation of the conservation law and yields the wrong shock speed.

A typical approach to embed the conservation property into the network model is to add a regularization term to the loss function \cite{cranmer2020lagrangian,greydanus2019hamiltonian,raissi2019physics}. Denoting the magnitude of the conserved quantity remainder in the system as $C(\mathbf{u})$, which we will derive in discrete form in Section \ref{sec:reg_conserv}, the regularization term can be constructed as
\begin{equation}\label{method:regularization}
	\mathcal{R}(\Theta) = \sum_{l=1}^L C(\mathbf{u}_{NN}(t_l;\Theta))^2.
\end{equation}
Here $\Theta$ denotes the network parameter set. In this way, the regularization term $\mathcal{R}(\Theta)$ penalizes the remainder of each conserved quantity in the network prediction. We will refer to the approach of regularizing the nCFN with the loss function \eqref{method:regularization} as nCFN-reg.  Our numerical examples will demonstrate that the solutions resulting from a non-conservative neural network are still unsatisfactory even when incorporating the regularization term, especially in long-term prediction.

\subsection{Conservative form network {(CFN)}}
\label{sec:conservation}
Motivated by classical results in numerical conservation laws, we propose a flux form network that seeks to preserve the conservation property.
Specifically, we seek to update the  cell average $\bar{u}_j(t_n)$ using the flux differences at the cell edges as
\begin{equation}\label{conserv_law}
	\frac{d}{dt}\bar{u}_j + \frac{1}{\Delta x} \Big(f_{j+1/2} - f_{j-1/2} \Big) = 0,
\end{equation}
where $f_{j+1/2}$ denotes the flux at the cell edge $x=x_{j+1/2}$. 
To approximate $f_{j+1/2}$ we define the {\em neural flux}  as
\begin{equation}\label{eq:neural_flux}
	f^{NN}_{j+1/2} = \mathcal{F}(\bar{u}^{n}_{j-p},\dots,\bar{u}_{j},\dots,\bar{u}_{j+q}),
\end{equation}
where $\mathcal{F}$ is a fully connected feed-forward neural network operator and the inputs $\bar{u}^{n}_{j-p},\dots,\bar{u}_{j},\dots,\bar{u}_{j+q}$ are neighboring cell averages centered at $x=x_{j-p},...,x_{j+q}$, respectively. For ease of presentation we denote the right hand side of \eqref{eq:neural_flux} as $\mathcal{F}_{p,q}(\bar{u}_{j})$.  Our implementation is also simplified  by using a symmetric stencil around $x_{j-1/2}$, so that $p = q-1$, although this is not required. To distinguish between the non-conservative flux forms, nCFN and nCFN-reg, we will refer to the conservative flux form network as CFN when discussing our numerical experiments.

\subsection{Time integration} \label{method:time_integration}
\begin{figure}[h!]
	\begin{center}
		\includegraphics[width=0.6\textwidth]{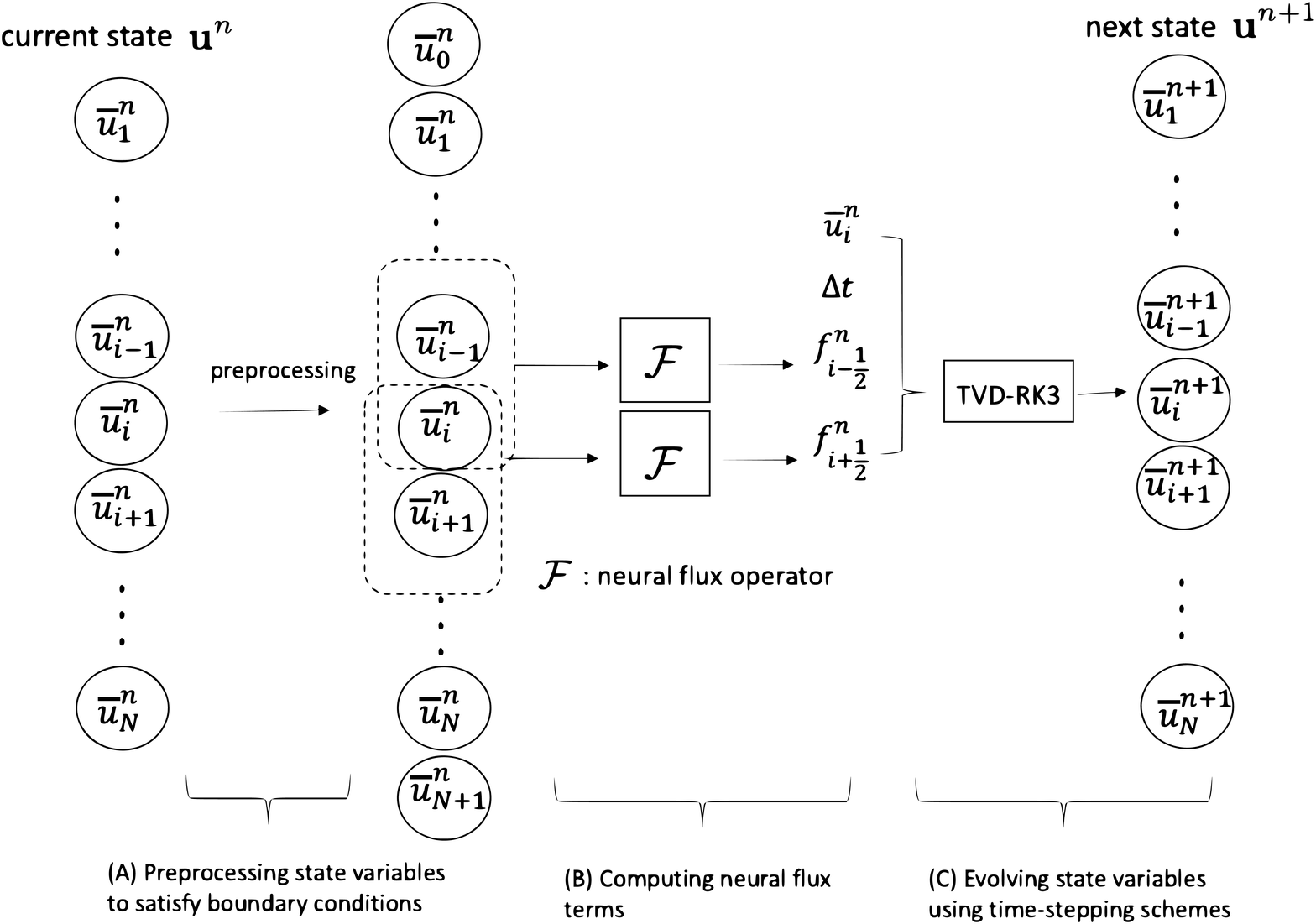}\quad\quad
		\caption{Evolution of the state variable $\mathbf{u}$ for one time step through neural network model: State variables are fed into a preprocessing layer to incorporate boundary conditions. The model then computes neural flux terms for each cell edge. The flux terms and state variables are fed into the time-integration method.}
		\label{method:diagram_evolution}
	\end{center}
\end{figure}

With the neural flux  \eqref{eq:neural_flux} in hand we now write the neural net form of \eqref{conserv_law} as
\begin{equation}\label{method:conserv}
	\frac{d}{dt}\bar{u}_j  + \frac{1}{\Delta x} \Big(\mathcal{F}_{p,q}(\bar{u}_{j}) - \mathcal{F}_{p,q}(\bar{u}_{j-1}) \Big) = 0.
\end{equation}
We solve \eqref{method:conserv} {as well as the non-conservative system in \eqref{method:non-conserv}} using the method-of-lines approach.  In our examples we will use the total variation diminishing Runge-Kutta (TVD-RK3) method, \cite{gottlieb1998total}, which for the generic system of ODEs in \eqref{eq:ODEmodel} is given by
  \begin{align}
                {\bf u}^{(1)}  &= {\bf u}^n + \Delta t \mathcal{F}({\bf u}^{n}),\nonumber\\
                {\bf u}^{(2)}  &= \frac{3}{4} {\bf u}^{n} + \frac{1}{4} {\bf  u}^{(1)} + \frac{1}{4} \Delta t {\mathcal F}({\bf u}^{(1)}), \nonumber\\
                {\bf u}^{n+1} &= \frac{1}{3} {\bf u}^n + \frac{2}{3}{\bf  u}^{(2)} + \frac{2}{3} \Delta t {\mathcal F}({\bf u}^{(2)}),
\label{eq:tvd-rk3}
        \end{align}
for integration between time steps $n$ and $n+1$.
Figure \ref{method:diagram_evolution} provides a diagram showing the evolution of the state variable ${\bf u}$ for one time step through the conservative neural network model.

\begin{remark}
\label{rem:TVD-RK3}
To the best of our knowledge there are no theoretical results regarding the choice of time integration schemes that guarantee stability for numerical PDEs using neural networks.  Thus, due to its theoretical stability guarantees for numerical conservation laws when using traditional solvers, we choose the TVD-RK3 method. Other time integration techniques may also be appropriate, and in some cases reduce computational cost or improve numerical accuracy.
\end{remark}


\subsection{Boundary Conditions}
For simplicity we will assume that the boundary conditions in \eqref{setup:eqn} are known. 
In particular, to satisfy the periodic boundary conditions in Example \ref{ex:burgers} we simply apply 
\begin{equation}
\bar{u}_{0} = \bar{u}_N, \quad
\bar{u}_{N+1} =  \bar{u}_1.
\label{eq:BCperiod}
\end{equation}
No flux boundary conditions are assumed in Examples \ref{ex:SWE} and \ref{ex:euler}. Since the solution profiles near the boundaries remain constant over time, we simply impose the boundary conditions for each variable in both examples as
\begin{equation}
\bar{u}_{0} = \bar{u}_1, \quad
\bar{u}_{N+1} =  \bar{u}_N.
\label{eq:BCzeroth}
\end{equation}
Higher order numerical boundary conditions can similarly be employed.



\subsection{The Recurrent Loss Function} \label{method:RNN} 
Given trajectory data in \eqref{setup:data}, where each trajectory has multiple measurements, we define the {\em recurrent loss function} as
\begin{equation} \label{method:loss}
	\mathcal{L}_{RNN}(\Theta) = \frac{1}{N_{traj}} \frac{1}{L}\sum_{k=1}^{N_{traj}} \sum_{l=1}^{L} \|\mathbf{u}^{(k)}_{NN}(t_l;{\Theta}) - \mathbf{u}^{(k)}(t_l)   \|^2,
\end{equation}
where 
\begin{equation}
	\label{eq:uNN}
	\mathbf{u}_{NN}^{(k)}(t_l;\Theta) = \underbrace{\mathcal{N} \circ \cdots \circ \mathcal{N}}_\text{$l$ times} (\mathbf{u}^{(k)}(t_0)).
\end{equation}

As already discussed, the network evolution operator $\mathcal{N}$  is designed to  predict the state value $\mathbf{u}(t_{l+1})$  from the current state value $\mathbf{u}(t_l)$, 
where  $\Theta$ again denotes network parameter set. In contrast to  \eqref{review:loss}, the recurrent loss function \eqref{method:loss} calculates loss over multiple time steps.  Using the recurrent loss approach has been found to improve  numerical stability \cite{chen2022deep}.  
As mentioned in Section \ref{sec:nonconservativenetwork}, one can embed the conservation property into the network by adding a regularization term in the loss function,
\begin{equation} \label{method:loss-reg}
		\mathcal{L}^{\lambda}_{RNN}(\Theta) = \mathcal{L}_{RNN}(\Theta)	+ \lambda^2 \mathcal{R}(\Theta),
\end{equation}
where $\mathcal{R}$ is defined in \eqref{method:regularization}. Following ML conventional notation we use $\lambda^2$ to denote the weighting parameter for the regularization term. 

\section{{Experiment design}}
\label{sec:numexamples}

Our results presented in Section \ref{sec:examples} will demonstrate that for the classical one-dimensional conservation laws studied in this investigation, using the CFN, for which the network learns the flux function of the unknown conservation law via \eqref{eq:neural_flux} from observation data, yields significantly better results than those obtained using either the nCFN in \eqref{eq:nonconservative} or the nCFN regularized by the loss function, nCFN-reg.  Below we provide the framework necessary to ensure the robustness and reliability of  our experimental results. To this end we consider details related both to data collection and training.

\subsection{Data Collection} \label{sec:data_setup}
To test our method we will consider examples of conservation laws for which the fluxes are known.\footnote{{Indeed in some of our examples the true solution is also known. However, since we randomly generate the initial conditions to obtain a set of $N_{traj}$ snapshots, we will simply consider the ``exact'' (reference)  solution to be the highly resolved numerical result.}}  We will use this information both to generate synthetic training data with which to train the DNN for the evolution process as well as to compute reference solutions to evaluate our results. Importantly we note that knowledge of the true system does not in any way facilitate the DNN model approximation. 

Our test problems range from idealistic, where we assume we have noise-free densely observed data for training and validation, to more difficult situations, where we consider  two cases: (1) the observable data are accurate (noiseless) but sparsely observed  and (2) the observable data are noisy but densely observed.
To mimic observed data in \eqref{setup:data} that would be available for training and validation, we numerically simulate the true underlying PDE model according to the observational settings provided in the particular case study for our examples.\footnote{Unless otherwise noted we use the CLAWPACK conservation laws package, \cite{clawpack}.}  We furthermore randomly sample the parameters in the initial conditions to obtain various trajectories of the observed data (see e.g. \eqref{eq:invis_IC1}).
The number of trajectories $N_{traj}$ and the trajectory length $L$ vary depending on the underlying properties of the PDE {(e.g.~time to shock formation)}.

\subsection{Network and Training Details}
As shown in  \figref{method:diagram_evolution}, the CFN consists of the  preprocessing layer, which ensures that the boundary conditions are satisfied,  and the neural flux operator, which computes the flux at cell edges. The neural flux operator is constructed using a fully connected feedforward neural network and is obtained by training the network hyperparameters (weights and biases) as the minimum of the recurrent loss function \eqref{method:loss}.   We employ the stochastic optimization method {\bf Adam} \cite{kingma2014adam} for this purpose.  The nCFN and nCFN-reg utilize the same preprocessing layer for the boundary conditions and each employs one fully connected feedforward network to learn the increment of the state variables.   For consistency all models are trained for 10,000 epochs with learning rate $10^{-4}$ in every example.  The same set of network structures is also employed.   Finally, we use the commonly chosen Rectified Linear Unit (ReLU) \cite{lecun2015deep} as the activation function.    These network and training details are summarized in Table \ref{table:NN_structure}.
\begin{table}[h!]
	\centering
	\begin{tabular}{ | m{1.9cm} | m{1.9cm} | m{2cm} |  m{2cm} | m{1.9cm} |}
		\hline
		model & (p,q) & hidden layers & hidden nodes & activations  \\
		\hline
		CFN & (2,3) & 5  & 64 & ReLU\\
		\hline
		nCFN & (3,3) & 5  & 64 & ReLU\\
		\hline
		nCFN-reg & (3,3) & 5  & 64 & ReLU\\
		\hline
	\end{tabular}
	\caption{Neural network architecture details for all examples. Note that for each model  $p$ and $q$ are chosen to provide symmetry ($p=q-1$ for CFN and $p=q$ for nCFN), although this is not a requirement.}
	\label{table:NN_structure}
\end{table}


It is also possible to tune the regularization parameter $\lambda^2$ for the nCFN-reg loss function in \eqref{method:loss-reg}.  Indeed, one can choose $\lambda^2 = \lambda^2(t)$, so that the influence of the regularization can fluctuate as the PDE evolves.  This would add considerable computational cost, however, and moreover, it is not readily apparent that employing standard approaches, such as the $L-$curve method or the discrepancy principle, \cite{Hansenbook},  are appropriate here.  Hence in our examples we choose the $\lambda^2$ from the values $10^{2(i-1)}$, $i = 1,\dots,4$, that yields the smallest residual solution after the first time step, and it remains fixed throughout the rest of the process.  In general we found that in Examples \ref{ex:burgers} and \ref{ex:SWE} that $\lambda^2 = 10^{-2}$ yielded the best results.  Example \ref{ex:euler}  (the Euler equations for gas dynamics) was considerably more sensitive to the choice of $\lambda^2$, likely due to the oscillatory nature of the solution.  In this case we refined our search to include $\lambda^2 = 5\times 10^{-2}$.   We therefore see that as an added advantage our new CFN approach  does not require regularization parameter tuning.

We emphasize that while our numerical experiments indicate that these parameter choices provide enough network complexity for each required learning task, we did not further try to optimize performance.
Moreover, as we want to ensure the robustness of our method, in our experiments we typically follow the common practice for learning system dynamics \cite{qin2019data,qin2021data,qin2021deep,wu2020data} and use the default values in Tensorflow or {other standard} choices for all hyperparameters. 

\subsection{Constructing the Regularization Term} \label{sec:reg_conserv}

The regularization term \eqref{method:regularization} is designed to promote conservation in the nCFN-reg method.  Below we show how this term is constructed for the scalar case. A straightforward extension can be made for systems.

We first expand \eqref{setup:law} to the physical domain of the problem, $(a,b)$, yielding
\begin{equation}\label{eq:law_Qoc}
	\int_{a}^{b} u(x,t_{l+1}) dx - \int_{a}^{b} u(x,t_{l}) dx = \int_{t_l}^{t_{l+1}}f(u(a,t))dt - \int_{t_l}^{t_{l+1}} f(u(b,t))dt, l = 0,\dots,L,
\end{equation}
where each $t_l$ denotes the time at which a data trajectory in \eqref{setup:data} is initially obtained.

Example \ref{ex:burgers} considers the inviscid Burgers equation with periodic boundary conditions.  In this case \eqref{eq:law_Qoc} simplifies to
\begin{equation}
\label{eq:massburgers}
	\int_{a}^{b} u(x,t_{l}) dx = \int_{a}^{b} u(x,t_0) dx, \quad l=0,\dots,L.
\end{equation}
For equations with non-periodic boundary conditions, \eqref{eq:massburgers} does not hold since in general $f(u(a,t)) \ne f(u(b,t))$.  Hence to construct \eqref{method:regularization}  we first define
\begin{align}
\label{eq:fluxterm}
	F_a^{(l)} = \frac{1}{\Delta t}\int_{t_{l}}^{t_{l+1}} f(u(a,t))dt,\quad
	F_b^{(l)} = \frac{1}{\Delta t}\int_{t_{l}}^{t_{l+1}} f(u(b,t))dt,\quad
l = 0,\dots,L-1,
\end{align}
and then use \eqref{eq:fluxterm} to approximate \eqref{eq:law_Qoc} as
\begin{equation}
	\label{eq:discrete_conserv_law}
	\sum_{j=1}^{N} \bar{u}_j(t_{l+1})\Delta x - \sum_{j=1}^{N} \bar{u}_j(t_{l})\Delta x = F_a^{(l)}{\Delta t} - F_b^{(l)}{\Delta t}, \quad l = 0,\dots L-1,
\end{equation}
which leads to
\begin{equation}\label{eq:discrete_conserv}
	\sum_{j=1}^{N} \left(\bar{u}_j(t_l) - \bar{u}_j(t_{0})\right)\Delta x = \sum_{k=1}^{l} \left(F_a^{(k-1)} - F_b^{(k-1)}\right){\Delta t},\quad l = 1,\dots,L.
\end{equation}
From here we define the (discrete) conserved quantity {remainder} at each $t_l$ as
\begin{equation}
\label{eq:conserve_u}
	C(\mathbf{u}(t_l)): = \left|\sum_{j=1}^{N} \left(\bar{u}_j(t_l) - \bar{u}_j(t_{0})\right)\Delta x - \sum_{k=1}^{l} \left(F_a^{(k-1)} - F_b^{(k-1)}\right){\Delta t}\right|,
\end{equation}
where $\mathbf{u}(t) = (\bar{u}_1(t), \dots \bar{u}_N(t))^{T}$.   It follows from \eqref{eq:discrete_conserv} that if the conservation property holds then $C(\mathbf{u}(t_l))=0$.  Regularization in \eqref{method:loss-reg} is therefore used to promote solutions that minimize \eqref{eq:conserve_u}.
In practice the network prediction of $\mathbf{u}_{NN}(t_l;\Theta)$  is used to calculate \eqref{eq:conserve_u}, directly yielding $\mathcal{R}(\Theta)$ in \eqref{method:regularization}.
We note that we will also be able to analyze the conservation properties of each of our numerical methods in Section \ref{sec:examples} by computing \eqref{eq:conserve_u} over the time domain of the solution.

\begin{remark}
It is important to point out that \eqref{eq:conserve_u} describes a best case scenario, where we have access to \eqref{eq:fluxterm}.  In order to construct the regularization term for the nCFN-reg in our experiments, we compute \eqref{eq:fluxterm} directly from the given flux terms in each example. This serves to demonstrate that even under ideal circumstances, regularizing the standard nCFN to promote conservation in the solution (nCFN-reg) is not as effective as constructing a conservative form network in the first place (CFN).
\end{remark}

\section{Numerical Examples}
\label{sec:examples}

We use three  well-studied one-dimensional examples of hyperbolic conservation laws to analyze our new conservative form network (CFN) approach.  We consider three different observational settings for each experiment: (i) an ideal case, where the observations are dense and noise-free; (ii) the situation where the observations are sparse but noise-free; and (iii) an environment for which the observations are dense but noisy. We compare the results of our new CFN approach to the more traditional non-conservative form network (nCFN) along with the regularized (nCFN-reg) version.

\subsection{Inviscid Burgers Equation}
\label{sec:burgers}

\begin{example}
\label{ex:burgers}
The inviscid Burgers equation is given by
\begin{equation}
\label{eq:burgers}
	u_t + (\frac{u^2}{2})_x = 0, \quad x \in (0,2\pi), \quad  t > 0,
\end{equation}
with periodic boundary conditions $u(0,t) = u(2\pi,t)$. We choose initial conditions
\begin{eqnarray}
u(x,0) &=& \alpha + \beta \sin (x),\nonumber\\\
\alpha &\sim& U[-\epsilon_s, \epsilon_s], \nonumber\\
\beta  &\sim& U[1-\epsilon_s, 1+\epsilon_s], 
\label{eq:burg_init}
\end{eqnarray}
where $\epsilon_s = 0.25$.
\end{example}
The $N_{traj}$ training data sets are generated by solving \eqref{eq:burgers} using the Engquist–Osher flux along with TVD-RK3 time integration based on the initial conditions
\begin{eqnarray}
	u^{(k)}(x,0) &=& \alpha ^{(k)} + \beta ^{(k)} \sin (x), \nonumber\\
\alpha^{(k)} &\sim& U[-\epsilon_s, \epsilon_s], \nonumber\\
\beta^{(k)} &\sim& U[1-\epsilon_s, 1+\epsilon_s]
\label{eq:invis_IC1}
\end{eqnarray}
for $k = 1,\dots,N_{traj}$ with $\epsilon_s = 0.25.$

In all of our experiments we set $N_{traj} = 200$. Each training trajectory has length $L = 20$  for either choice of recurrent loss function, \eqref{method:loss} {or \eqref{method:loss-reg}}.

To check the {\em robustness} of our predictions for Example \ref{ex:burgers}, we run our experiments for 50 choices of fixed $\alpha$ and $\beta$ and compare the three methods, CFN, nCFN, and nCFN-reg. Our reference solution is calculated using the Engquist-Osher flux term on a fine grid, with $\Delta x = \frac{2\pi}{1024}$ in \eqref{method:conserv}. All figures use the initial value with $\alpha = 0.06342$ and $\beta = 1.17322$ for illustration. Other choices for $\alpha$ and $\beta$ yield comparable results.




%

\subsubsection*{Case I: Dense and Noise-free Observations}

We first consider an idealized environment for which the observations are dense and noise-free. Specifically we choose $N = 512$, yielding $\Delta x = \frac{2\pi}{512}$, so that our solution is well-resolved. We also choose a constant time step $\Delta t$ for all experiments so as not to complicate our analysis.  In this regard we observe that the maximum wave speed for Burgers equation, $|u(x,t)|$, can be determined for all $t$ using \eqref{eq:invis_IC1} as $\max\{|u|\} = 1+2\epsilon_s$. We therefore set $\Delta t = 0.005$ to satisfy the CFL condition with $\#CFL = 0.9$.

We note that the training time domain $[0,L\Delta t]$ with $L=20$ contains only smooth solution snapshots. Since each DNN model requires the training data to include both smooth and discontinuous solution profiles to learn the long term dynamics of Example \ref{ex:burgers},  a larger trajectory length $L$ is needed.  As choosing a larger $L$ would significantly increase computational costs  we employ a sub-sampling technique to generate training data from observed snapshots of the solution onto  an extended domain. The same  sub-sampling technique is used for Example \ref{ex:euler}.  The details are provided below.

We define a new parameter $M > L$ as the extended length of each trajectory. The snapshots of the solution, \eqref{setup:data}, are obtained for each of the $N_{traj}$ trajectories at times  $t = m\Delta t$, $m = 1,\dots, M$. In our experiments we choose $M = 300$ which yields the total training time domain as $[0,1.5]$. We then sub-sample each of the $k = 1,\dots,N_{traj}$ by randomly selecting a start time value, $t_0^{(k)}$, from the set $\{\mu \Delta t\}_{\mu = 0}^{M-L}$. Each sub-sampled trajectory of length $L$ is then built consecutively from the snapshot solutions. That is, each trajectory is comprised of the solutions in \eqref{setup:data} at sequential times $t_0^{(k)} + l \Delta t$, $l = 1,\dots ,L$. In this way we can train over a longer period of time without increasing the expense of network training. This approach, of course, requires that more initial observations are available.

\begin{figure}[!h]
	\centering
	\begin{subfigure}[b]{0.24\textwidth}
		\includegraphics[width=\textwidth]{%
				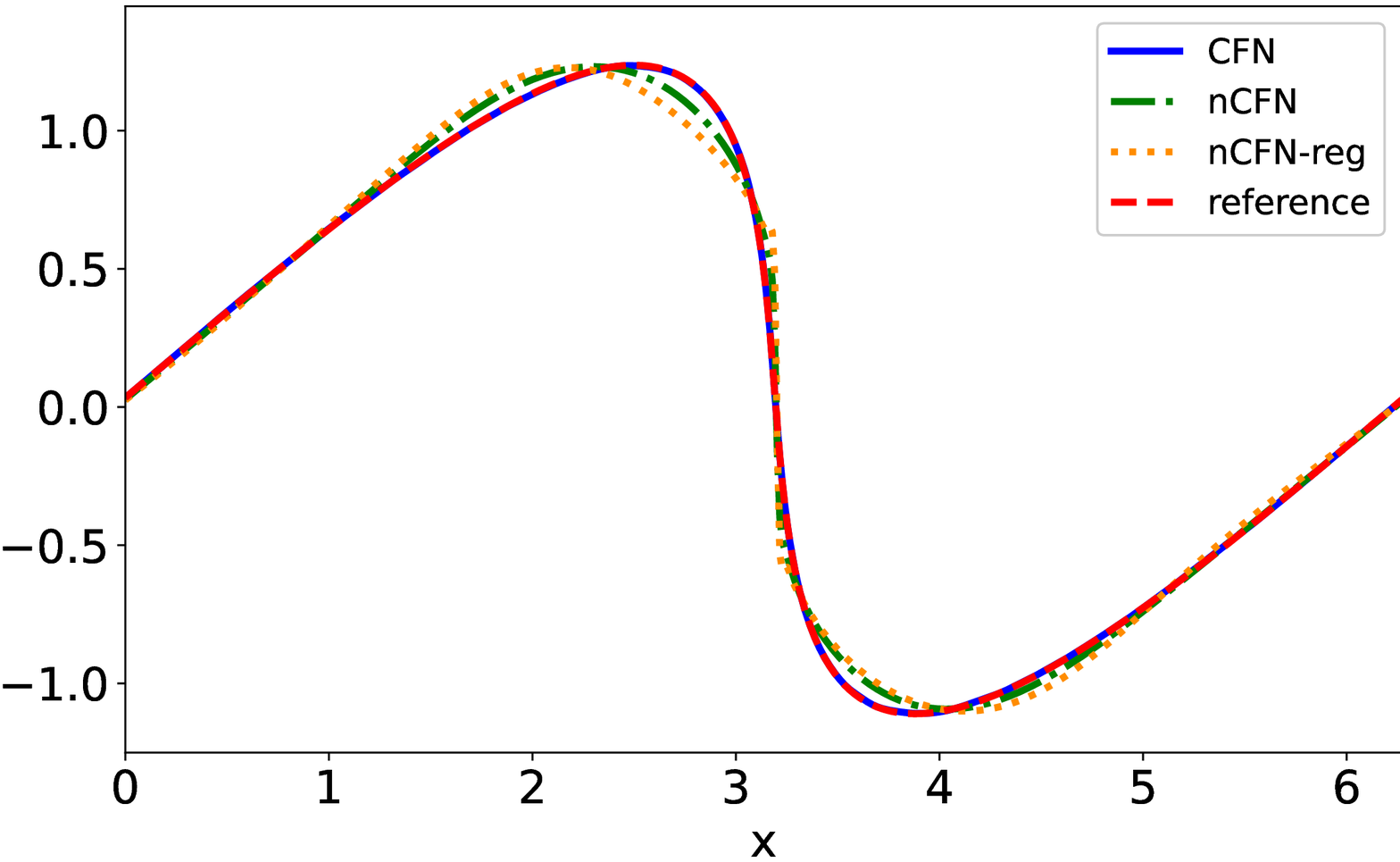}
		\caption{$t = 0.75$}
		\label{fig:burgt1case1}
	\end{subfigure}%
	~
\begin{subfigure}[b]{0.24\textwidth}
		\includegraphics[width=\textwidth]{%
			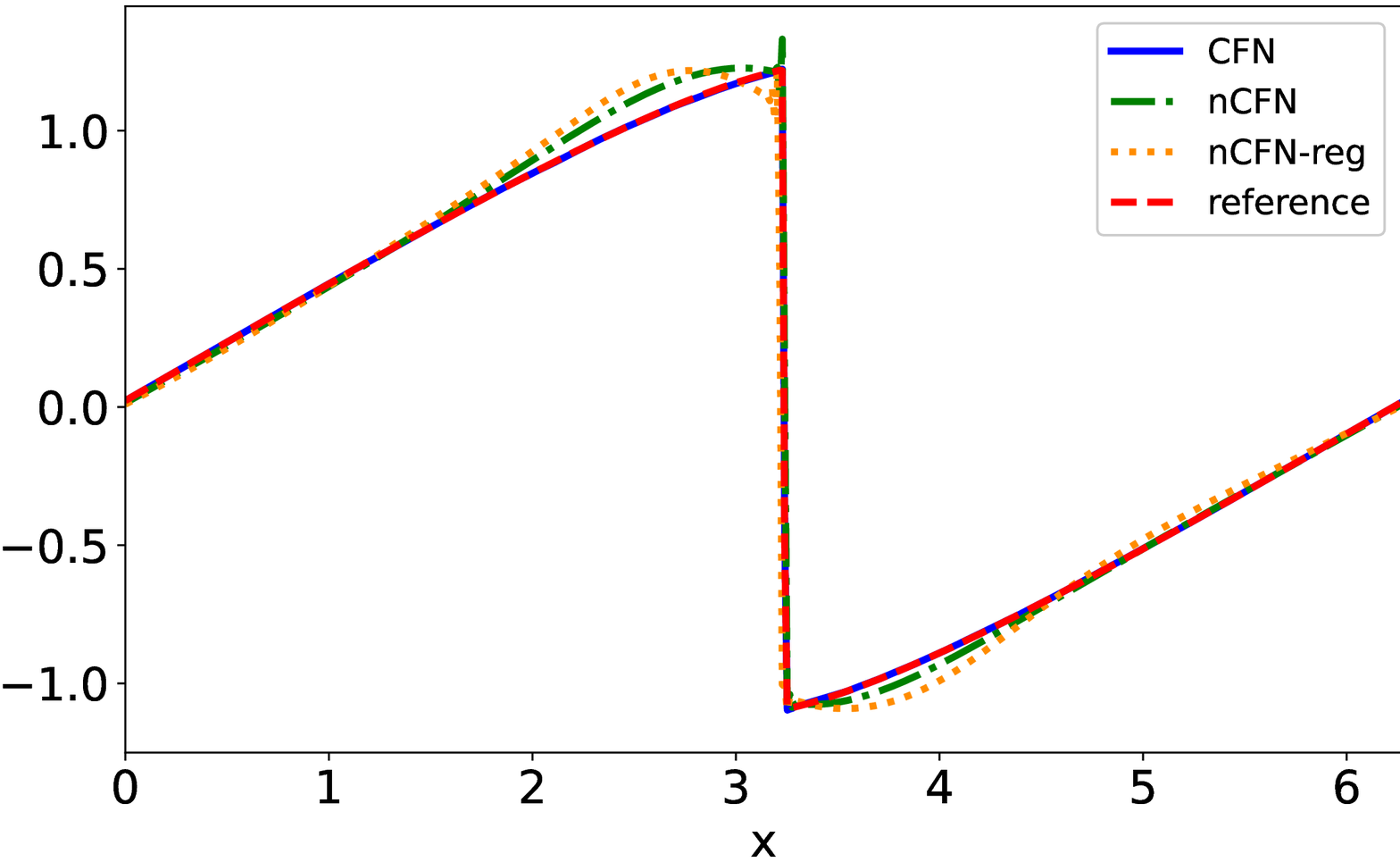}
		\caption{$t = 1.5$}
		\label{fig:burgt2case1}
	\end{subfigure}%
	~
\begin{subfigure}[b]{0.24\textwidth}
		\includegraphics[width=\textwidth]{%
			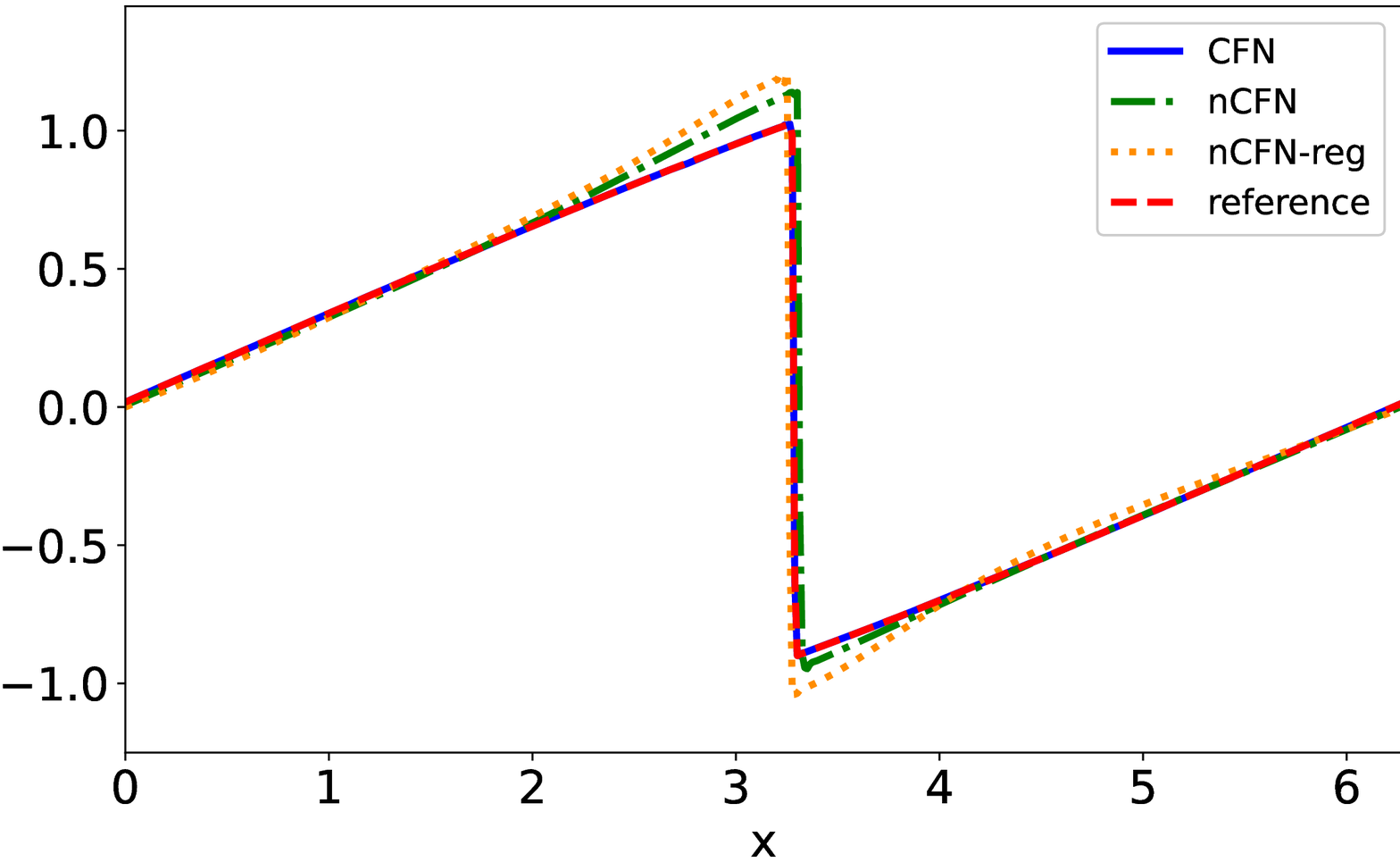}
		\caption{$t = 2.25$}
		\label{fig:burgt3case1}
	\end{subfigure}%
	~
\begin{subfigure}[b]{0.24\textwidth}
		\includegraphics[width=\textwidth]{%
			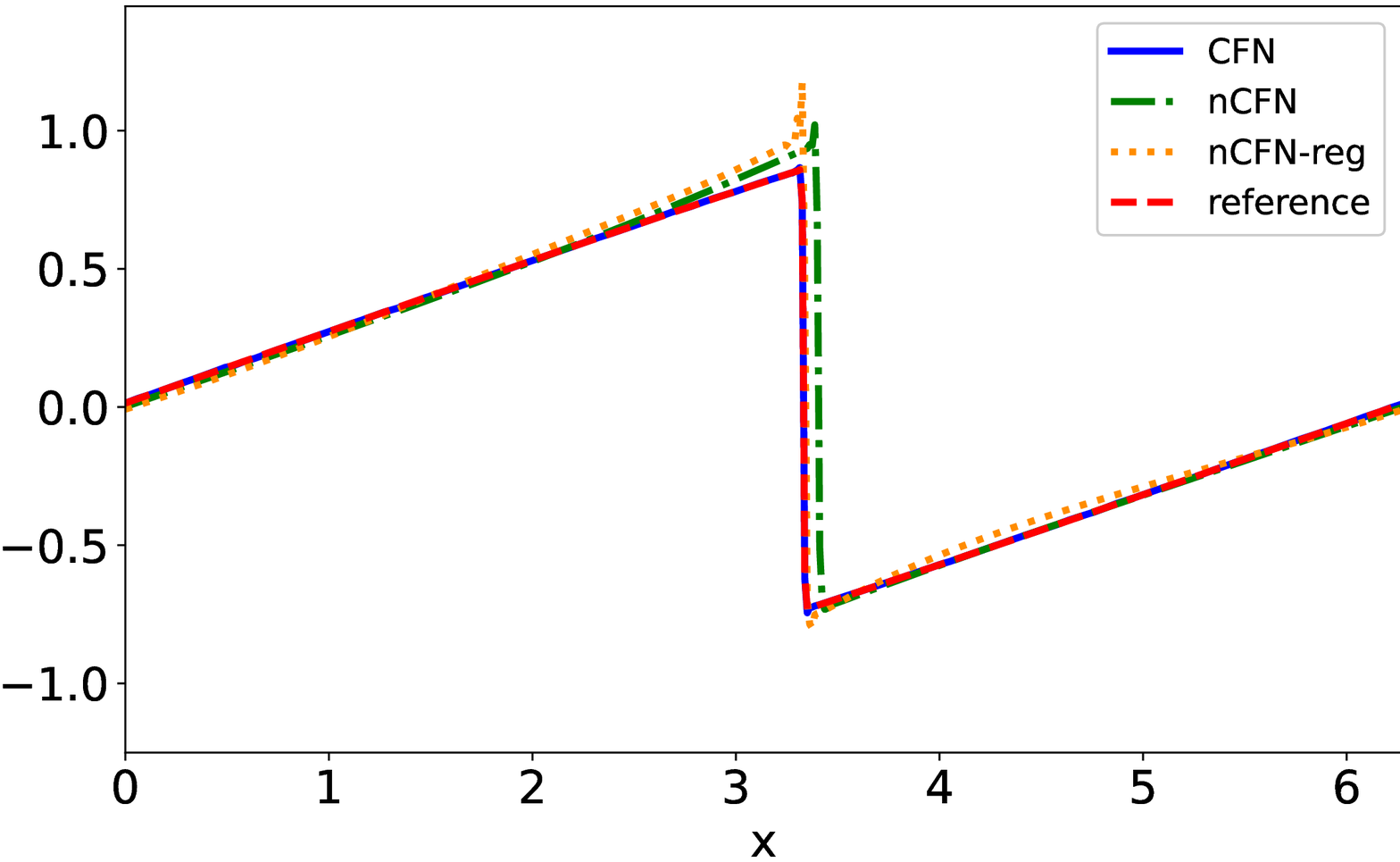}
		\caption{$t = 3$}
		\label{fig:burgt4case1}
	\end{subfigure}%
	\caption{Comparison of the reference solution to Example \ref{ex:burgers} with the trained DNN model predictions at different times for dense ($N=512$) and noise-free observations.}
	\label{fig:burg_ideal}
\end{figure}
\begin{figure}[h!]
	\centering
			\includegraphics[width=0.3\textwidth]{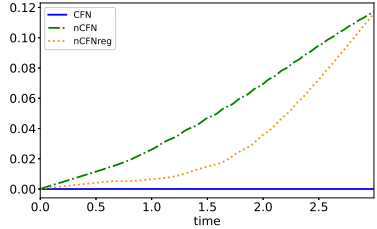}
			\caption{Discrete conserved quantity remainder $C(\mathbf{u})$ in \eqref{eq:conserve_u} of the network predictions for Example \ref{ex:burgers}.}
			\label{fig:conserv_error}
\end{figure} 
\figref{fig:burg_ideal} presents the solution to Example \ref{ex:burgers} for this ideal case at four different times, within and beyond training time domain $[0,1.5]$. Observe that the three methods capture the solution profiles and predict the correct shock propagation speed within the training domain (shown for $t=0.75$ and $t=1.5$).  The nCFN and nCFN-reg results are less accurate, and do not appear to be completely resolved. Beyond the training domain ($t>1.5$), only the CFN and nCFN-reg methods yield the correct shock propagation speed (\figref{fig:burgt3case1}, \figref{fig:burgt4case1}). The nCFN-reg solution develops a non-physical overshoot near the shock location. This behavior is further observed in \figref{fig:conserv_error}, where the conserved quantity remainder $C({\bf u})$  obtained by \eqref{eq:conserve_u} is displayed for each method. Clearly the CFN produces the only conservative method.


\begin{remark}
We also compared our results to those obtained using the method in \cite{chen2022deep} which has a global design and beyond the fully-connected layers also  contains additional disassembly and assembly layers.  This structure inherently means that the method has significantly more parameters to tune and also requires more training when compared to our CFN approach, which has a local flux structure.  In particular the set of training data provided in all of our case studies, including the idealized environment,  leads to  overfitting in the training process and fails to yield conservation.  Additional training data will lead to more comparable results, although there is no guarantee that they will ultimately satisfy the entropy conditions for the weak formulation of the PDE. The shock speed of propagation may be incorrect and the solution may exhibit non-physical oscillations near the shock.\footnote{Indeed, a primary motivation in \cite{chen2022deep} is to learn the dynamics of generic PDEs on unstructured grids, and the data in our examples are collected on structured grids.} 
\end{remark}

\subsubsection*{Case II: Sparse and noise-free observations}
\begin{figure}[!hbtp]
	\centering
\begin{subfigure}[b]{0.24\textwidth}
		\includegraphics[width=\textwidth]{%
				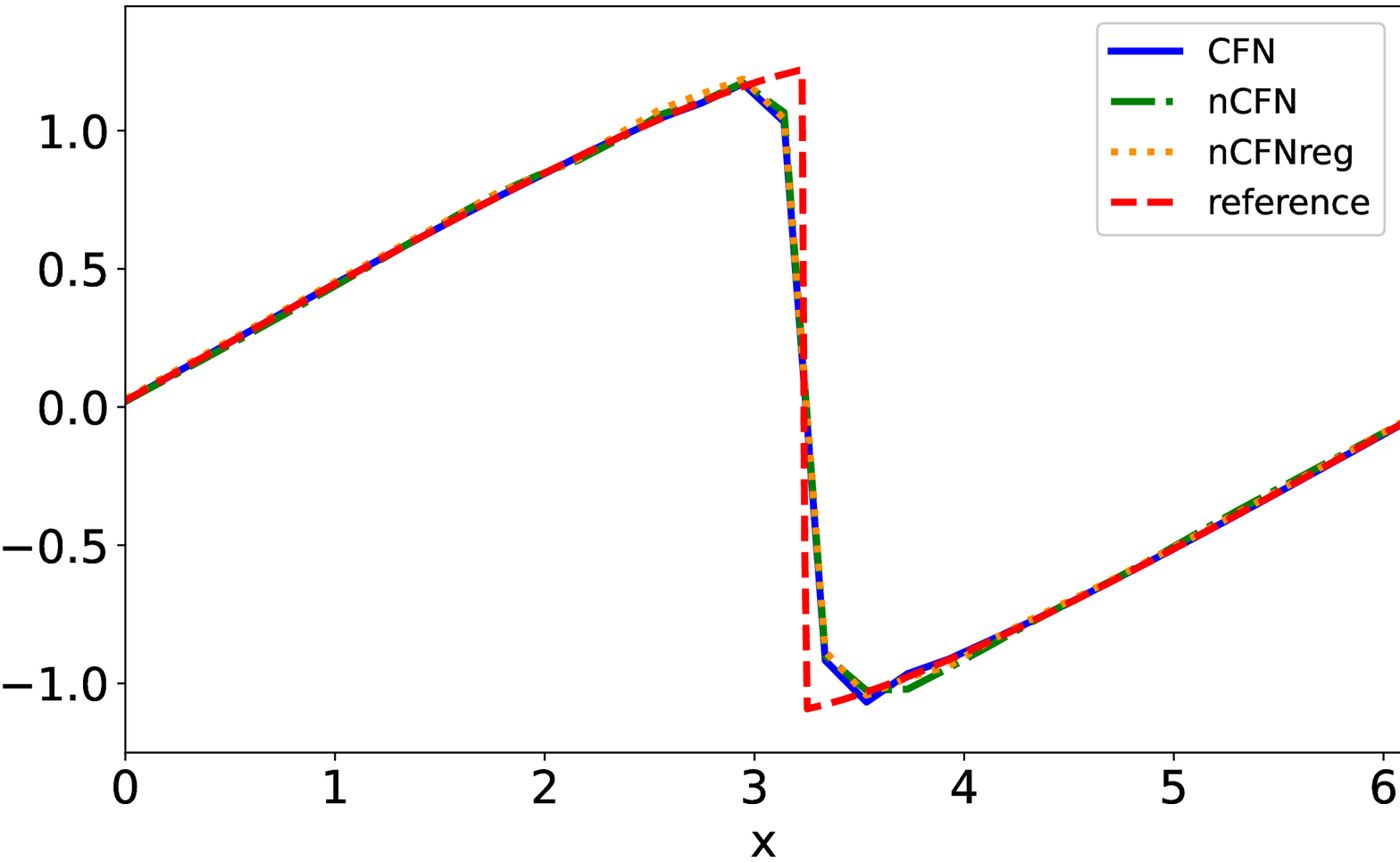}
		\caption{ $N = 32$, $t = 1.5$}
		\label{fig:burgN32_t1_case2}
	\end{subfigure}%
	~
	\begin{subfigure}[b]{0.24\textwidth}
		\includegraphics[width=\textwidth]{%
				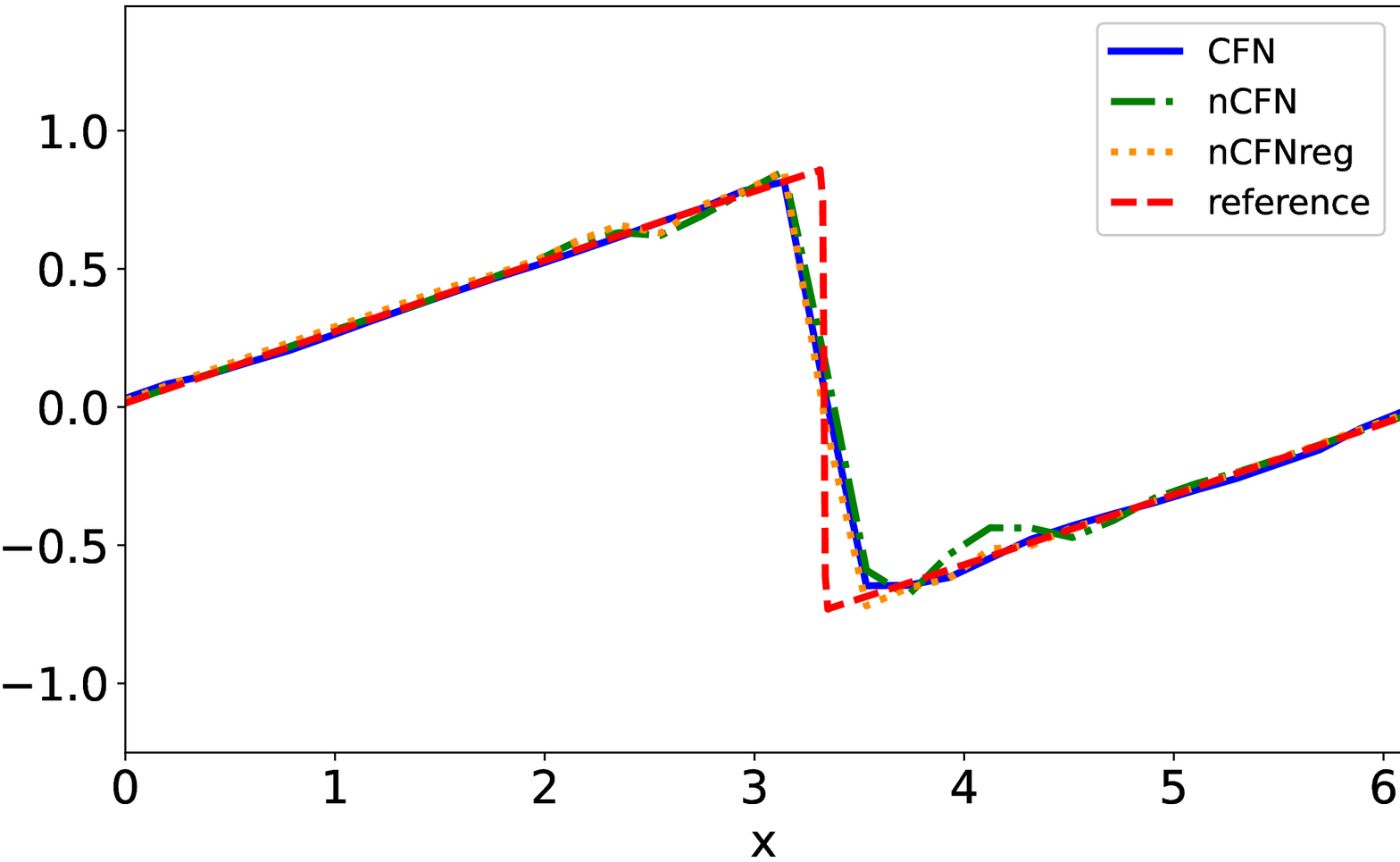}
		\caption{ $N = 32$, $t = 3$}
		\label{fig:burgN32_t2_case2}
	\end{subfigure}%
\begin{subfigure}[b]{0.24\textwidth}
		\includegraphics[width=\textwidth]{%
				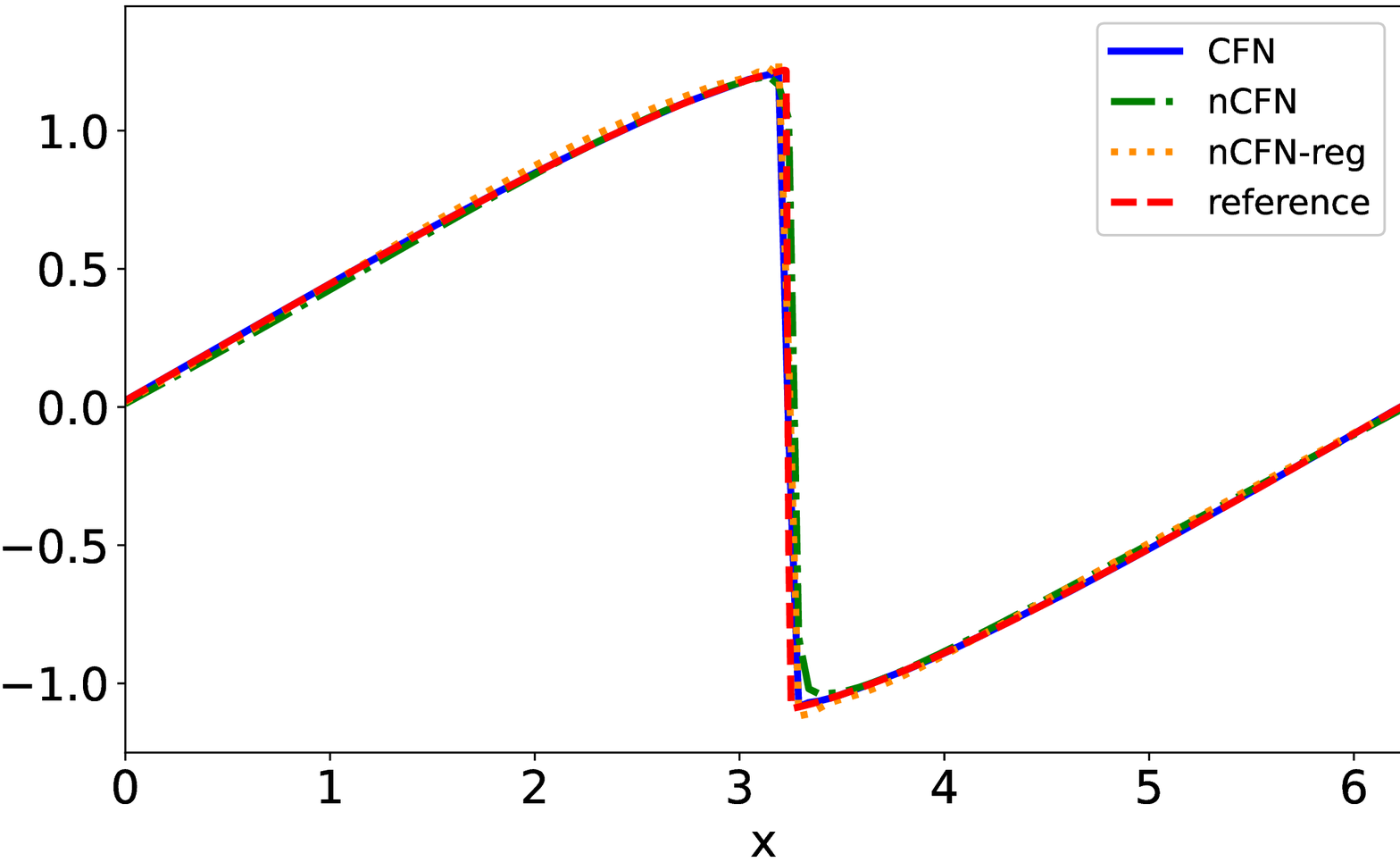}
			\caption{ $N = 128$, $t = 1.5$}
		\label{fig:burgN128_t1_case2}
	\end{subfigure}%
	~
	\begin{subfigure}[b]{0.24\textwidth}
		\includegraphics[width=\textwidth]{%
			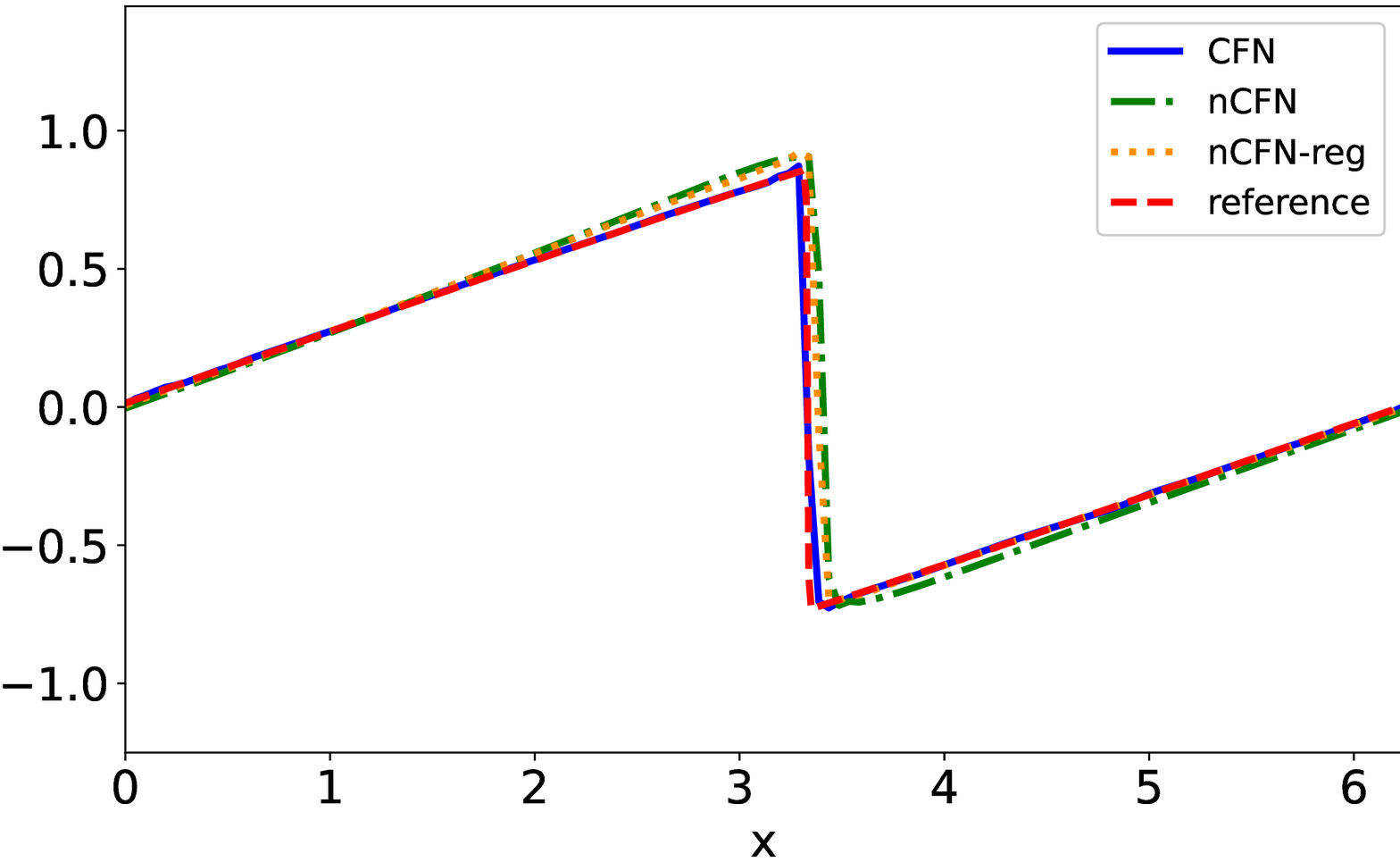}
		\caption{ $N = 128$, $t = 3$}
		\label{fig:burgN128_t2_case2}
	\end{subfigure}%
	\caption{Comparison of the reference solution to Example \ref{ex:burgers} with the trained DNN model predictions at different sparsity levels using noise-free observations.}
	\label{fig:burg_sparse}
\end{figure}
In this case the training data are obtained by solving Example \ref{ex:burgers} on a coarse grid. Specifically, for $\Delta x = \frac{2\pi}{N}$ we choose $N = 32,128$. Once again we fix the time step as $\Delta t  = 0.005$.
\figref{fig:burg_sparse} compares the results using the CFN, nCFN and nCFN-reg for different sparsity levels at times in ($t=1.5$) and out of ($t=3$) the training time domain. Both time instances are after the shock forms. Observe that for each choice of $N$ only the CFN captures the correct shock propagation speed. {\figref{fig:burg_sparse_ptErr} displays the pointwise error at different sparsity levels for each method when $t=3$.  {It is apparent that the width of the interval containing the error resulting from shock shrinks (as expected) with increased resolution for all three methods.  However, neither the nCFN nor the nCFN-reg demonstrate convergence.}

\begin{figure}[!h]
	\centering
	\begin{subfigure}[b]{0.322\textwidth}
		\includegraphics[width=\textwidth]{%
			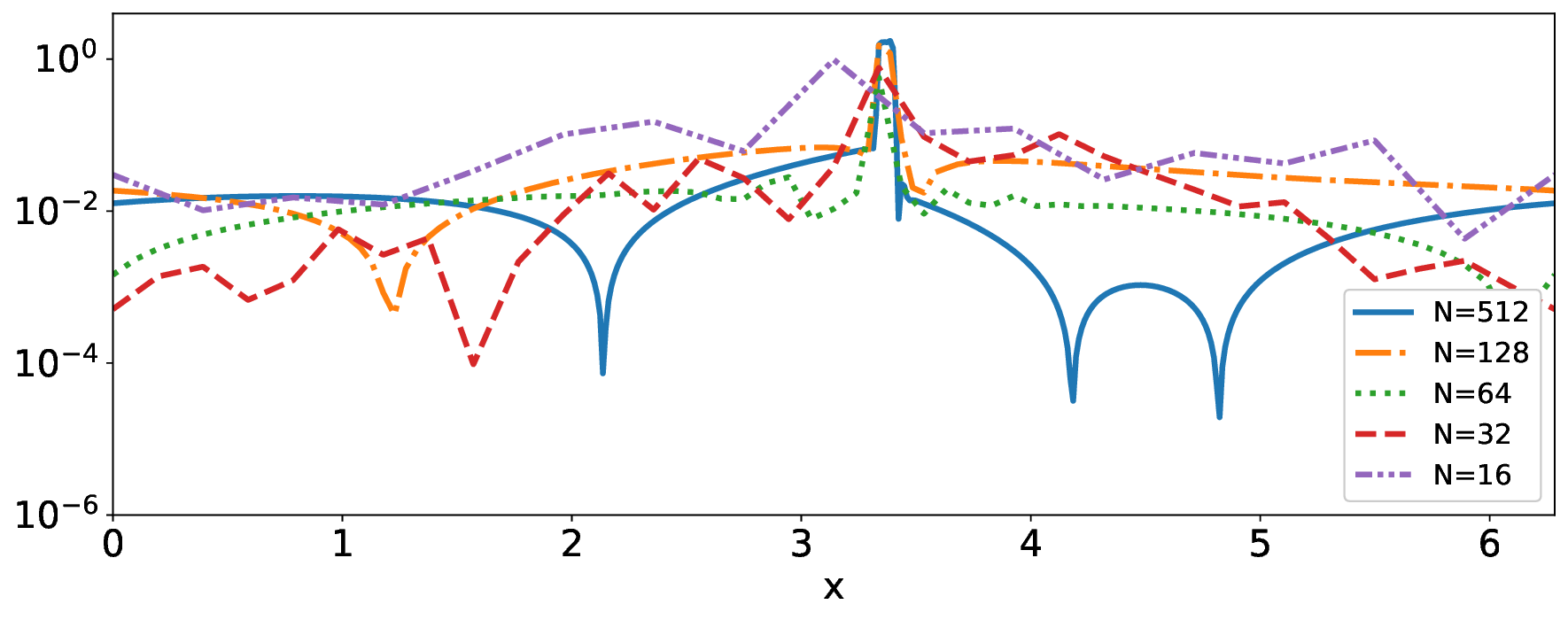}
		\caption{nCFN}
		\label{fig:ptErr_nCFN_sparse}
	\end{subfigure}%
	\begin{subfigure}[b]{0.322\textwidth}
		\includegraphics[width=\textwidth]{%
			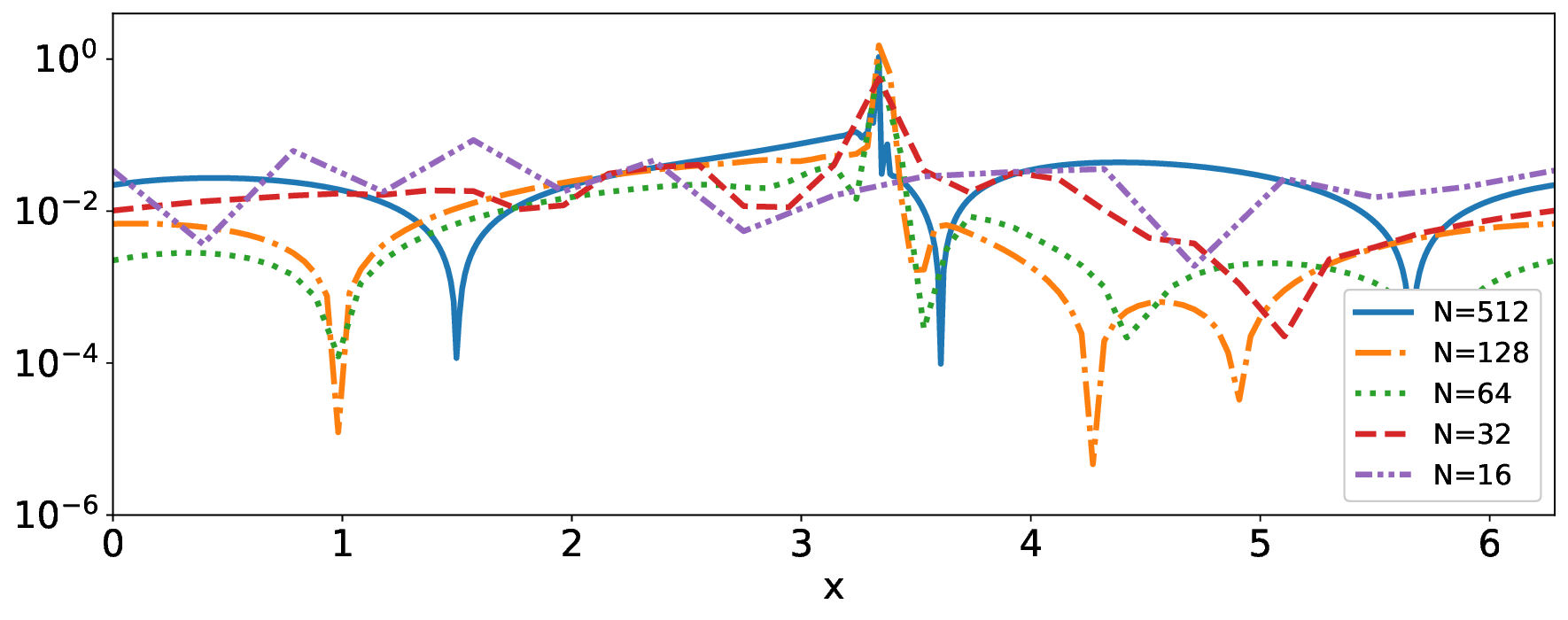}
		\caption{nCFN-reg}
		\label{fig:ptErr_nCFNreg_sparse}
	\end{subfigure}%
	~
	\begin{subfigure}[b]{0.322\textwidth}
		\includegraphics[width=\textwidth]{%
			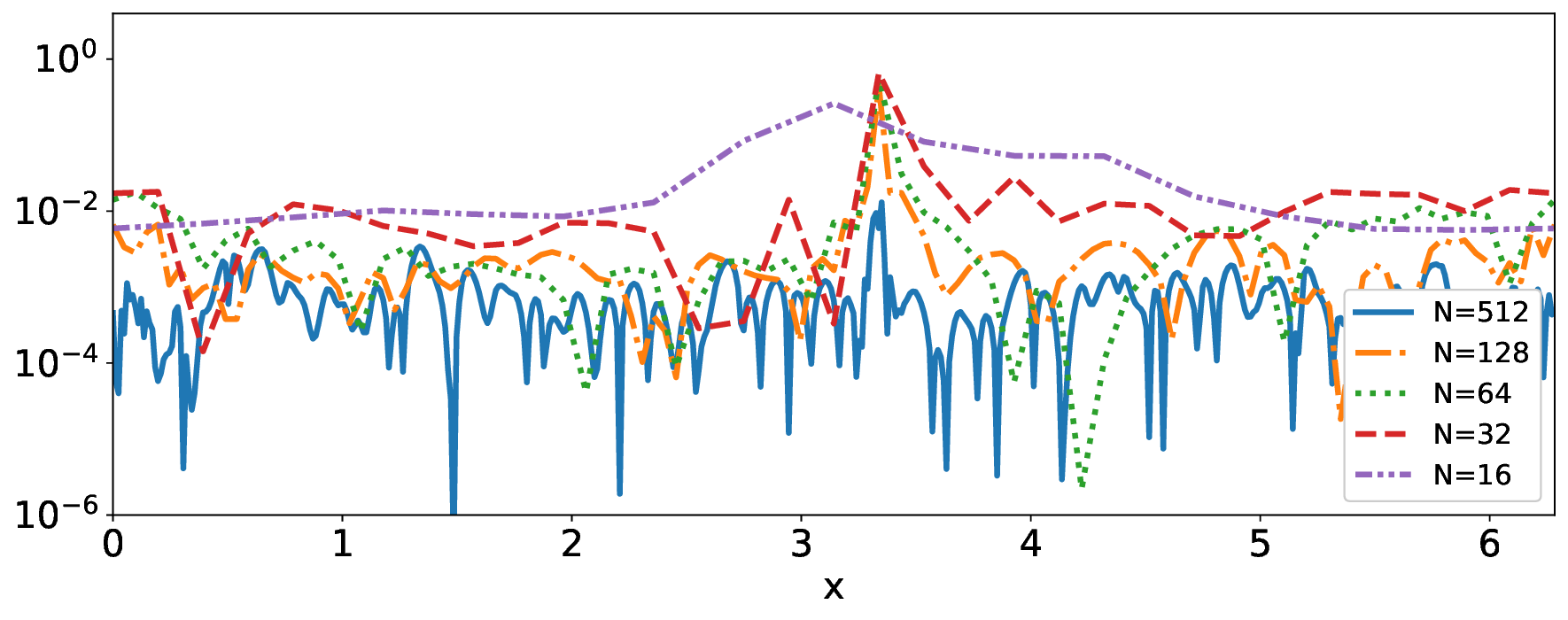}
		\caption{ CFN}
		\label{fig:ptErr_CFN_sparse}
	\end{subfigure}%
	\caption{Log-scale absolute error of the trained DNN model predictions to Example \ref{ex:burgers} at different sparsity levels when $t=3$. No observation error.}
	\label{fig:burg_sparse_ptErr}
\end{figure}

\subsubsection*{Case III: Dense and noisy observations}
In this testing environment the observations in \eqref{setup:data} now contain noise and are given by
\begin{equation}
\label{eq:data_noise_dense}
  \tilde{\mathbf{u}}^{(k)}(t_l) = \mathbf{u}^{(k)}(t_l) + \boldsymbol{\epsilon}_l^{(k)},\quad l = 1,\dots,L, \quad k = 1,\dots, N_{traj}.
\end{equation}
Here $\boldsymbol{\epsilon}_l^{(k)}$ is  i.i.d. Gaussian with zero mean and variance  $\sigma^2$. We test various $\sigma$ values scaled from the absolute value mean of the solution, {$\overline{|u|}$},
\begin{equation}\label{eq:noiselevel}
\sigma=a{\overline{|u|}},\quad  a\geq 0
\end{equation}
where the mean $\overline{u}$ is taken over the spatiotemporal domain.
We consider noise levels of 100\% , 50\%, 20\% and 10\%, that correspond to $\alpha=1,0.5,0.2,$ and $0.1$ respectively.
\begin{figure}[!h]
	\centering
	\begin{subfigure}[b]{0.24\textwidth}
		\includegraphics[width=\textwidth]{%
			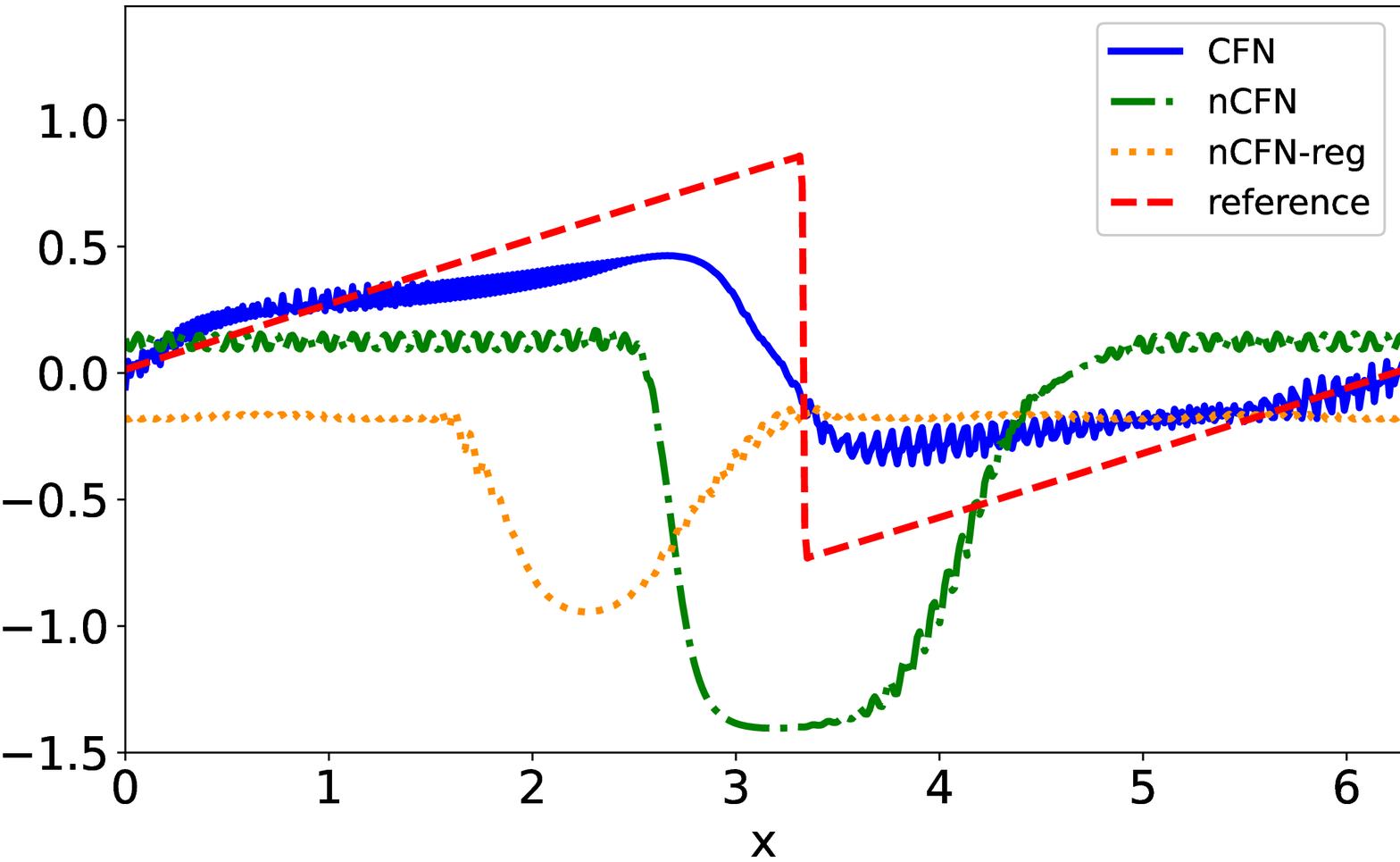}
		\caption{100\% noise }
		\label{fig:burgSNR1case2}
	\end{subfigure}%
	~
	\begin{subfigure}[b]{0.24\textwidth}
		\includegraphics[width=\textwidth]{%
			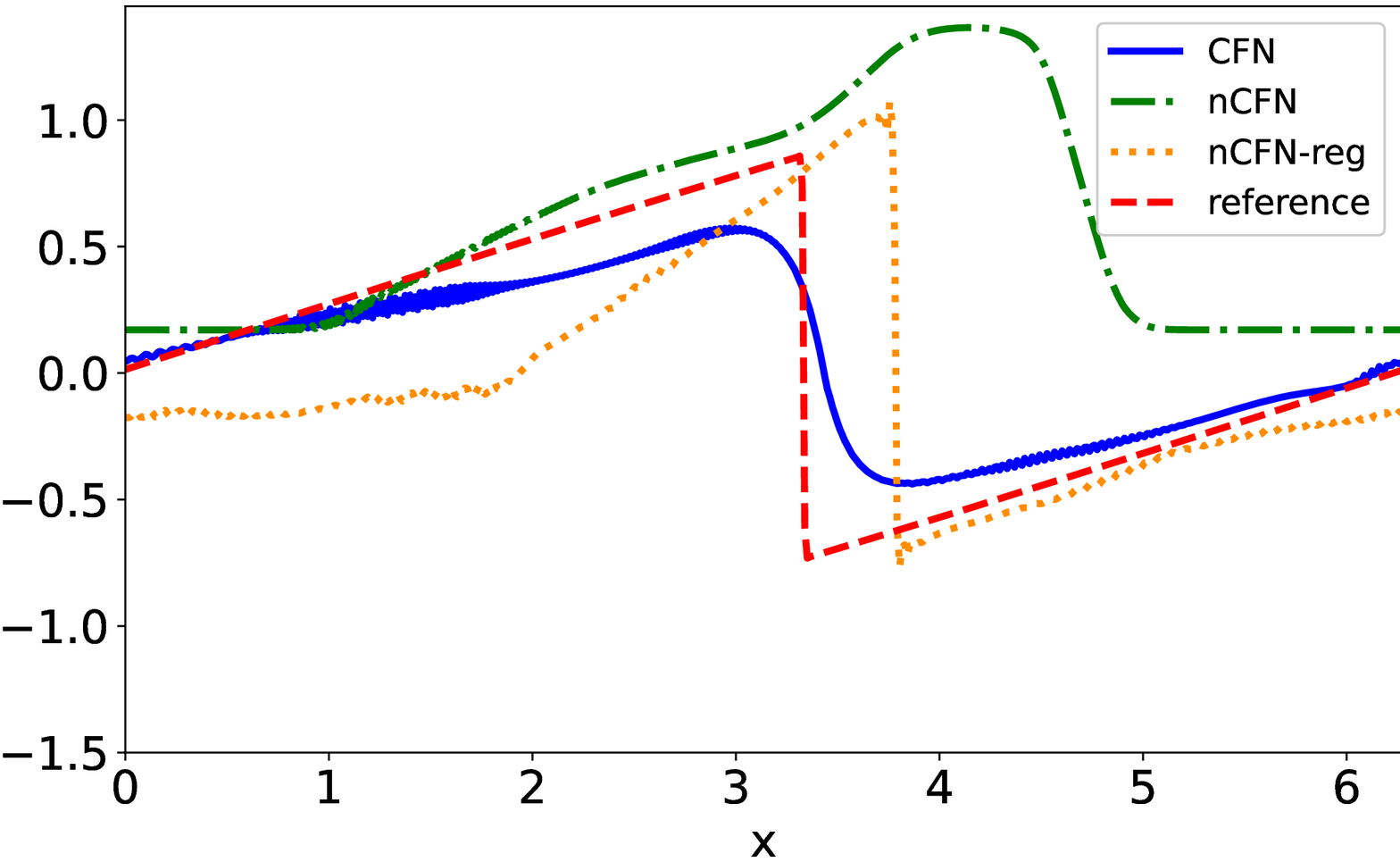}
		\caption{50\% noise }
		\label{fig:burgSNR2case2}
	\end{subfigure}%
	~
	\begin{subfigure}[b]{0.24\textwidth}
		\includegraphics[width=\textwidth]{%
			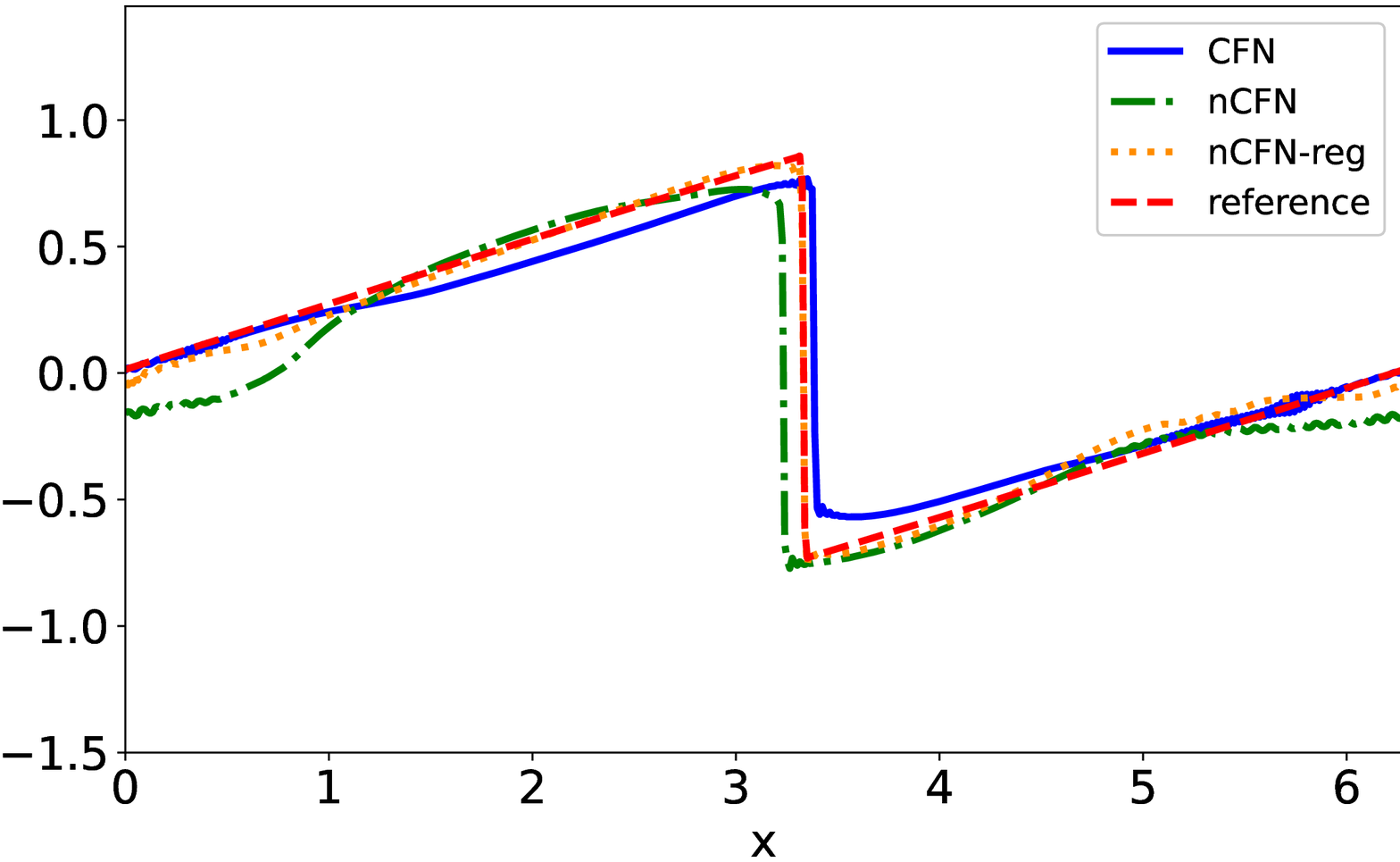}
		\caption{20\% noise }
		\label{fig:burgSNR5case2}
	\end{subfigure}%
	~
	\begin{subfigure}[b]{0.24\textwidth}
		\includegraphics[width=\textwidth]{%
			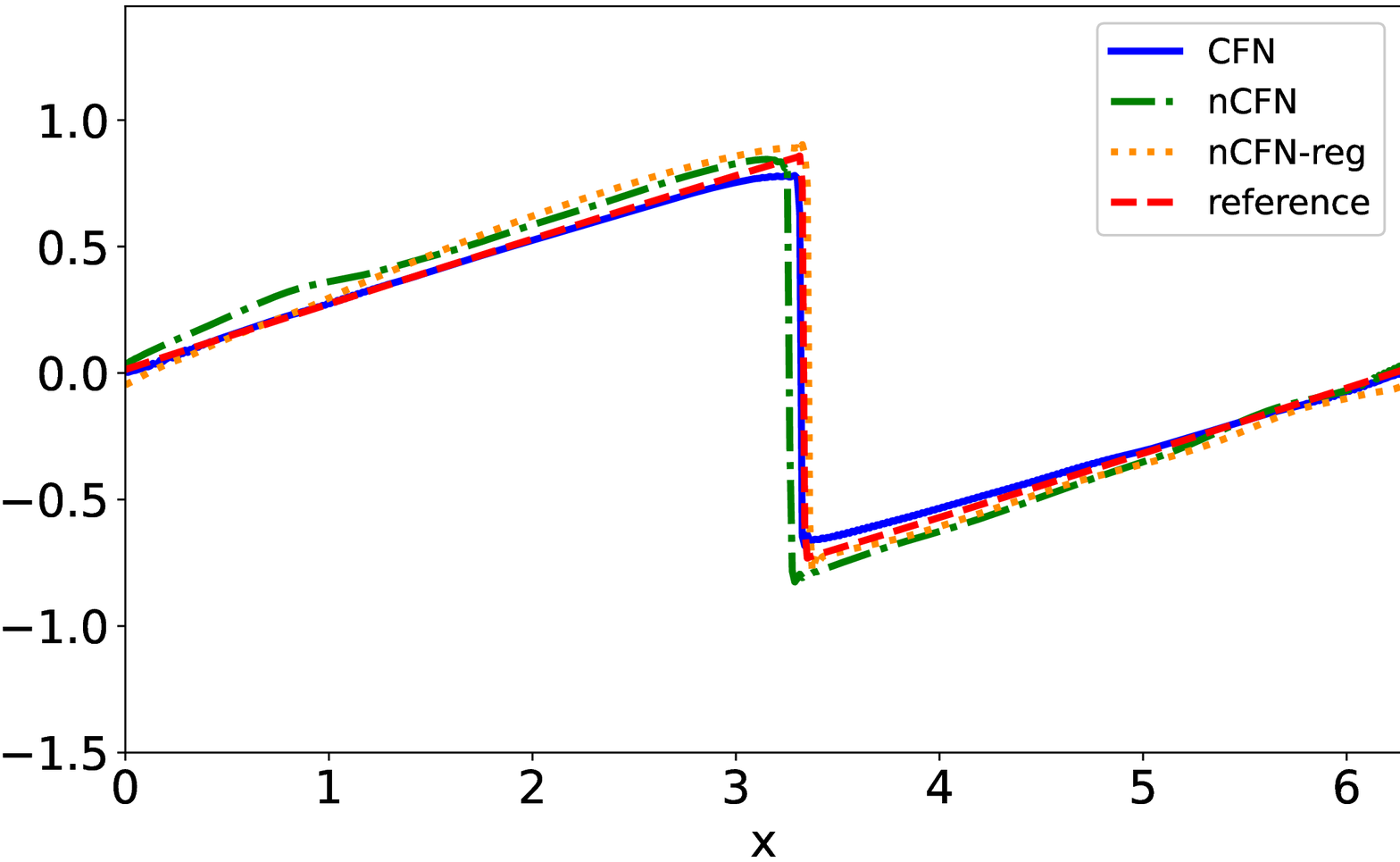}
		\caption{10\% noise}
		\label{fig:burgSNR10case2}
	\end{subfigure}%
	\caption{Comparison of the reference solution to Example \ref{ex:burgers} with the trained DNN model predictions at different noise levels using dense ($N=512$) observations when $t=3$.}
	\label{fig:burg_noisy}
\end{figure}

\begin{figure}[!htbp]
	\centering
	~
	\begin{subfigure}[b]{0.322\textwidth}
		\includegraphics[width=\textwidth]{%
			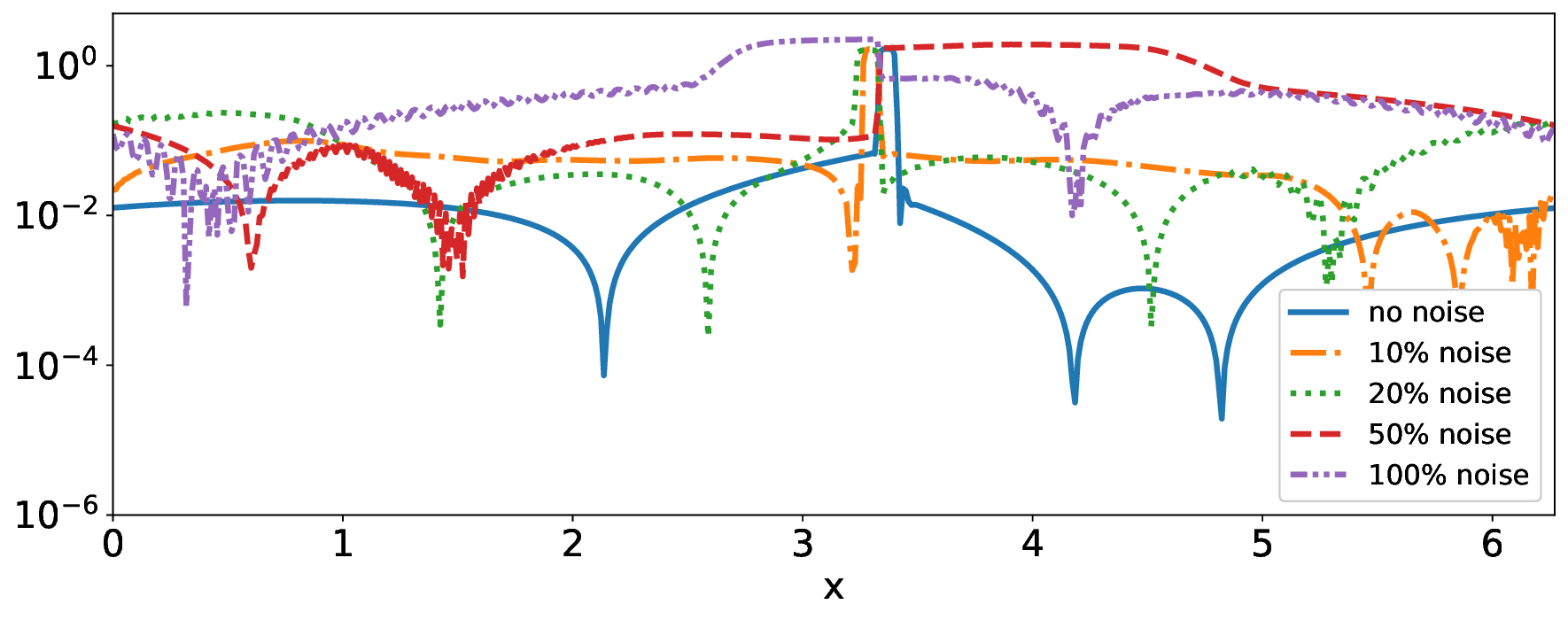}
		\caption{nCFN}
		\label{fig:ptErr_nCFN_noisy}
	\end{subfigure}%
	\begin{subfigure}[b]{0.322\textwidth}
		\includegraphics[width=\textwidth]{%
			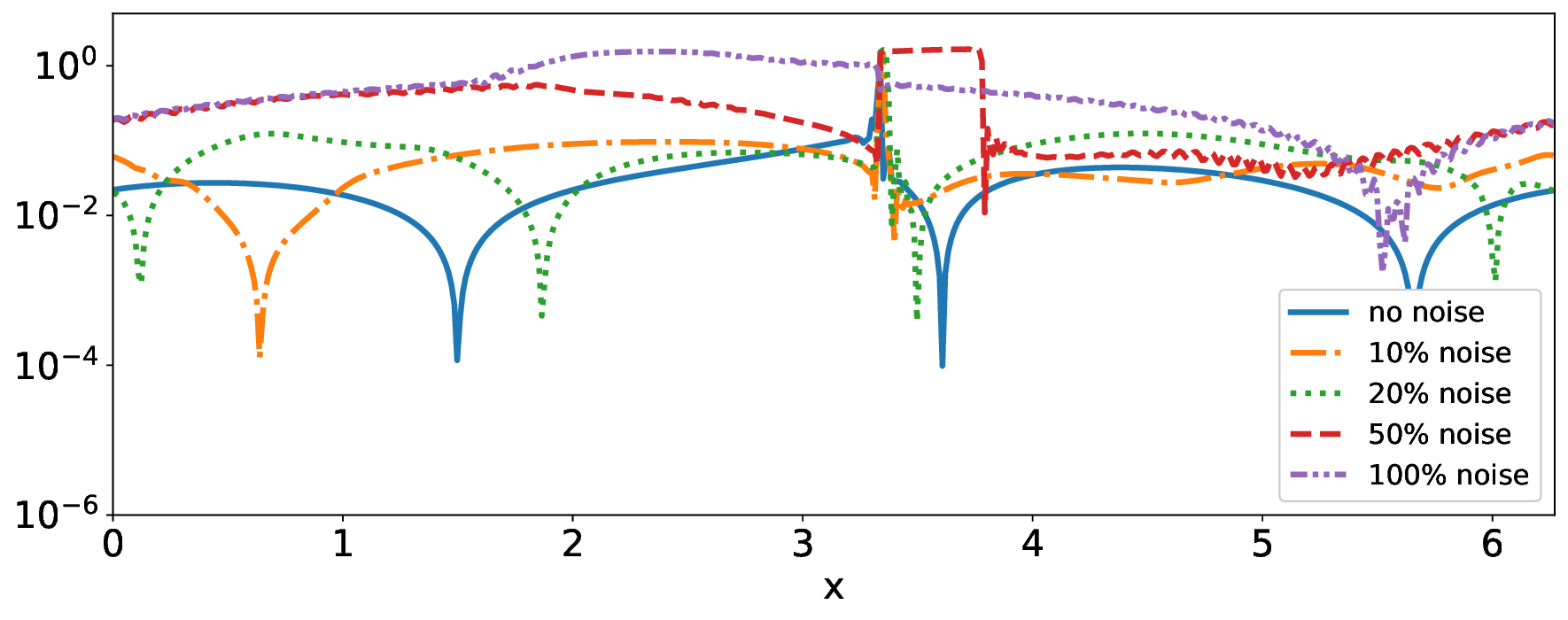}
		\caption{nCFN-reg}
		\label{fig:ptErr_nCFNreg_noisy}
	\end{subfigure}%
	\begin{subfigure}[b]{0.322\textwidth}
		\includegraphics[width=\textwidth]{%
			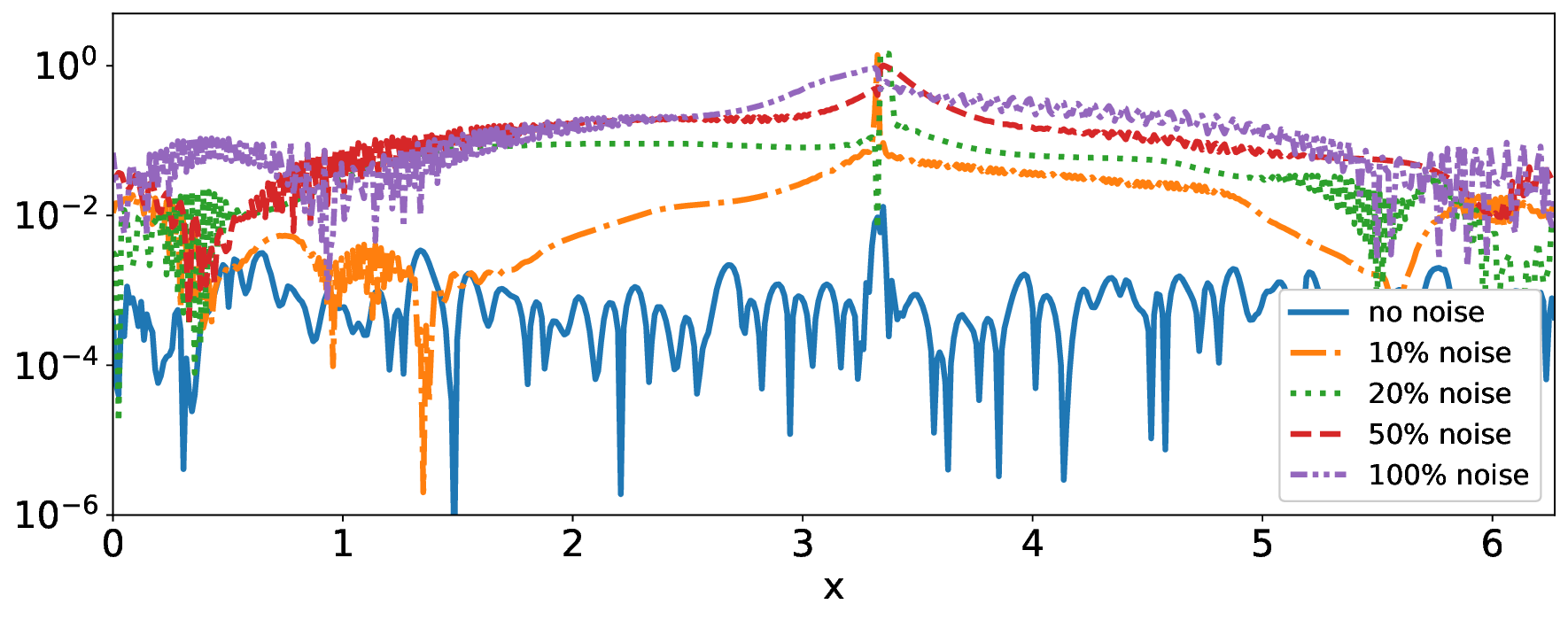}
		\caption{ CFN}
		\label{fig:ptErr_CFN_noisy}
	\end{subfigure}%
	\caption{Log-scale absolute error of the trained DNN model predictions to Example \ref{ex:burgers} using different noise levels using dense ($N=512$) observations when $t=3$.}
	\label{fig:burg_noisy_ptErr}
\end{figure}

The solution to Example \ref{ex:burgers} is presented in \figref{fig:burg_noisy} at time $t=3$ after the shock forms. Learning the underlying dynamics is challenging when the observations are noisy since the non-physical oscillatory behavior caused by the noise can influence the training process {(overfitting)}. Indeed  \figref{fig:burgSNR1case2} demonstrates that none of the solutions corresponding to any of the three training networks can capture the shock in high noise environments. However, the CFN method is the only network that captures the rough profile of the underlying solution. In contrast, solutions resulting from both nCFN and nCFN-reg deviate significantly from the reference solution. While all methods improve as the amount of noise decreases, CFN and nCFN-reg yield overall better results. In \figref{fig:burgSNR5case2} and \figref{fig:burgSNR10case2}, we observe that CFN and nCFN-reg capture the correct shock propagation speed (with some magnitude error). {\figref{fig:burg_noisy_ptErr} displays the pointwise errors for each of the methods  at time $t=3$, which are consistent to what is observed in \figref{fig:burgt4case1}.  That is, even in the ``ideal'' case, neither the nCFN nor the nCFN-reg can be adequately resolved.  Adding small amount of noise {which is comparable to the error already incurred therefore} does not affect the results.  {For the same reason, small amounts of noise can reduce the accuracy in the CFN case (since it is larger than the error produced for Case I).}  As noise is increased, the results for the nCFN and nCFN-reg method become meaningless -- ${\mathcal O}(1)$ in much of the domain.  The largest interval width of error surrounding the discontinuity is again seen in the nCFN case, which concurs with the results shown in \figref{fig:burg_noisy}.}  
Figures comparing the discrete conserved quantity remainder, \eqref{eq:conserve_u}, of each method are omitted for  Cases  II and III since the methods all generate the same general behavior pattern as what is shown for Case I in  \figref{fig:conserv_error}.


\subsection{Shallow water equations}
When combined with initial conditions given in \eqref{eq:SWE_IC}, Example \ref{ex:SWE} is known as the dam break problem which over time admits both shock and contact discontinuities.
\label{sec:SWE}
\begin{example}
	\label{ex:SWE}
	Consider the system of equations
	\begin{align}
		h_t + (vh)_x & = 0,\nonumber\\
		(hv)_t + (hv^2+\frac{1}{2}gh^2)_x & = 0,
	\label{eq:SWE}
	\end{align}
for $t > 0$ and $x\in(-5,5)$.  Here we use $g = 1$. We assume no flux boundary conditions and initial conditions given by
\begin{align}
	\label{eq:SWE_IC}
	h(x,0) &= \begin{cases}
		h_l, & \text{if $x\leq x_0$}, \\
		h_r, & \text{otherwise},
	\end{cases}\quad\quad
	 v(x,0) =\begin{cases}
		v_l, & \text{if $x\leq x_0$},\\
		v_r, & \text{otherwise},
	\end{cases}
\end{align}
where
\begin{align*}
	h_l &\sim U[2-\epsilon_{h_l}, 2+\epsilon_{h_l}], \quad \epsilon_{h_l}=0.2,\\
	h_r &\sim U[1-\epsilon_{h_r}, 1+\epsilon_{h_r}], \quad\epsilon_{h_l}=0.1,\\
	v_l, v_r, x_0  &\sim U[-\epsilon, \epsilon],\hspace{.3in} \epsilon=0.1.
\end{align*}
\end{example}

Example \ref{ex:SWE} describes the one-dimensional dam break problem  in which the initial heights of the water, $h_l$ and $h_r$, are different on each side of the dam, located at $x_0$ in our numerical experiments. After the dam breaks, a rarefaction wave forms and travels to the left of the dam, while a shock wave starts to propagate on the right.  The training data are observed at different time intervals up until time $t = 0.1$  and then used to train each of the three networks to predict the long term dynamics.

The $N_{traj} = 200$ training data sets of length $L = 20$ are generated by solving \eqref{eq:SWE} using CLAWPACK (HLLE Riemann Solver) for initial conditions given by
\begin{align*}
	h^{(k)} (x,0) &= \begin{cases}
		h_l^{(k)} , &\text{if $x\leq x_0^{(k)} $},\\
		h_r^{(k)} , & \text{otherwise},
	\end{cases}\quad\quad
	v^{(k)} (x,0) =\begin{cases}
		v_l^{(k)} , & \text{if $x\leq x_0^{(k)} $},\\
		v_r^{(k)} , & \text{otherwise},
	\end{cases}
\end{align*}
where
$$h_l^{(k)} \sim U[2-\epsilon_{h_l}, 2+\epsilon_{h_l}], \quad 
	h_r^{(k)} \sim U[1-\epsilon_{h_r}, 1+\epsilon_{h_r}], \quad
	v_l^{(k)}, v_r^{(k)}, x_0^{(k)}  \sim U[-\epsilon, \epsilon],$$
with  $\epsilon_{h_l}=0.2$, $\epsilon_{h_r} = .1$, $\epsilon = .1$ and $k = 1,\dots, N_{traj}$.
We obtain a reference solution using CLAWPACK using $N=1024$ so that $\Delta x = \frac{10}{1024}$. All figures shown for Example \ref{ex:SWE} correspond to \eqref{eq:SWE_IC} with $h_l = 3.5691196, h_r = 1.17867352, v_l = -0.06466697, v_r = -0.04519738, x_0 = 0.00383271$. While some parameter choices yield comparable solutions for each method, the CFN consistently outperforms the other techniques.

\subsubsection*{Case I: Dense and noise-free observations}
\begin{figure}[!h]
	\centering
\begin{subfigure}[b]{0.24\textwidth}
		\includegraphics[width=\textwidth]{%
			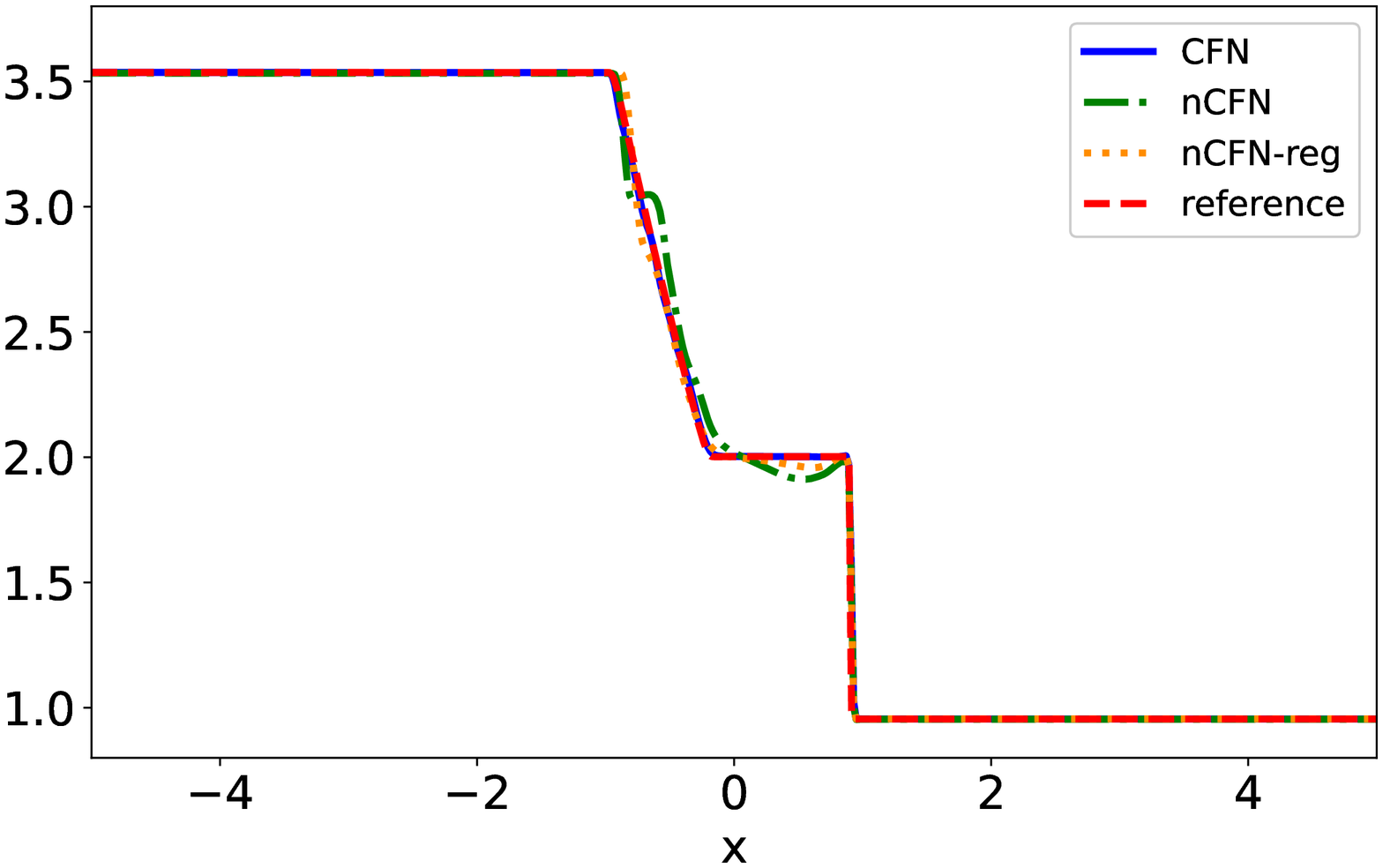}
		\caption{ $h, t=0.5$}
		\label{fig:SWEt1case1_h}
	\end{subfigure}%
	~
\begin{subfigure}[b]{0.24\textwidth}
	\includegraphics[width=\textwidth]{%
		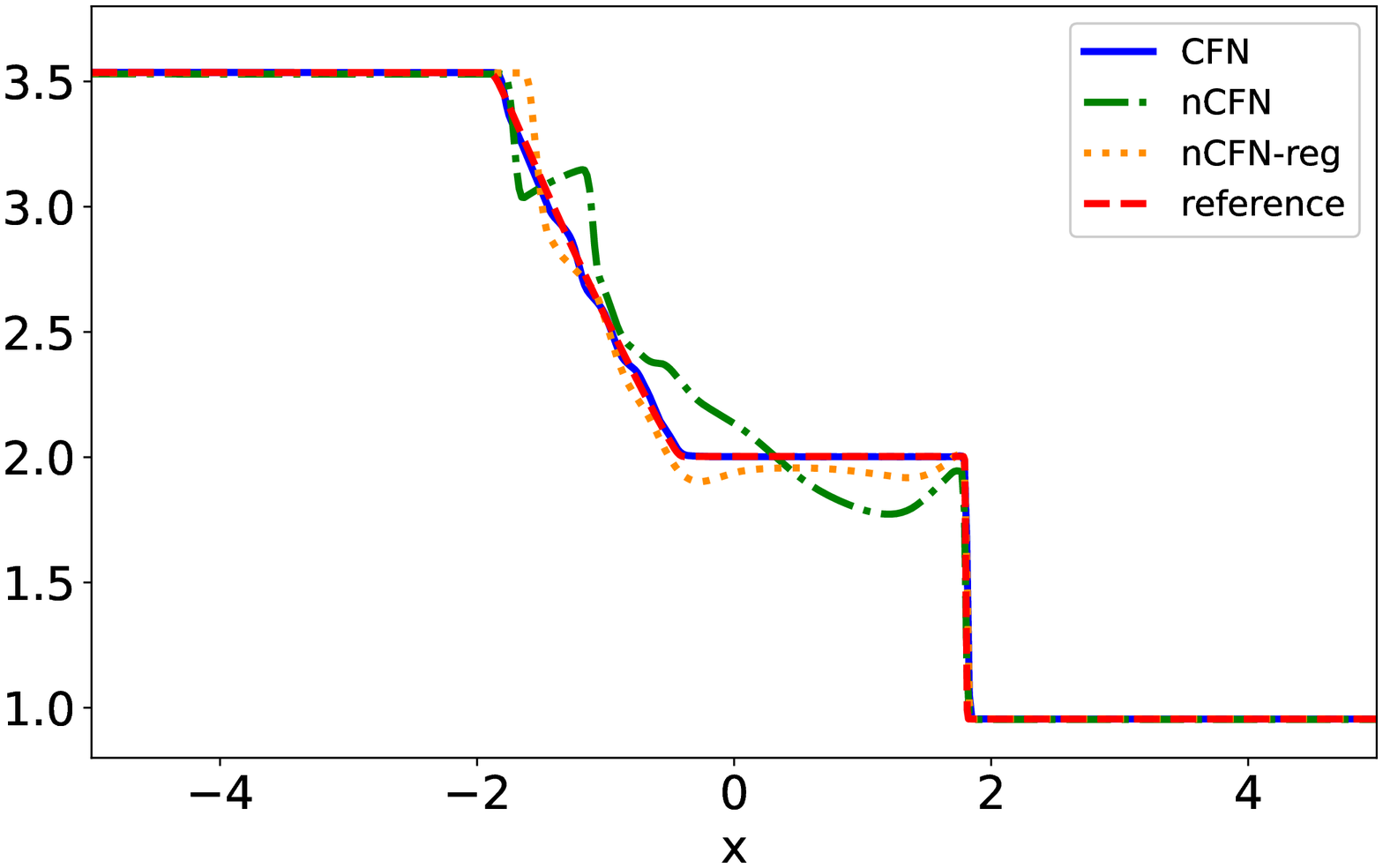}
	\caption{ $h, t=1$}
	\label{fig:SWEt2case1_h}
\end{subfigure}%
	~
\begin{subfigure}[b]{0.24\textwidth}
	\includegraphics[width=\textwidth]{%
		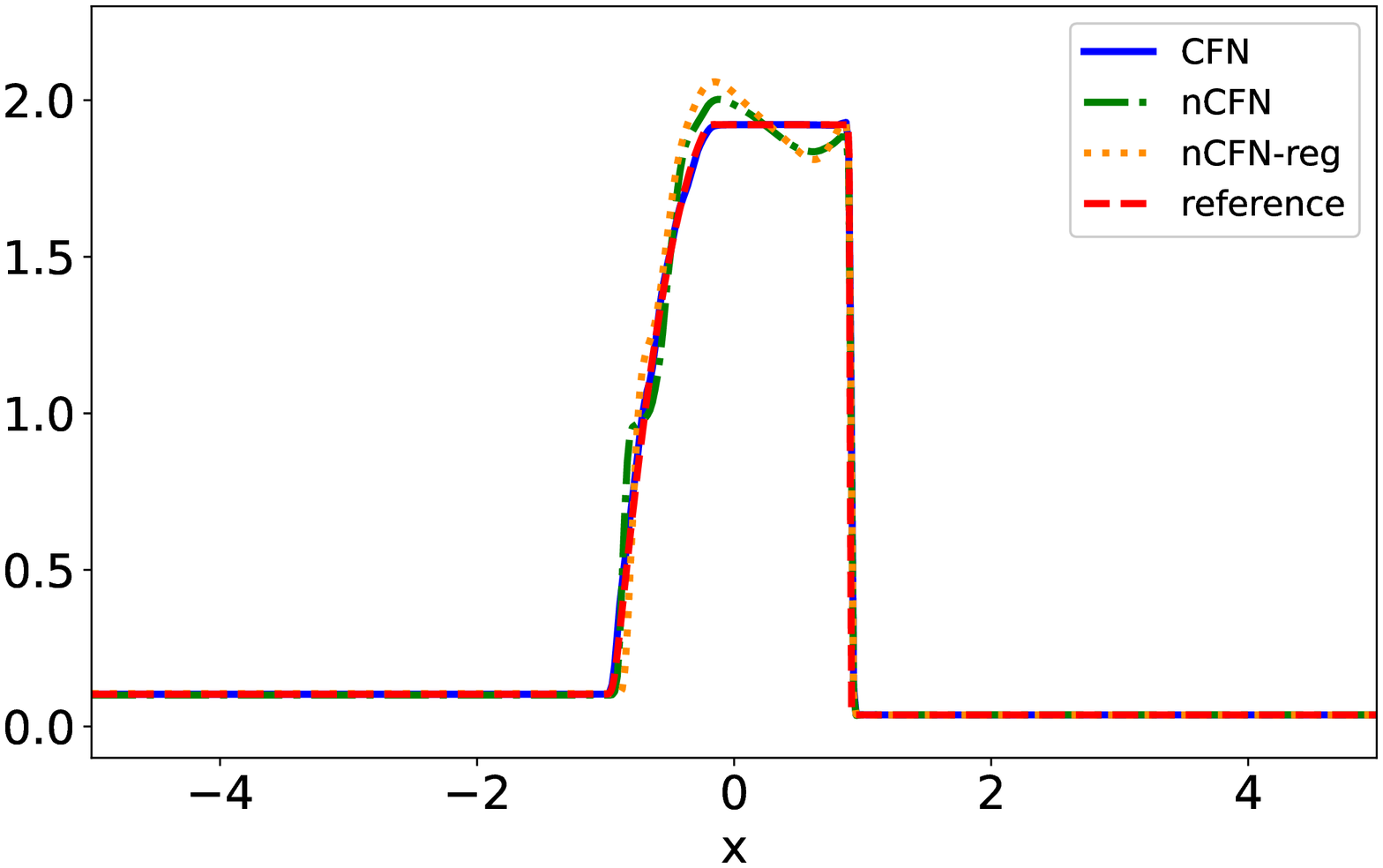}
	\caption{ $hv, t=0.5$}
	\label{fig:SWEt1case1_hv}
\end{subfigure}%
	\begin{subfigure}[b]{0.24\textwidth}
		\includegraphics[width=\textwidth]{%
		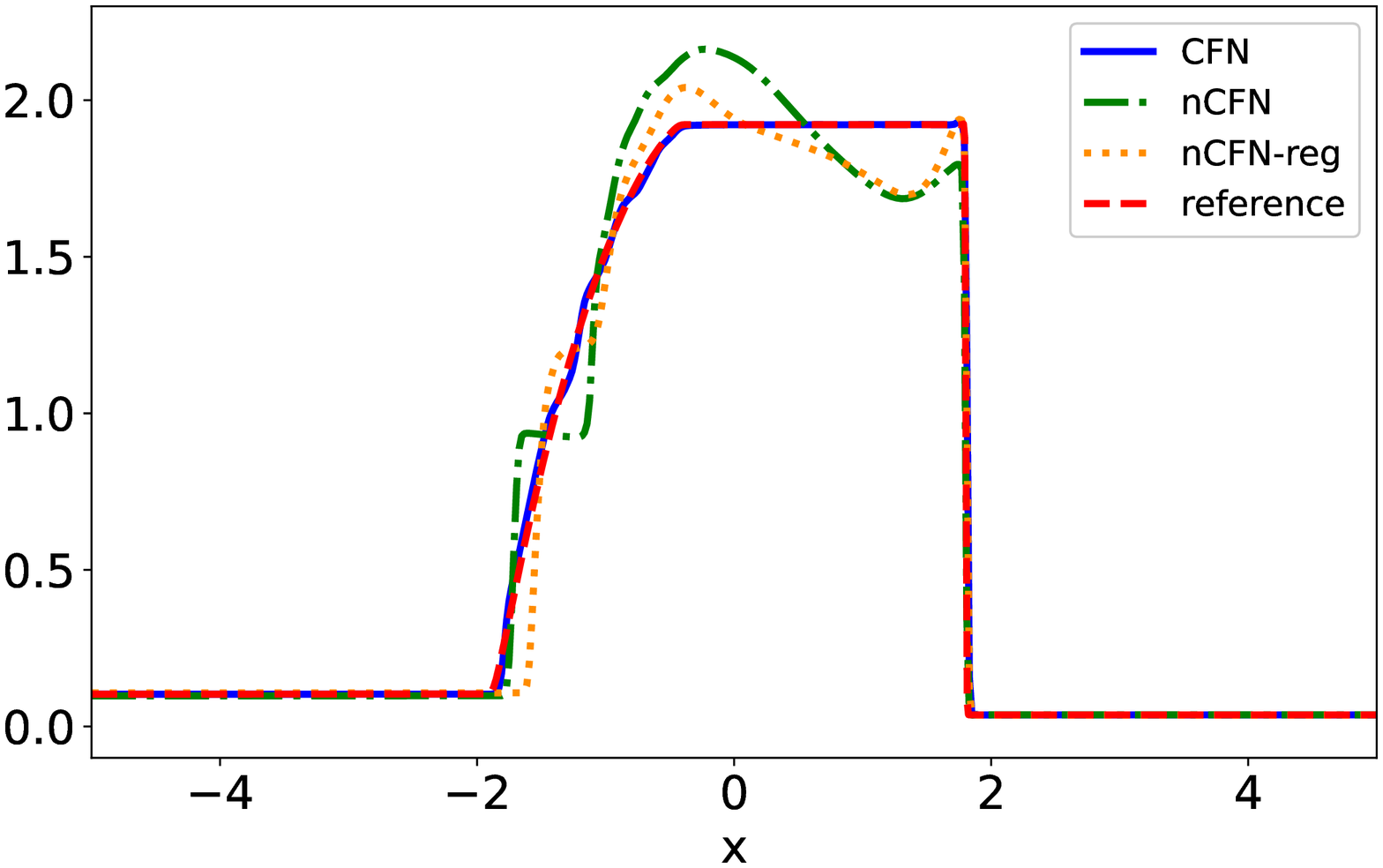}
		\caption{ $hv, t=1$}
		\label{fig:SWEt2case1_hv}
	\end{subfigure}%
	\caption{Comparison of the references solution to Example \ref{ex:SWE} and the trained DNN model predictions at different times for dense and noise-free observations.}
	\label{fig:SWE_solu}
\end{figure}
In the ideal environment we set $N=512$ so that $\Delta x = \frac{10}{512}$. We numerically impose the no flux boundary conditions using \eqref{eq:BCzeroth}. CLAWPACK is employed to simulate solutions  up to time $t=.1$ with data collections at time instances $t_l = l\Delta t$ for $l=1,\dots, 20$.  The time step  $\Delta t = 5 \times 10^{-3}$ is chosen to satisfy $\Delta t \leq \min\{\Delta t_{CLAW}\}$, where CLAWPACK determines $\Delta t_{CLAW}$ to guarantee stability for the solution in the given time domain.  

\begin{figure}[!h]
	\centering
	\begin{subfigure}[b]{0.30\textwidth}
		\includegraphics[width=\textwidth]{%
			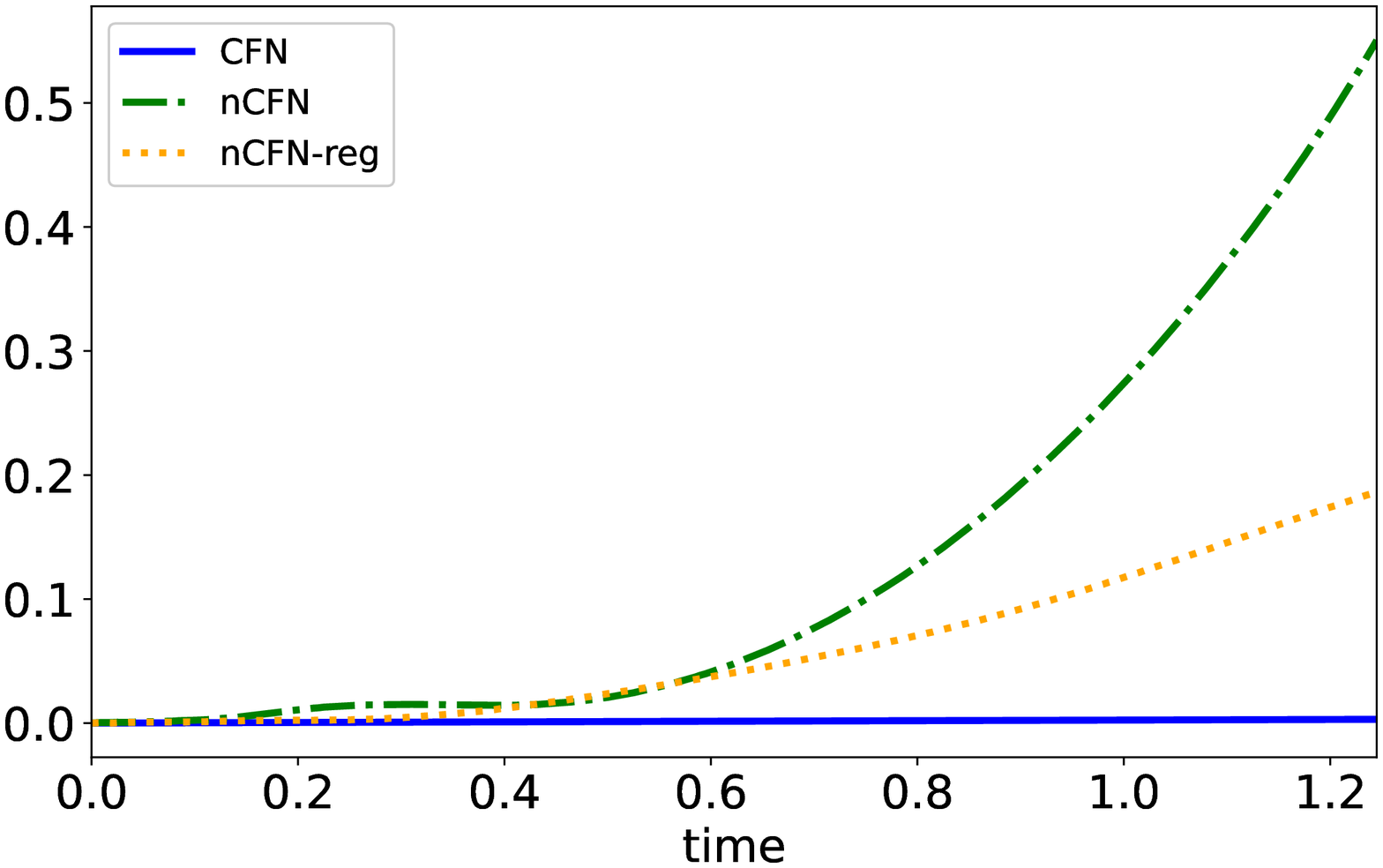}
		\caption{ $C(\mathbf{h})$}
		\label{fig:SWE_conserv_height}
	\end{subfigure}%
	~
	\begin{subfigure}[b]{0.30\textwidth}
	\includegraphics[width=\textwidth]{%
		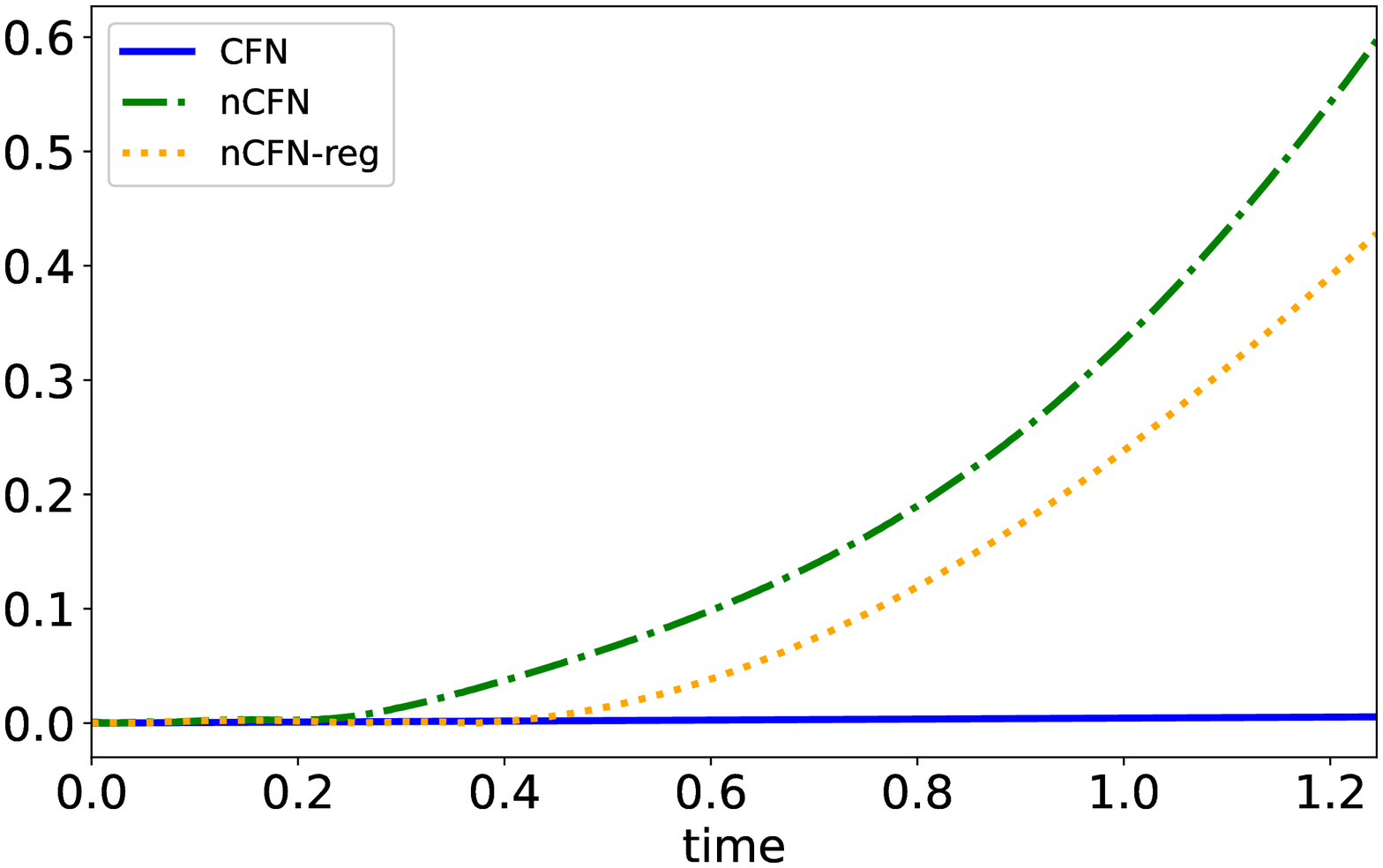}
		\caption{ $C(\mathbf{h}\odot\mathbf{v})$}
	\label{fig:SWE_conserv_momentumt}
\end{subfigure}%
			\caption{Discrete conserved quantity remainder given by \eqref{eq:conserve_u} of each method for Example  \ref{ex:SWE}.  (a) $C(\mathbf{h})$ and (b) $C(\mathbf{h}\odot\mathbf{v})$.}
	\label{fig:SWE_conserv}
\end{figure}
\figref{fig:SWE_solu} compares the numerical solutions at time $t=0.5$ and $t = 1$, both of which extend past the training time.  While all methods capture the main features of the solution at $t = 0.5$, it is evident that the CFN yields the most accurate results.  The errors in both the nCFN and nCFN-reg solutions are significantly larger when $t = 1$, and the rarefaction wave structure is not discernible in the nCFN case. We determine the conservation of each method by calculating \eqref{eq:conserve_u} for $\mathbf{h}$ and $\mathbf{h}\odot\mathbf{v}$, where $\odot$ denotes elementwise multiplication, and show the results in \figref{fig:SWE_conserv}.  As in the case for Burgers equation, only the CFN method is conservative.

\subsubsection*{Case II: Sparse and noise-free observations}
To simulate this environment we use CLAWPACK to solve Example \ref{ex:SWE} on coarser grids, respectively $N = 64$ and $N= 128$, to obtain the training data collected at $t_l = l\Delta t$, $l = 1,\dots,L$, where $L = 20$.  For consistency we again choose $\Delta t = 0.005$ so that the training trajectory final time is $t = 0.1$.
\begin{figure}[!h]
	\centering
\begin{subfigure}[b]{0.24\textwidth}
		\includegraphics[width=\textwidth]{%
			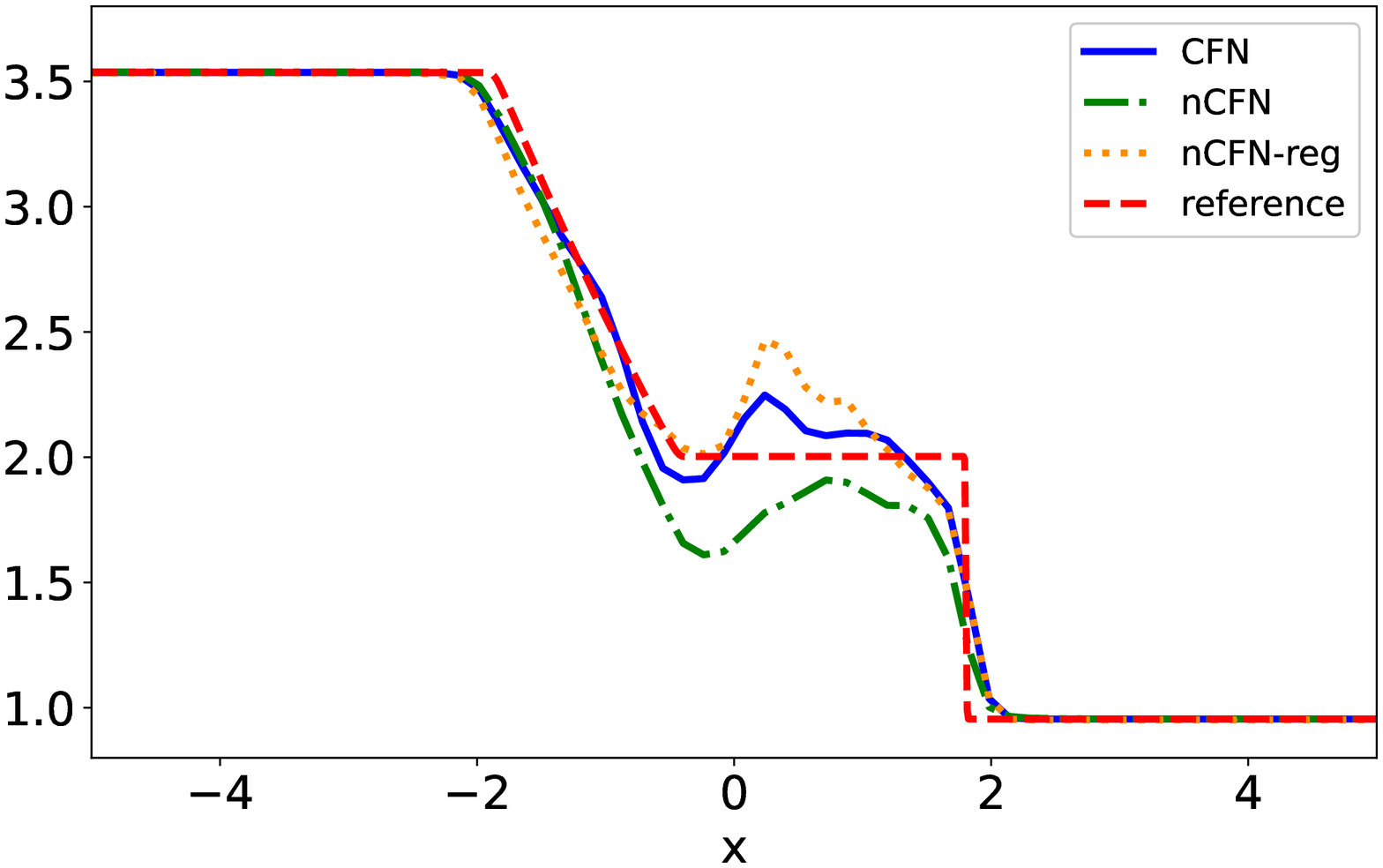}
		\caption{$h, N = 64$}
		\label{fig:SWE_Nx64_h}
	\end{subfigure}%
\begin{subfigure}[b]{0.24\textwidth}
		\includegraphics[width=\textwidth]{%
			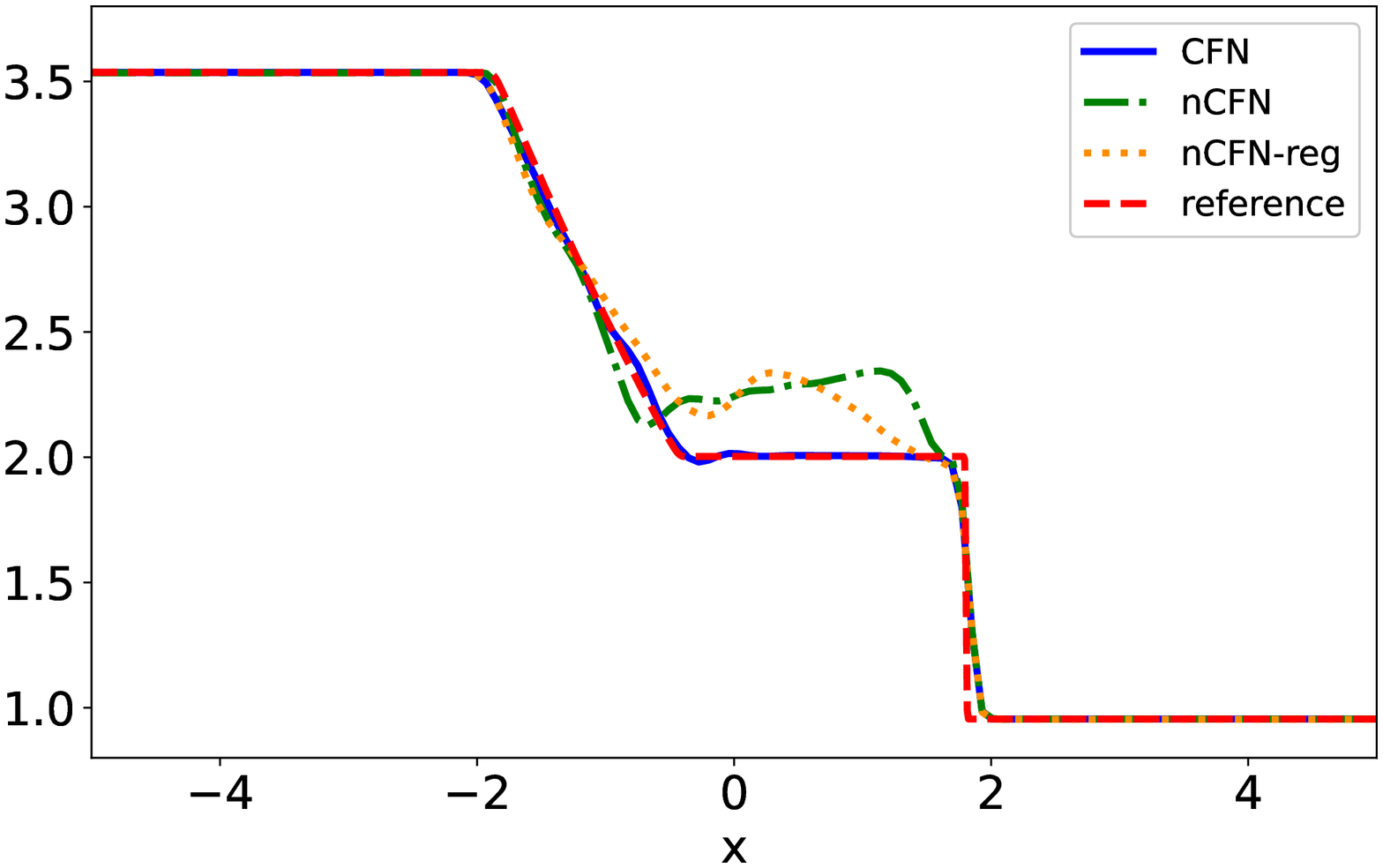}
		\caption{$h, N = 128$}
		\label{fig:SWE_Nx128_h}
	\end{subfigure}%
	~
	\begin{subfigure}[b]{0.24\textwidth}
		\includegraphics[width=\textwidth]{%
		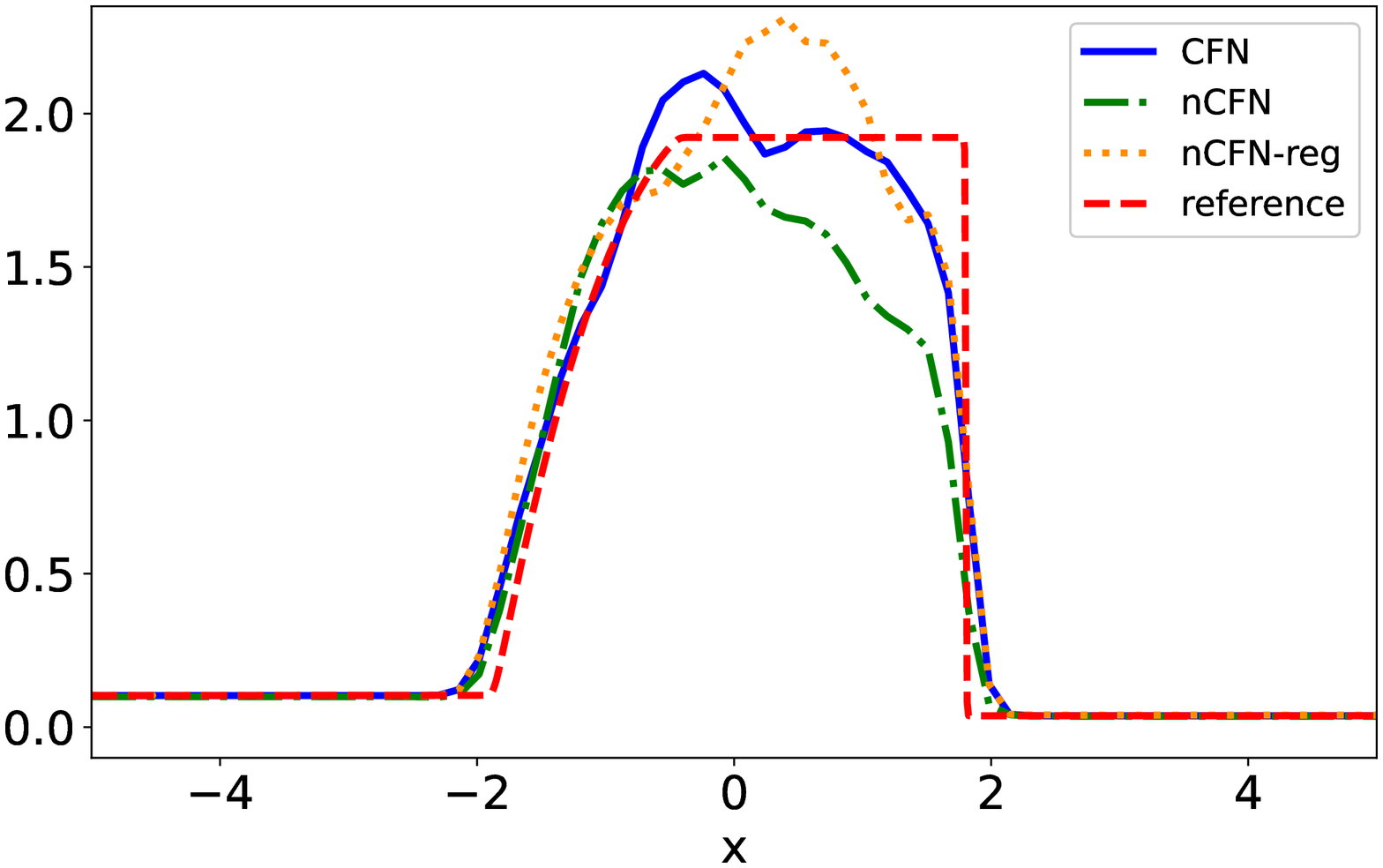}
		\caption{$hv, N = 64$}
		\label{fig:SWE_Nx64_hv}
	\end{subfigure}%
	~
\begin{subfigure}[b]{0.24\textwidth}
		\includegraphics[width=\textwidth]{%
		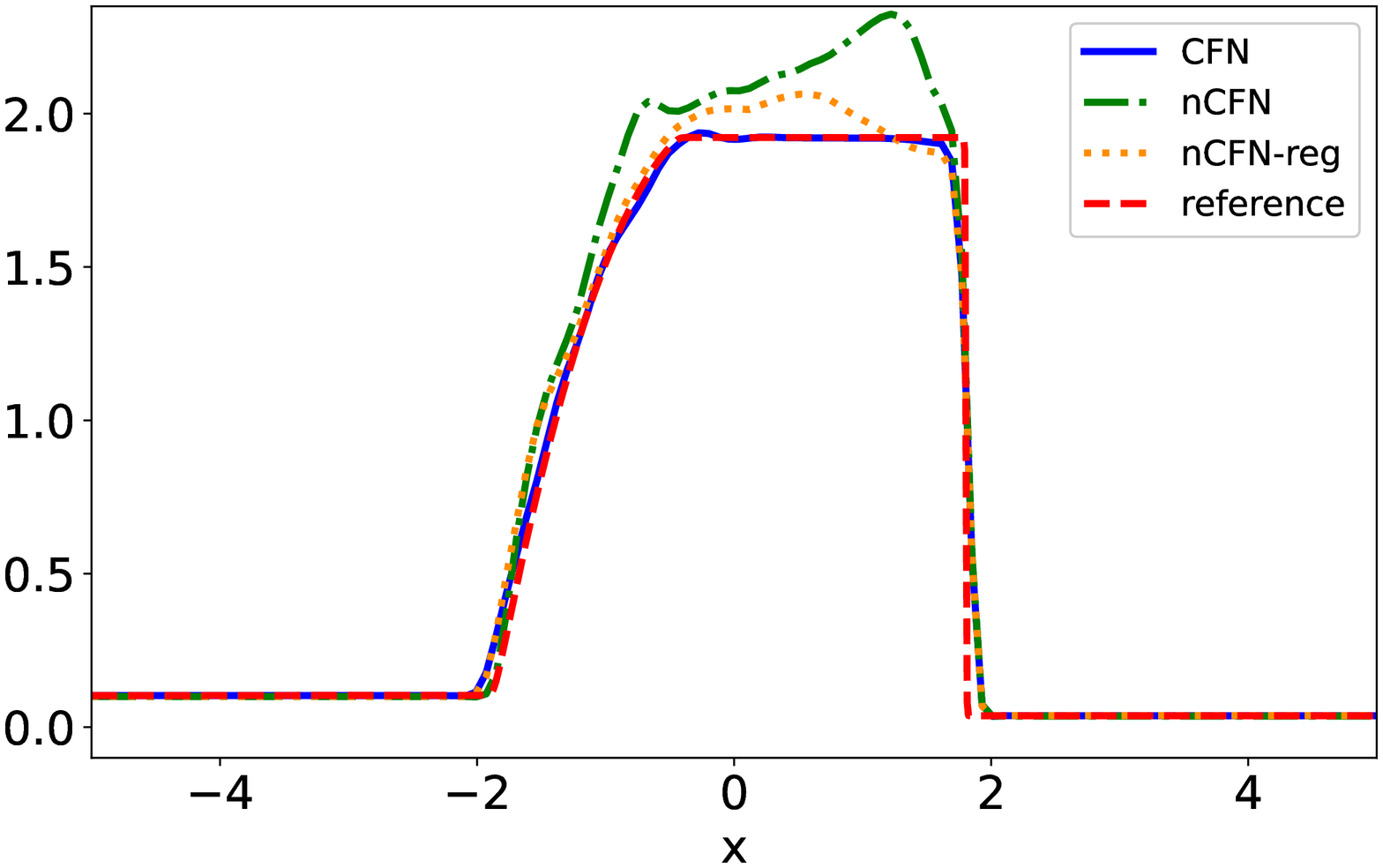}
		\caption{$hv, N = 128$}
		\label{fig:SWE_Nx128_hv}
	\end{subfigure}%
	~
	\caption{Comparison of the reference solutions to Example \ref{ex:SWE} and the trained DNN model predictions using noise free observations at  time $t = 1$ for $N = 64$ and $128$.}
	\label{fig:SWE_sparse}
\end{figure}

\figref{fig:SWE_sparse} compares the solutions using CFN, nCFN and nCFN-reg with the reference solution (again defined as the CLAWPACK solution with $N=1024$) at time $t=1$.  For  $N=64$ it is apparent that none of the methods are able to accurately learn the system dynamics, and large fluctuations are particularly noticeable in the region between the rarefaction and shock wave.   For $N = 128$ we observe improvement for all models.  The CFN clearly yields the  most accurate results,  and is the only method able to capture the structure of the rarefaction wave.  This is not surprising since we already observed in \figref{fig:SWE_solu} that $N = 512$ did not provide enough resolution, even for the nCFN-reg case.  Thus we see the importance of training the network using the flux form.

\subsubsection*{Case III: Dense and noisy observations}
In this case the training data are given by
\begin{eqnarray}
	\label{eq:swe_noise_dense}
	\tilde{\mathbf{h}}^{(k)}(t_l) &=& \mathbf{h}^{(k)}(t_l) + \boldsymbol{\epsilon}_l^{(k)}, \nonumber\\
       \widetilde{(\mathbf{h}^{(k)}\odot\mathbf{v}^{(k)})}(t_l) &=& (\mathbf{h}^{(k)}\odot\mathbf{v}^{(k)})(t_l) + \boldsymbol{\eta}_l^{(k)},
\end{eqnarray}
for $l=1,\dots,L$  and $k = 1,\dots, N_{traj}.$
The added noise $\boldsymbol{\epsilon}_l^{(k)}$ and $\boldsymbol{\eta}_{l}^{(k)}$ are  i.i.d. Gaussian with zero mean and variance determined using various noise values based on the mean of $\mathbf{u}$ \eqref{eq:noiselevel}.  We again consider noise levels corresponding to 100\%, 50\%, 20\%, and 10\%.

The solution profiles for height and momentum in Example \ref{ex:SWE} obtained using the different network constructions are shown in \figref{fig:SWE_height_noisy}. It is apparent that all three methods yield significant diffusion in high noise environments.  It is noteworthy that when the amount of noise is at 20\%, both the nCFN and nCFN-reg methods produce solutions that seem to increase (rather than diffuse) energy, suggesting that these methods are learning noise-related dynamics.  In this regard, the CFN method appears to be the most robust, meaning that along with the overall improved quality of the solution with decreasing amounts of noise,  the solution itself behaves consistently as  a function of the noise level, with less diffusion apparent as the amount of noise decreases. 
\begin{figure}[!h]
\centering
\begin{subfigure}[b]{0.24\textwidth}
	\includegraphics[width=\textwidth]{%
		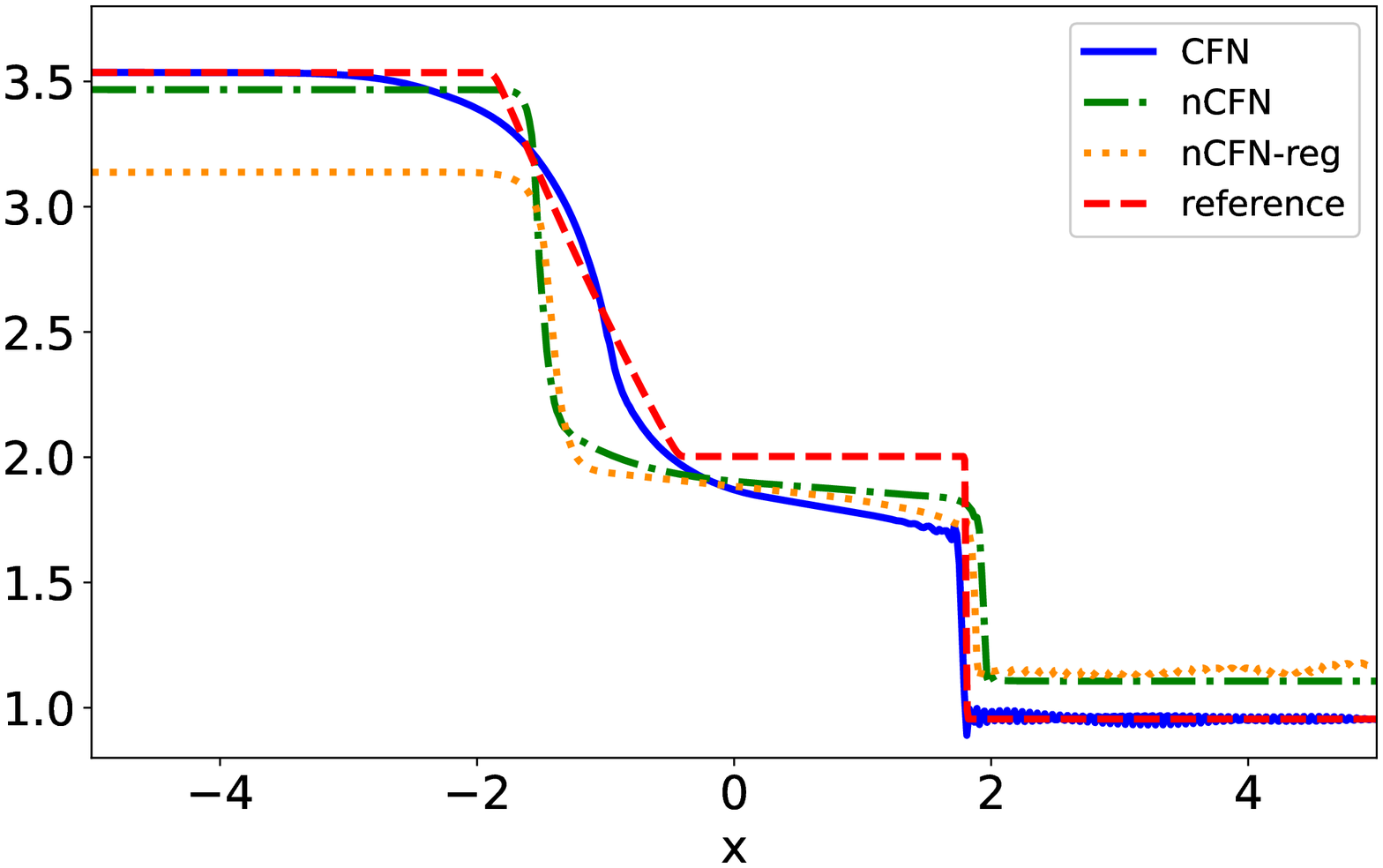}
	\caption{100\% noise}
	\label{fig:SWEsnr1case2_h}
\end{subfigure}%
~
\begin{subfigure}[b]{0.24\textwidth}
	\includegraphics[width=\textwidth]{%
		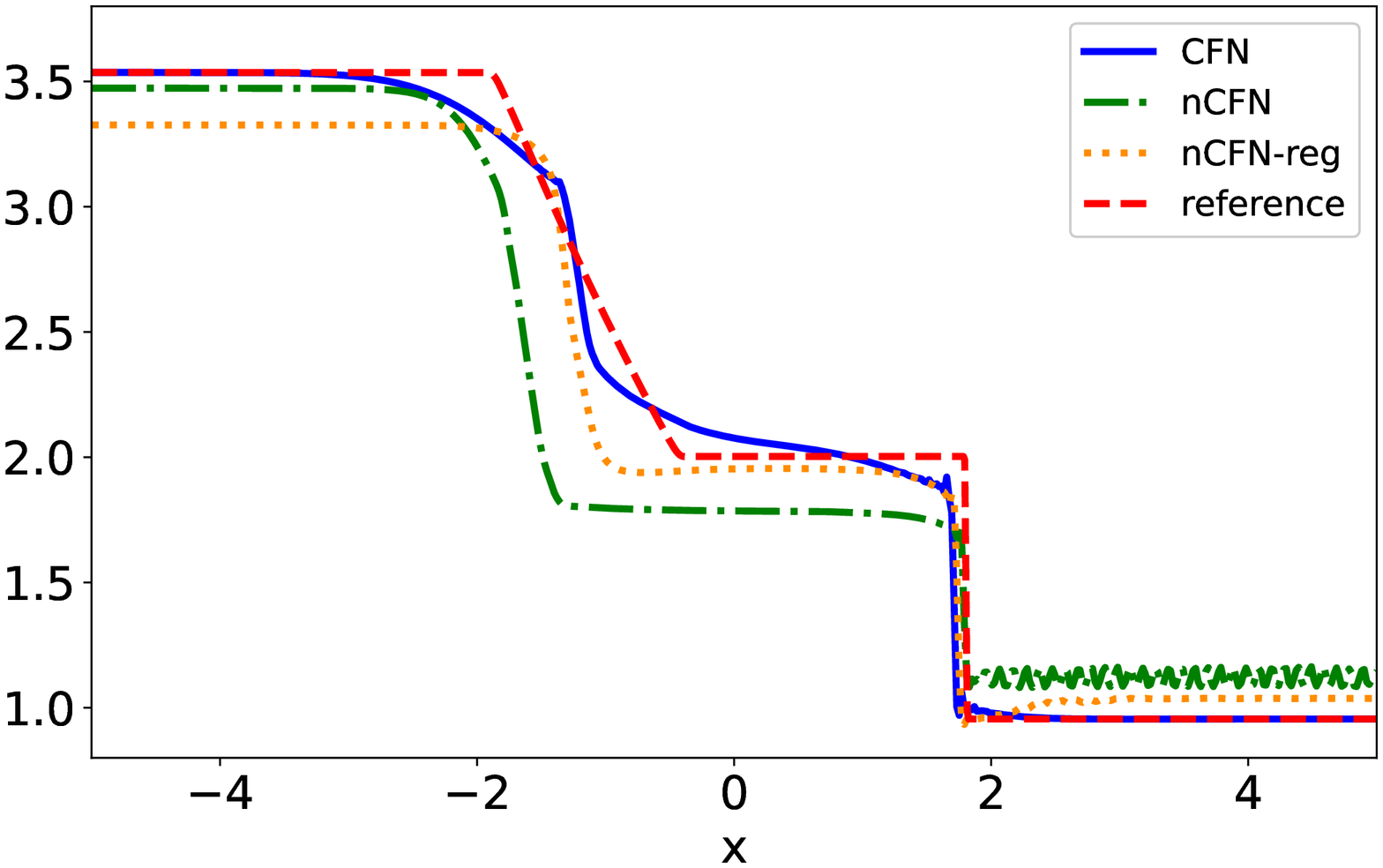}
	\caption{50\% noise}
	\label{fig:SWEsnr2case2_h}
\end{subfigure}%
\begin{subfigure}[b]{0.24\textwidth}
	\includegraphics[width=\textwidth]{%
		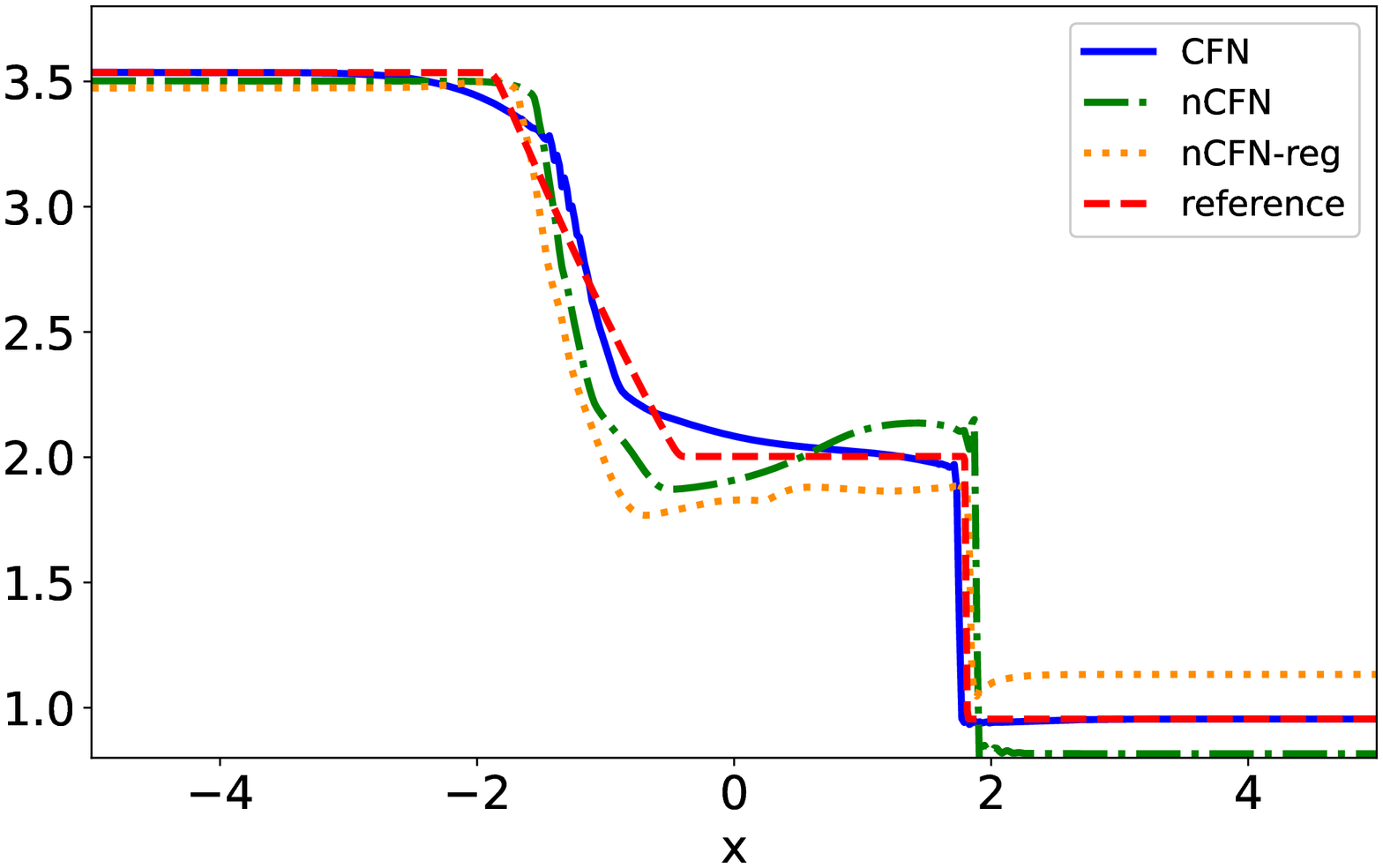}
	\caption{20\% noise}
	\label{fig:SWEsnr5case2_h}
\end{subfigure}
~
\begin{subfigure}[b]{0.24\textwidth}
	\includegraphics[width=\textwidth]{%
	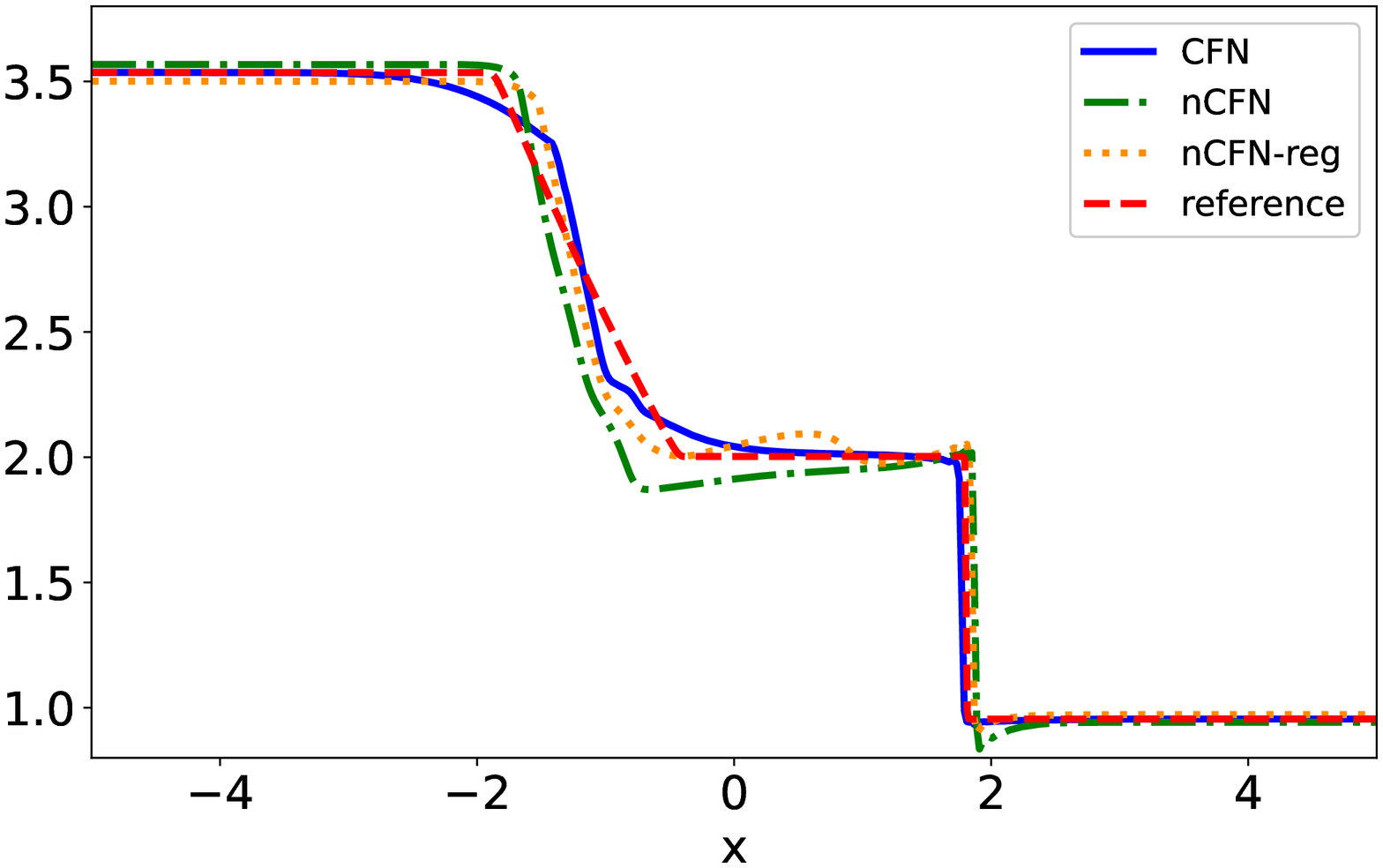}
	\caption{10\% noise}
	\label{fig:SWEsnr10case2_h}
\end{subfigure}%
\\
\begin{subfigure}[b]{0.24\textwidth}
	\includegraphics[width=\textwidth]{%
		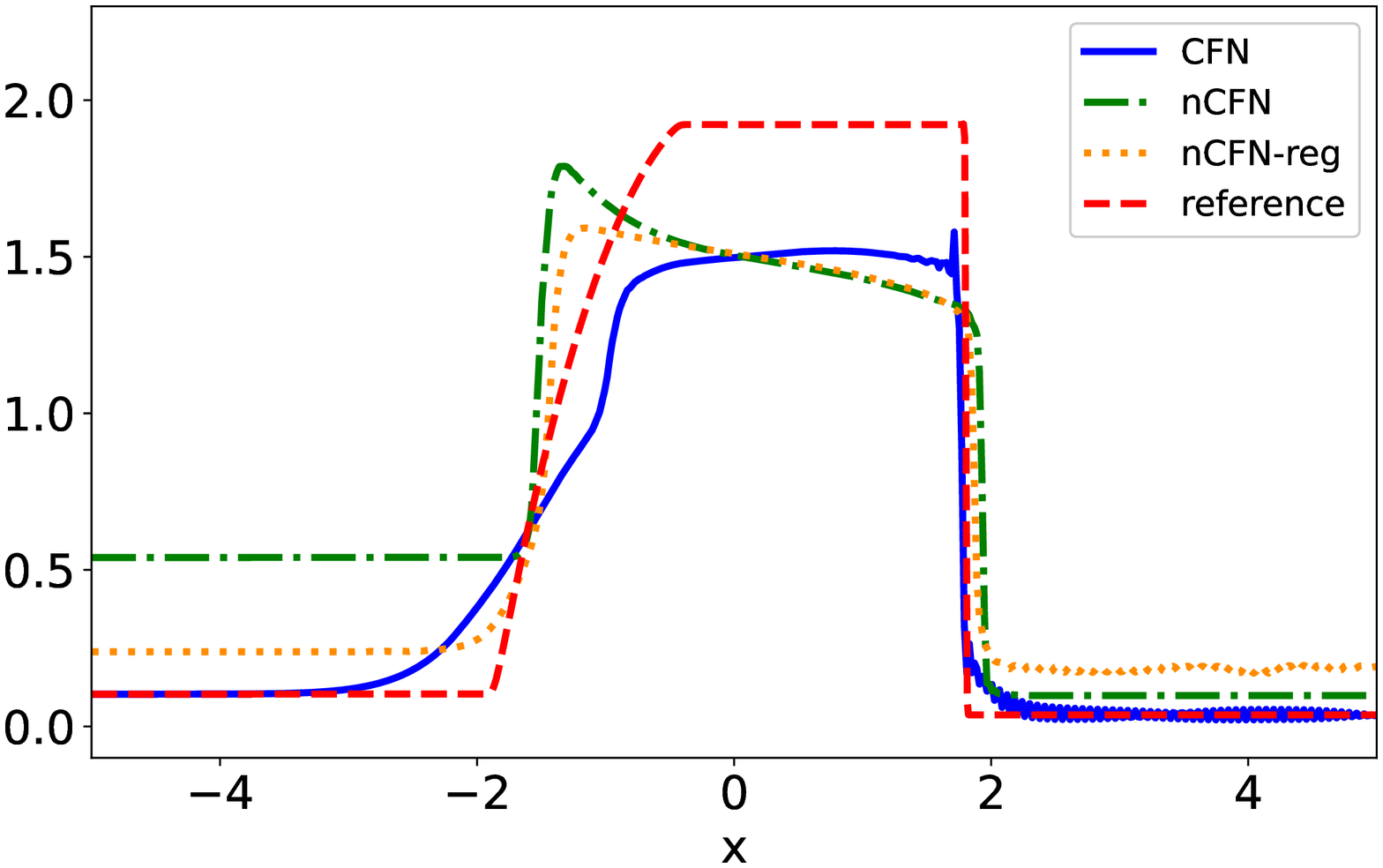}
	\caption{100\% noise }
	\label{fig:SWEsnr1case2_hv}
\end{subfigure}%
~
\begin{subfigure}[b]{0.24\textwidth}
	\includegraphics[width=\textwidth]{%
		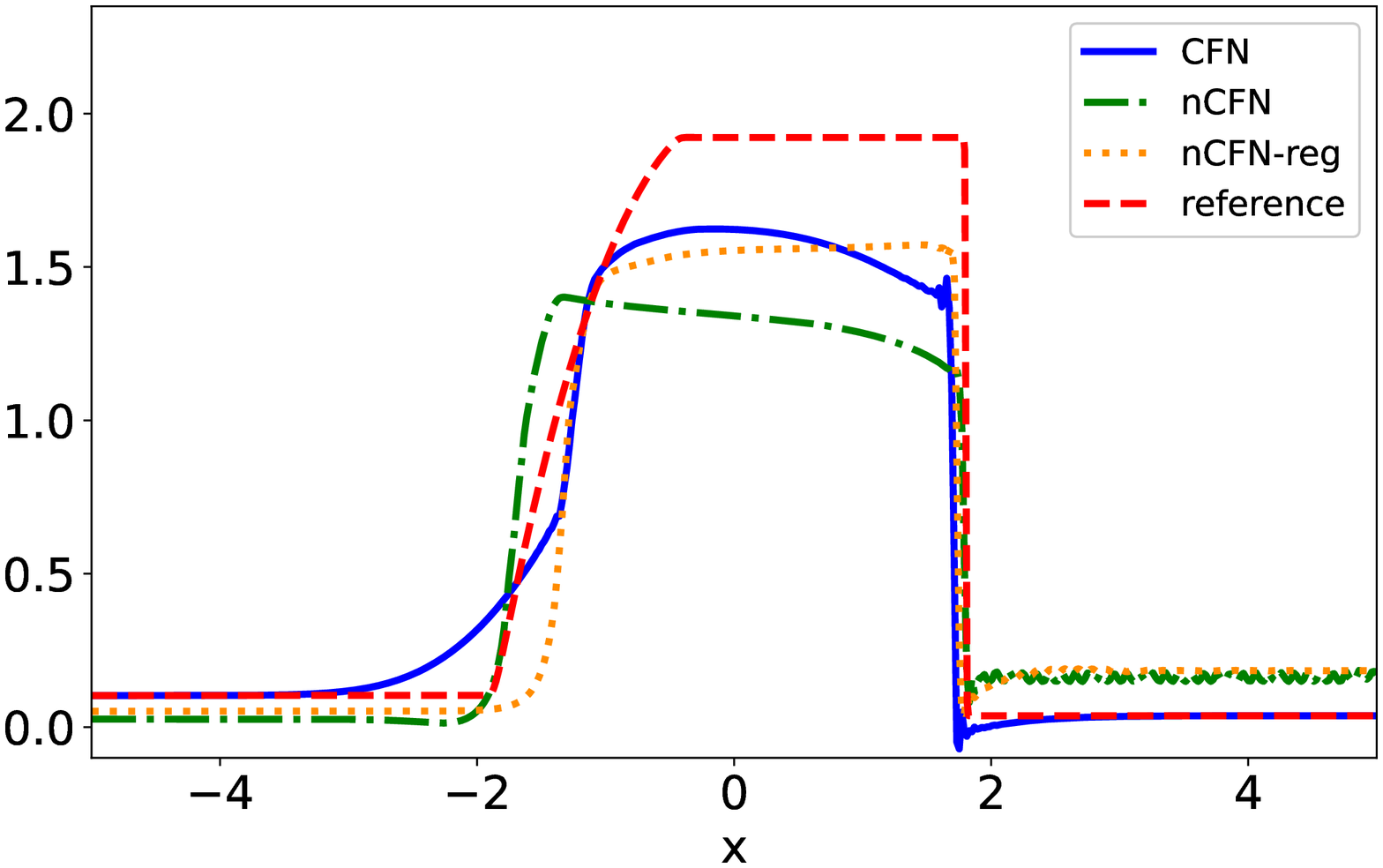}
	\caption{50\% noise}
	\label{fig:SWEsnr2case2_hv}
\end{subfigure}%
~
\begin{subfigure}[b]{0.24\textwidth}
	\includegraphics[width=\textwidth]{%
		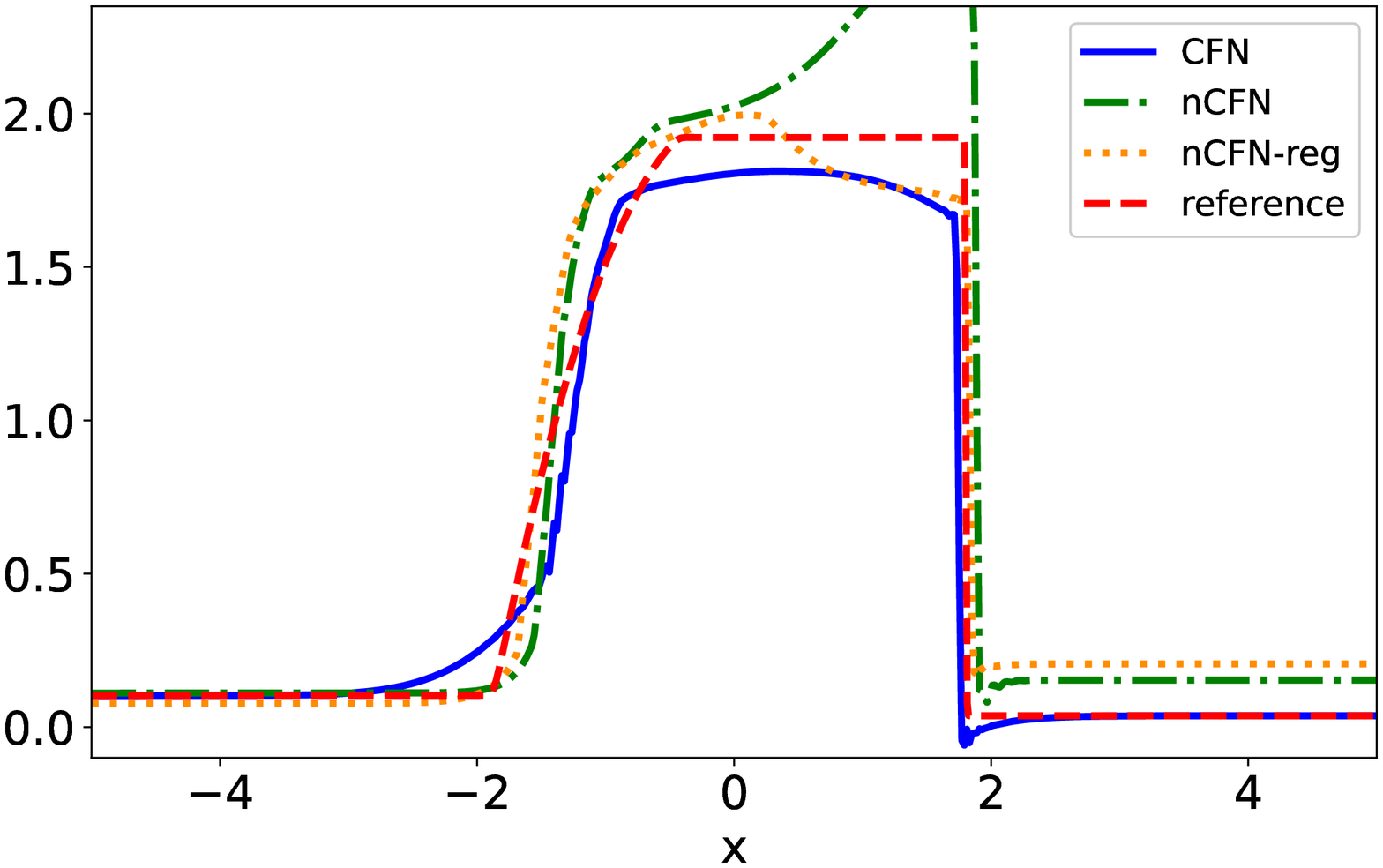}
	\caption{20\% noise }
	\label{fig:SWEsnr5case2_hv}
\end{subfigure}%
~
\begin{subfigure}[b]{0.24\textwidth}
	\includegraphics[width=\textwidth]{%
	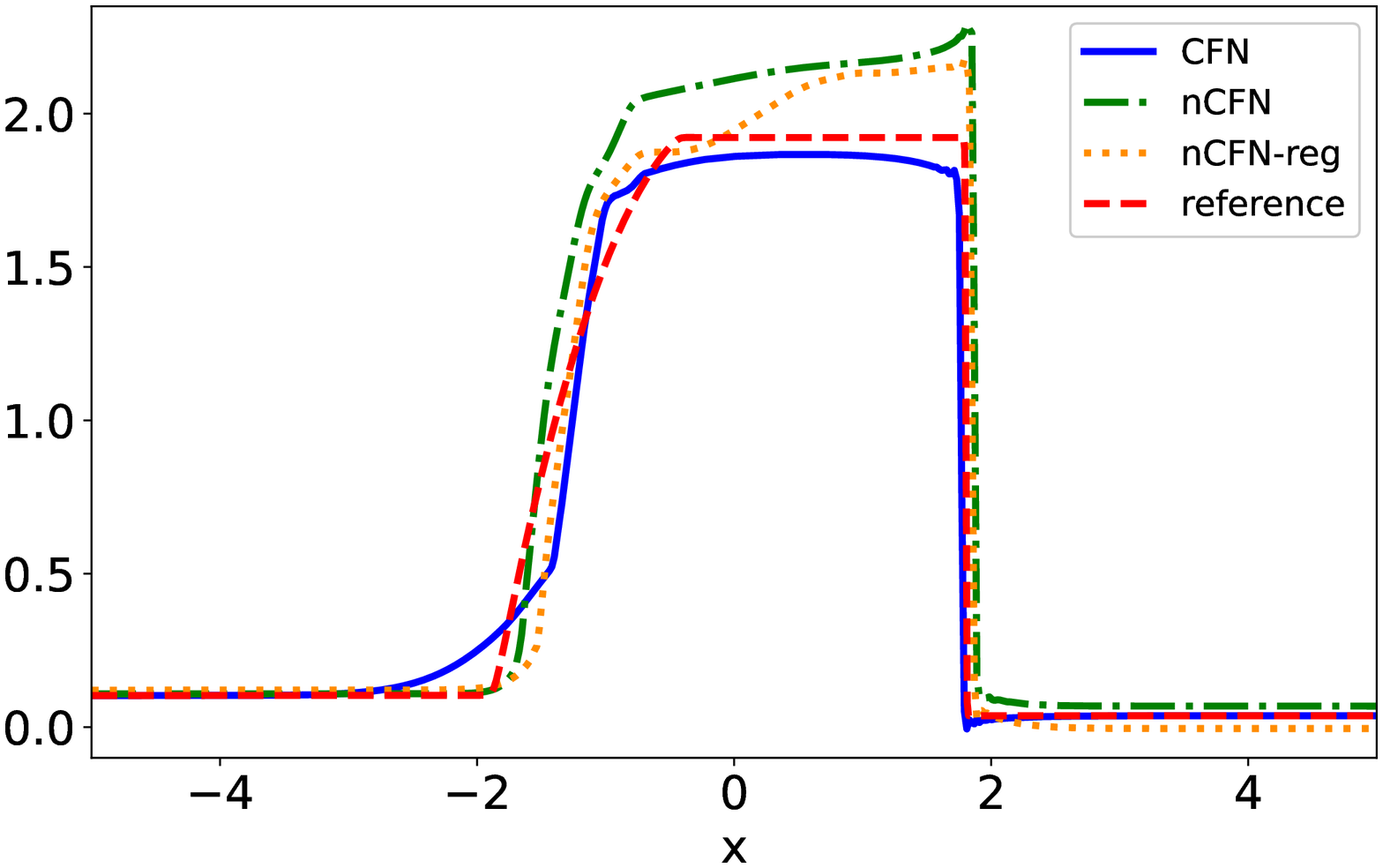}
	\caption{10\% noise}
	\label{fig:SWEsnr10case2_hv}
\end{subfigure}%
\caption{Comparison of the reference solution for height $h$ (top) and  momentum $hv$ (bottom) in Example \ref{ex:SWE} and the trained DNN model predictions at $t=1$ for dense ($N=512$) and noisy observations.}
\label{fig:SWE_height_noisy}
\end{figure}


We again omit figures comparing the discrete conserved quantity remainder, \eqref{eq:conserve_u}, of each method for  Cases  II and III since the methods all generate the same general behavior pattern as what is shown for Case I in \figref{fig:SWE_conserv}.


\subsection{Euler equation}
\label{sec:euler}
As a final example we  consider the Euler equations for gas dynamics, specifically the Shu-Osher problem \cite{SHU1988439}.  The problem is challenging since the resulting shock wave impacts a sinusoidally-varying density field yielding more complex structures than apparent in Examples \ref{ex:burgers} and \ref{ex:SWE}.

\begin{example}\label{ex:euler}
	Consider the system of equations for $t > 0$ given by
	\begin{align}
		\rho_t + (\rho u)_x &= 0,  \nonumber\\
		(\rho u)_t + (\rho u^2 + p)_x & = 0, \nonumber\\
		E_t + (u(E+p))_x &= 0,
	\end{align}\label{eq:Euler}
in the domain $(-5,5)$. We assume no flux boundary conditions and initial conditions given by\footnote{{The usual Shu-Osher problem does not use  $\rho(x,0) = 1 + \varepsilon \sin(5x)  e^{-(x-x_1)^4}$ for $x$ in the right part of the domain.  We include this term to ``flatten'' the solution at the boundary so that we can apply \eqref{eq:BCzeroth} without introducing an artificial boundary layer.}}
	\begin{align*}
			\rho(x,0)&=\begin{cases}
		\rho_l, & \text{if $x\leq x_0$},\\
		1 + \varepsilon \sin(5x), & \text{if $x_0 < x \le x_1$},\\
		1 + \varepsilon \sin(5x)  e^{-(x-x_1)^4} & \text{otherwise},\\
		\end{cases}\quad
		u(x,0)=\begin{cases}
		u_l, \quad \text{if $x\leq x_0$},\\
		0, \quad \text{otherwise},\\
		\end{cases}\\
		p(x,0) &=\begin{cases}
			p_l, \quad \text{if $x\leq x_0$},\\
			p_r, \quad \text{otherwise},
		\end{cases}\quad
		E(x,0)  = \frac{p_0}{\gamma - 1} + \frac{1}{2} \rho(x,0) u(x,0)^2.
	\end{align*}
The parameters are given by
\begin{align}
\rho_l &\sim U[\hat{\rho}_l(1-\epsilon),\quad \hat{\rho}_l(1+\epsilon)], \quad\hat{\rho}_l=3.857135,\nonumber\\
\varepsilon &\sim U[\hat{\varepsilon}(1-\epsilon), \quad \hat{\varepsilon}(1+\epsilon)], \quad {\hat{\varepsilon}} = 0.2,\nonumber\\
p_l &\sim U[\hat{p}_l(1-\epsilon), \quad \hat{p}_l(1+\epsilon)], \quad \hat{p}_l=10.33333,\nonumber\\
p_r &\sim U[\hat{p}_r(1-\epsilon), \quad \hat{p}_r(1+\epsilon)], \quad \hat{p}_r=1,\nonumber\\
u_l &\sim U[\hat{u}_l(1-\epsilon), \quad \hat{u}_l(1+\epsilon)], \quad \hat{u}_l=2.62936,\nonumber\\
x_0 &\sim U[\hat{x}_0(1-\epsilon), \quad \hat{x}_0(1+\epsilon)], \quad \hat{x}_0 = -4,
\label{eq:init_dist}
\end{align}
with $\epsilon = .1$, $x_1 = 3.29867$ and $\gamma = 1.4$. We note that  $\hat{\rho}$, $\hat{p}_l$, $\hat{u}$ are the same values as those used in the CLAWPACK Shu-Osher example.
\end{example}

The $k = 1,\dots,N_{traj}$ training sets are generated by solving \eqref{eq:Euler} using CLAWPACK (HLLE Riemann Solver) for initial conditions given by
	\begin{align*}
	\rho^{(k)}(x,0)&=\begin{cases}
	\rho_l^{(k)}, & \text{if $x\leq x_0^{(k)} $},\\
	1 + \varepsilon^{(k)} \sin(5x), & \text{if $x_0^{(k)} < x \leq x_1$},\\
	1 + \varepsilon^{(k)} \sin(5x)  e^{-(x-x_1)^4}, & \text{otherwise},\\
	\end{cases}\\
	u^{(k)}(x,0)&=\begin{cases}
		u_l^{(k)}, \quad \text{if $x\leq x_0^{(k)}$},\\
		0, \quad \text{otherwise},\\
	\end{cases}
 p^{(k)}(x,0) =\begin{cases}
		p_l^{(k)}, \quad \text{if $x\leq x_0^{(k)}$},\\
		p_r^{(k)}, \quad \text{otherwise},
	\end{cases}\\
	E^{(k)}(x,0) &= \frac{p^{(k)}_0}{\gamma - 1} + \frac{1}{2} \rho^{(k)}(x,0) u^{(k)}(x,0)^2.
\end{align*}
The corresponding parameters are given in \eqref{eq:init_dist} (written without the superscript $k$) and
the boundary conditions are imposed using \eqref{eq:BCzeroth} in all experiments.
The reference solution is obtained using CLAWPACK with $N=1024$ so that $\Delta x = \frac{10}{1024}$.

To train over a longer period of time without increasing the computational cost we once again employ the same sub-sampling technique used in Example \ref{ex:burgers} to generate training data from observed snapshots of the solution on an extended domain. As before we set $M=300$ as the extended length of each trajectory. The snapshots of the solution, \eqref{setup:data}, are obtained via CLAWPACK for each of the $N_{traj} = 300$ trajectories at times  $t = m\Delta t$, $m = 1,\dots, M$, where $\Delta t = 0.002$ (chosen to satisfy the CFL condition). Each sub-sampled trajectory of length $L$ is then built consecutively from the snapshot solutions. That is, each trajectory is comprised of the solutions in \eqref{setup:data} at sequential times $t_0^{(k)} + l \Delta t$, $l = 1,\dots ,L$.

\subsubsection*{Case I: Dense and noise-free observations}
\begin{figure}[h!]
	\centering
	\begin{subfigure}[b]{0.24\textwidth}
		\includegraphics[width=\textwidth]{%
			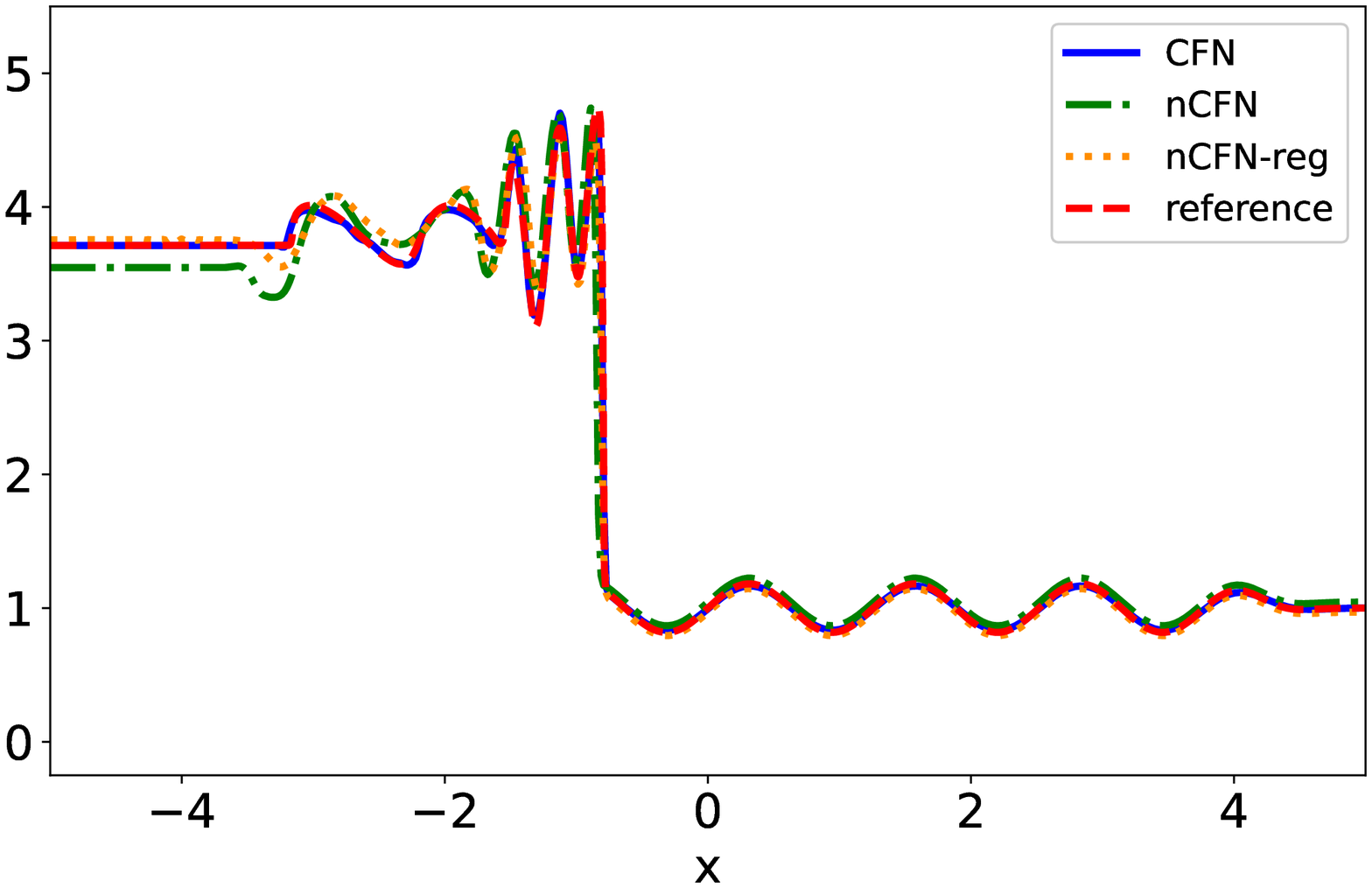}
		\caption{$\rho, t=0.8$}
		\label{fig:euler_ideal_rho_t1}
	\end{subfigure}%
	~
	\begin{subfigure}[b]{0.24\textwidth}
		\includegraphics[width=\textwidth]{%
			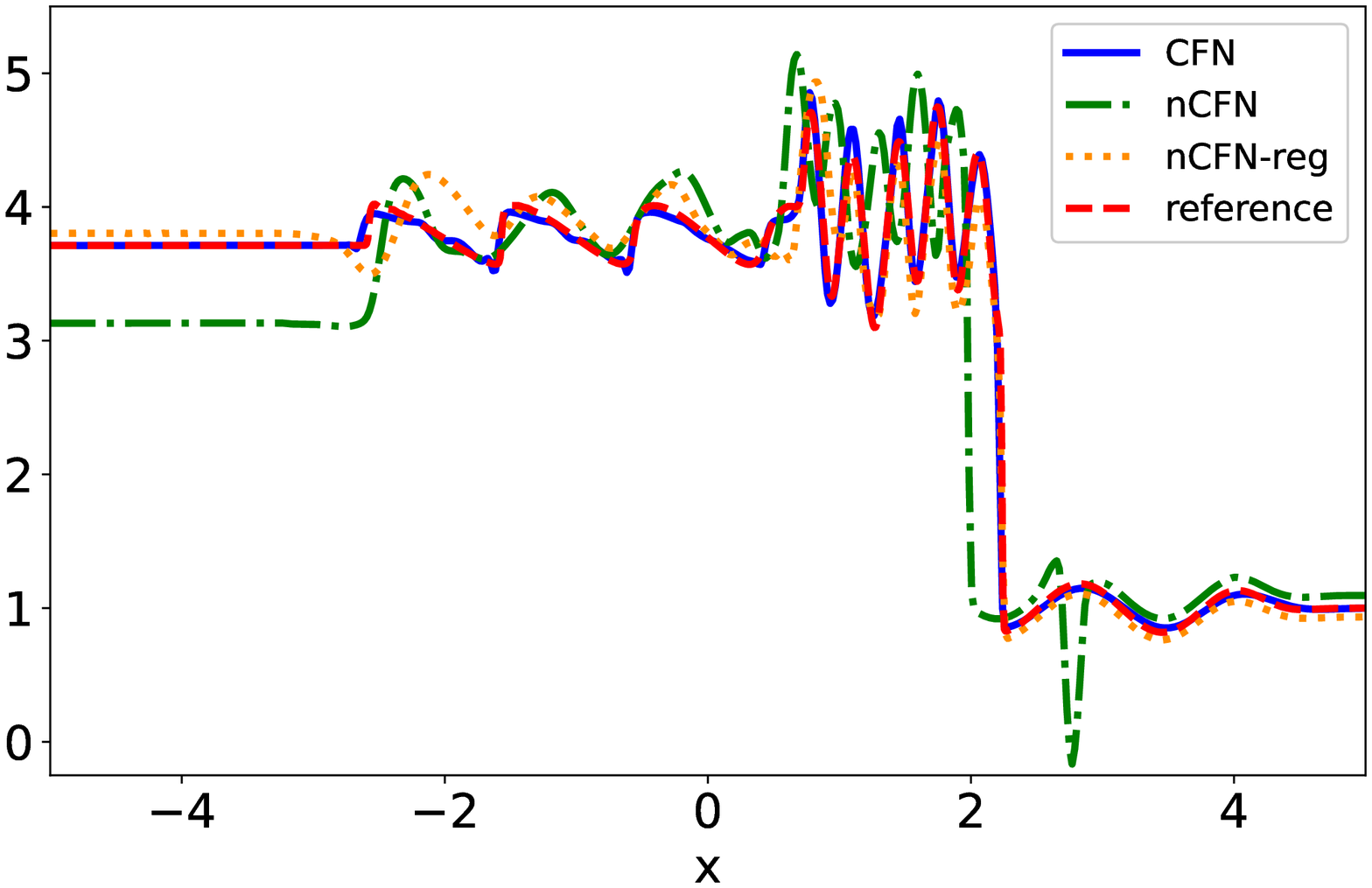}
		\caption{$\rho, t=1.6$}
		\label{fig:euler_ideal_rho_t2}
	\end{subfigure}%
	\begin{subfigure}[b]{0.24\textwidth}
		\includegraphics[width=\textwidth]{%
			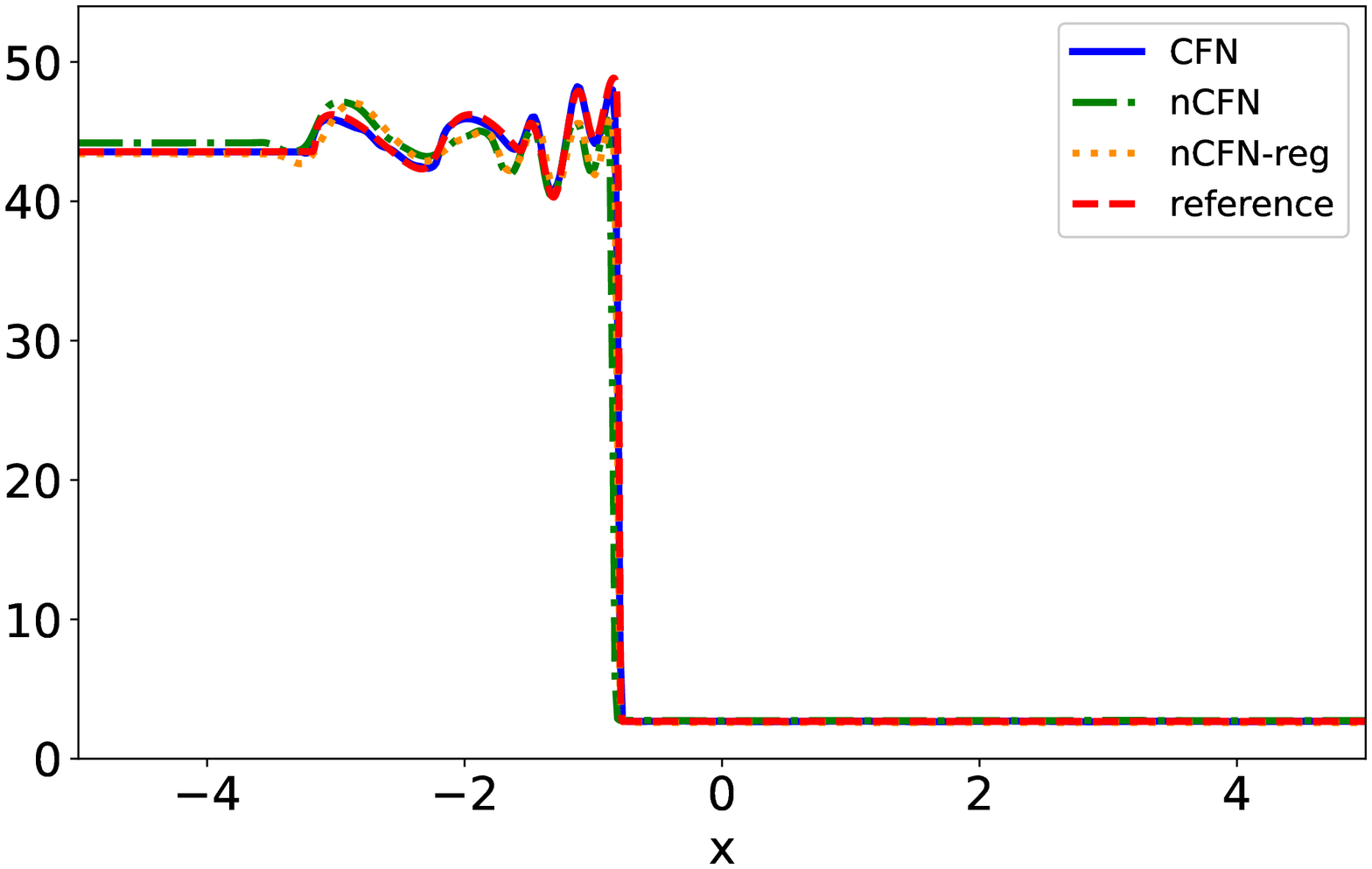}
		\caption{$E, t=0.8$}
		\label{fig:euler_ideal_E_t1}
	\end{subfigure}%
	~
	\begin{subfigure}[b]{0.24\textwidth}
	\includegraphics[width=\textwidth]{%
		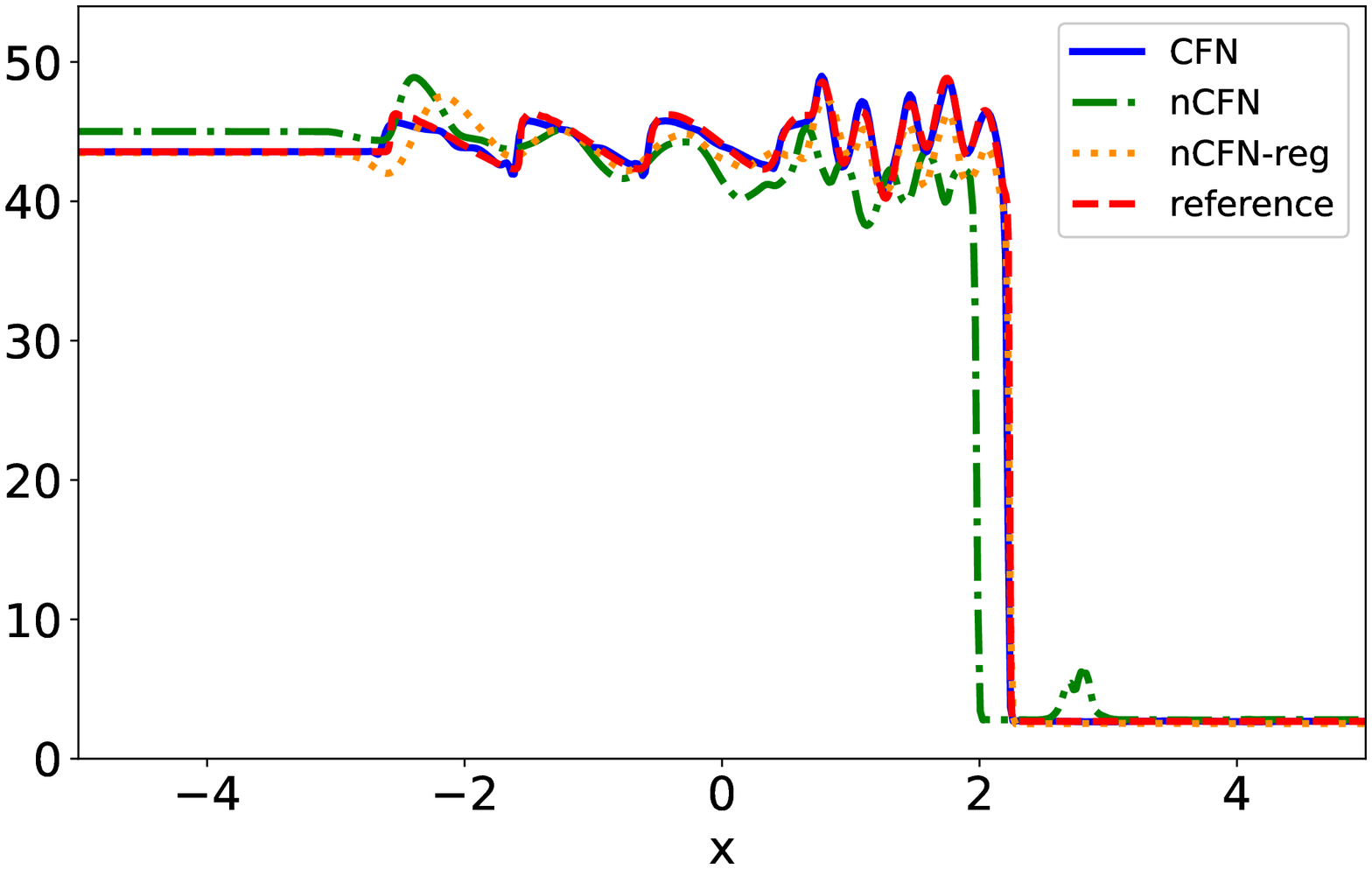}
	\caption{$E, t=1.6$}
	\label{fig:euler_ideal_E_t2}
\end{subfigure}%
	\caption{Comparison of the reference density $\rho$ and energy $E$ solutions to Example \ref{ex:euler} and the trained DNN model predictions at different times for dense ($N=512$) and noise-free observations.}
	\label{fig:euler_ideal}
\end{figure}

We first consider an idealized environment for which the observations are dense and noise-free. Specifically we choose $N = 512$, yielding $\Delta x = \frac{10}{512}$, so that our solution is well-resolved.
\figref{fig:euler_ideal} shows the solutions $\rho$ and $E$ at times $t=0.8$ and $t=1.6$, both of which are outside of training time domain $[0,0.6]$. Observe that as the shock wave interacts with the density field, the solution exhibits oscillations to the left side of the shock front. It is evident that only the CFN network produces a solution that captures the oscillatory features of the solution. By contrast, the nCFN solution exhibits significant errors with non-physical oscillations to the right of the shock.  It moreover produces the wrong shock front location at $t = 1.6$.  The results for nCFN-reg are somewhat improved, but still do not accurately capture the dynamics of the system.

\begin{figure}[h!]
	\centering
	\begin{subfigure}[b]{0.30\textwidth}
		\includegraphics[width=\textwidth]{%
			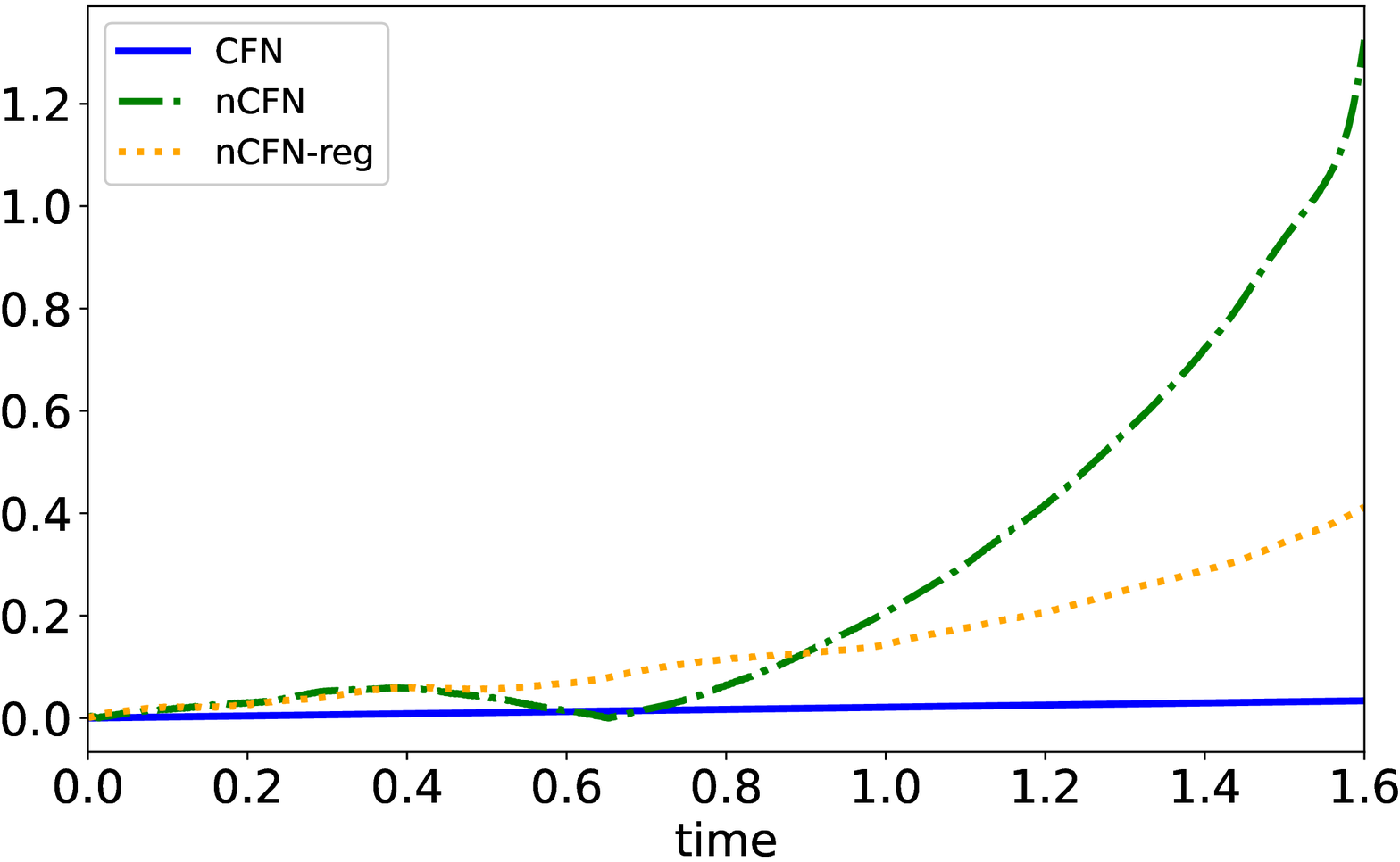}
		\caption{ $C(\boldsymbol{\rho})$}
		\label{fig:euler_conserv_density}
	\end{subfigure}%
	~
	\begin{subfigure}[b]{0.30\textwidth}
		\includegraphics[width=\textwidth]{%
			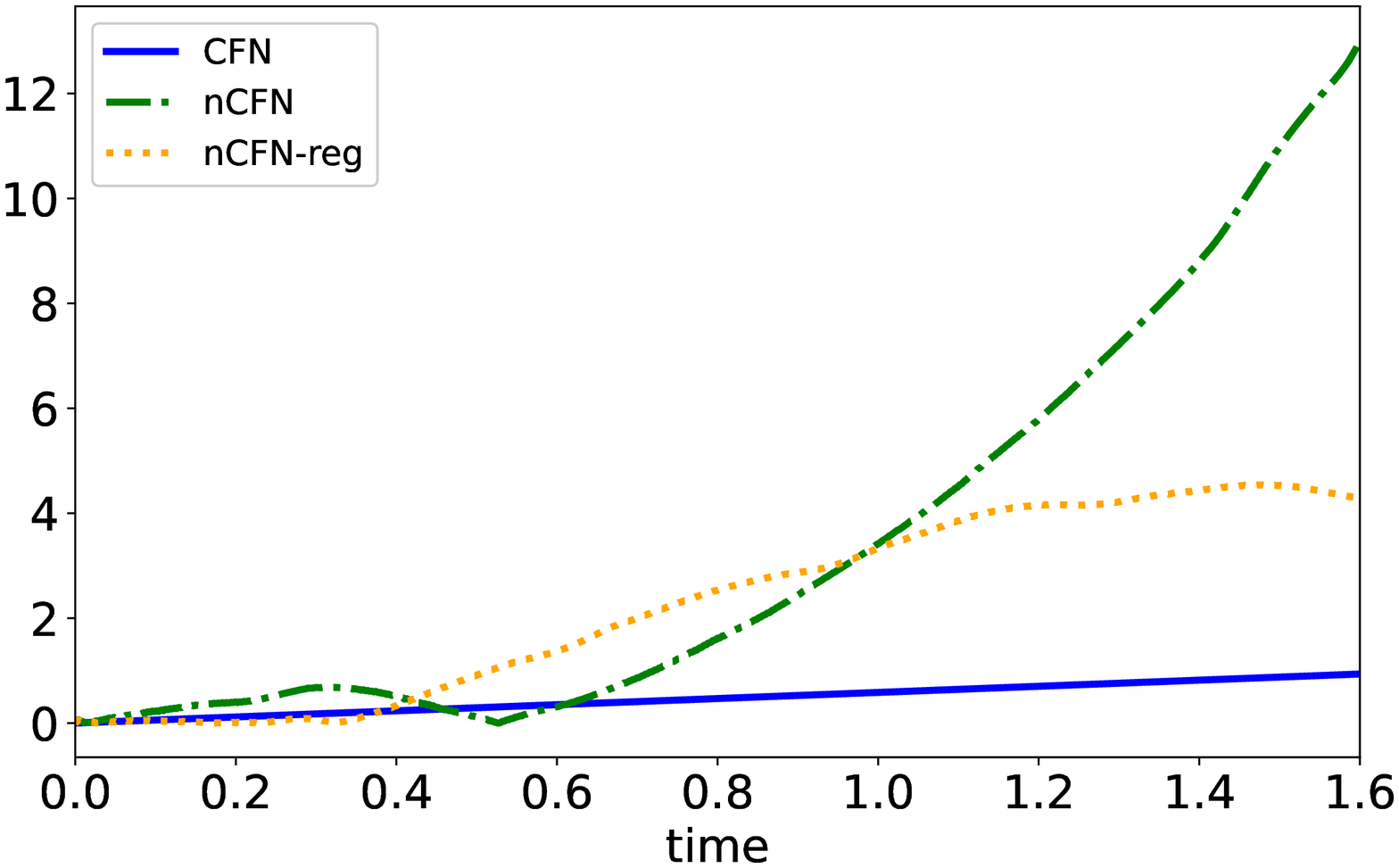}
		\caption{ $C(\mathbf{E})$}
		\label{fig:euler_conserv_energy}
	\end{subfigure}%
	\caption{Discrete conserved quantity remainder given in \eqref{eq:conserve_u}  of each method for Example  \ref{ex:euler}.  (a) $C(\boldsymbol{\rho})$ and (b) $C(\mathbf{E})$.}
	\label{fig:euler_conserv}
\end{figure}

\figref{fig:euler_conserv} confirms our observations in  \figref{fig:euler_ideal}.  Specifically, we see that none of the methods are conservative, with the error increasing more rapidly in the nCFN and the nCFN-reg cases.  The error corresponding to the CFN appears to grow linearly with time, indicating long term stability when considering classical numerical conservation laws analysis.

\subsubsection*{Case II: Sparse and noise-free observations}
\begin{figure}[!h]
	\centering
		\begin{subfigure}[b]{0.24\textwidth}
		\includegraphics[width=\textwidth]{%
			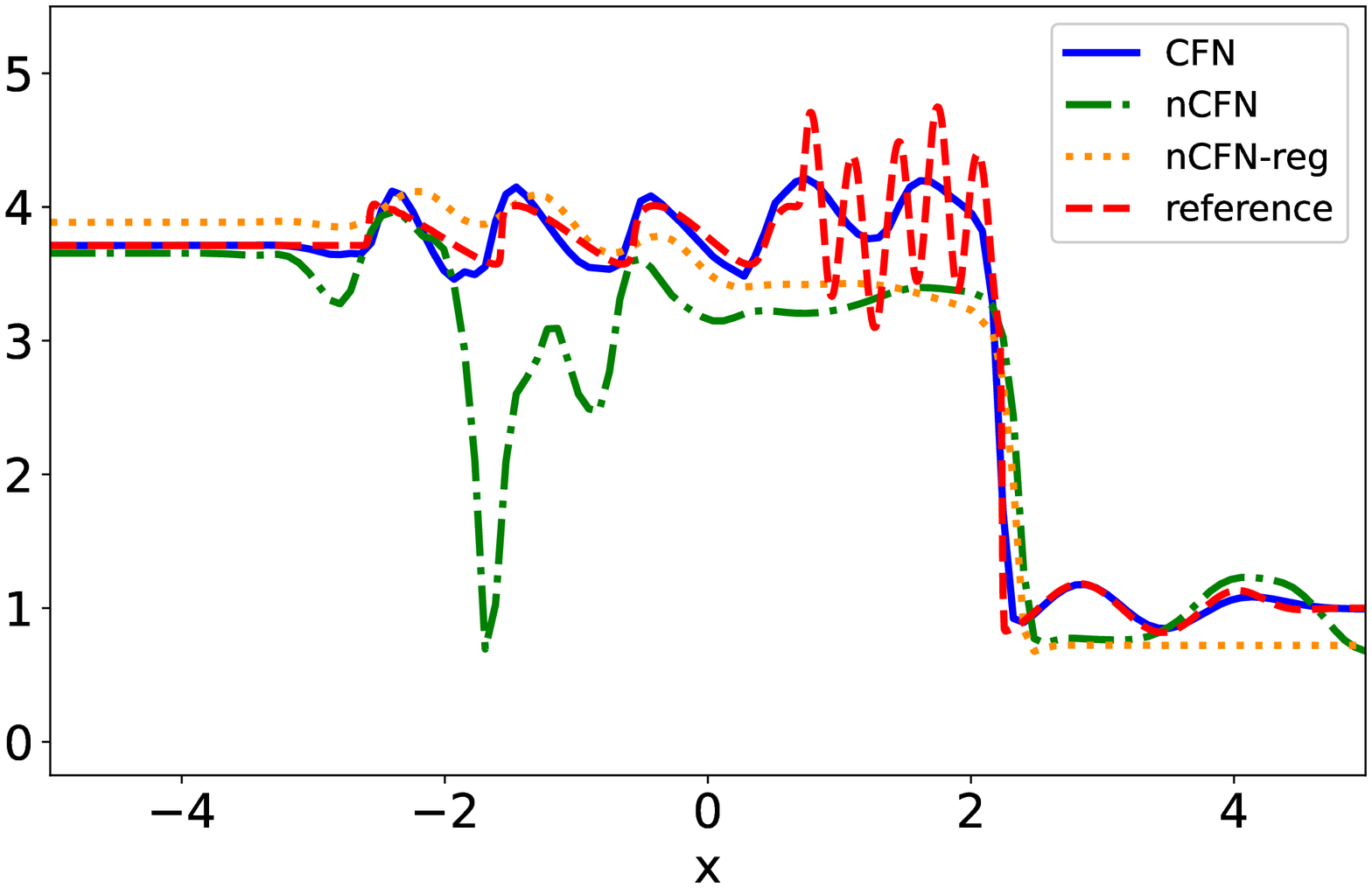}
		\caption{$\rho, N=128$}
		\label{fig:euler_sparse_rho_N128}
	\end{subfigure}%
~
	\begin{subfigure}[b]{0.24\textwidth}
		\includegraphics[width=\textwidth]{%
			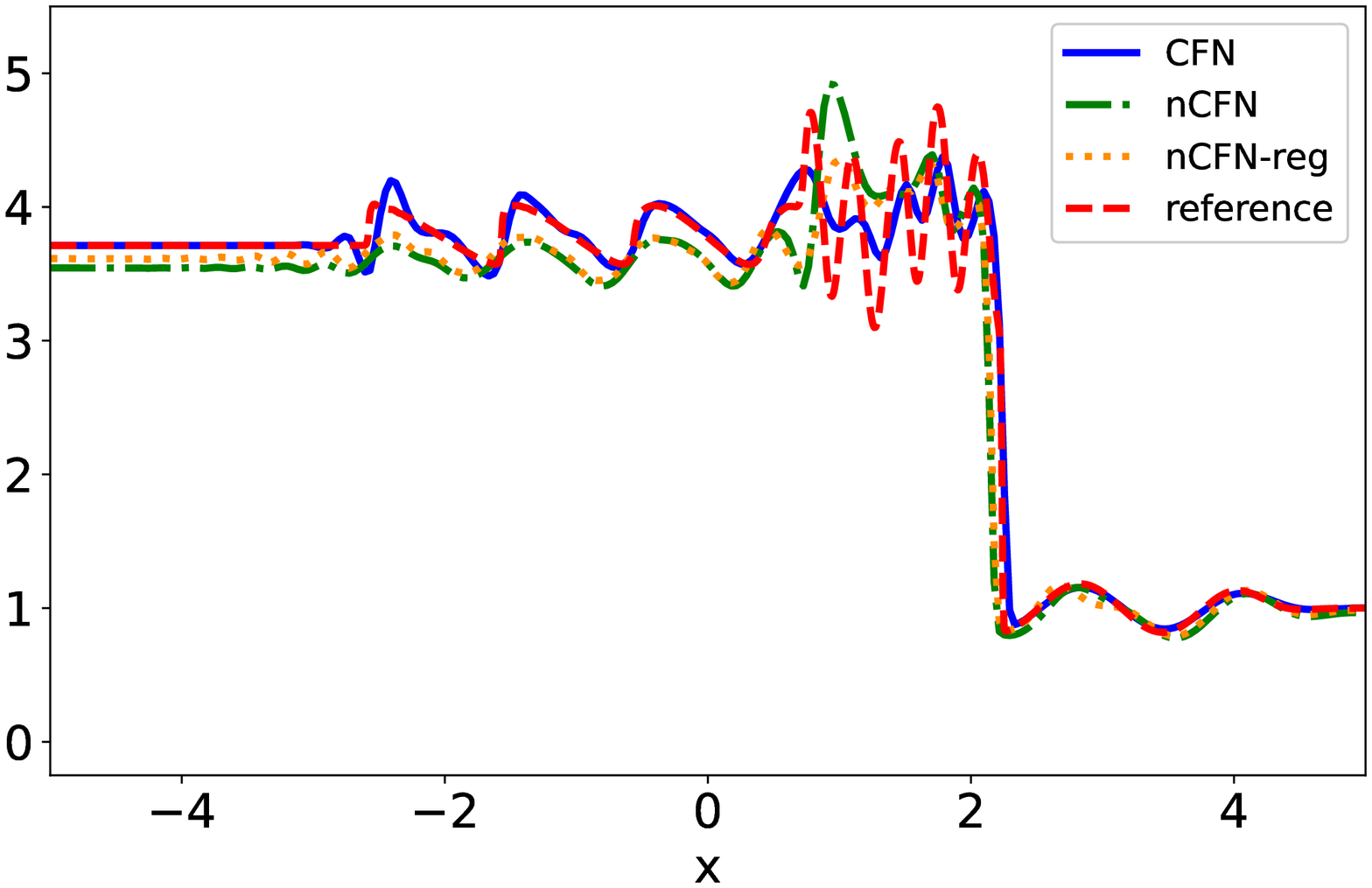}
		\caption{$\rho, N=256$}
		\label{fig:euler_sparse_rho_N256}
	\end{subfigure}%
\begin{subfigure}[b]{0.24\textwidth}
	\includegraphics[width=\textwidth]{%
		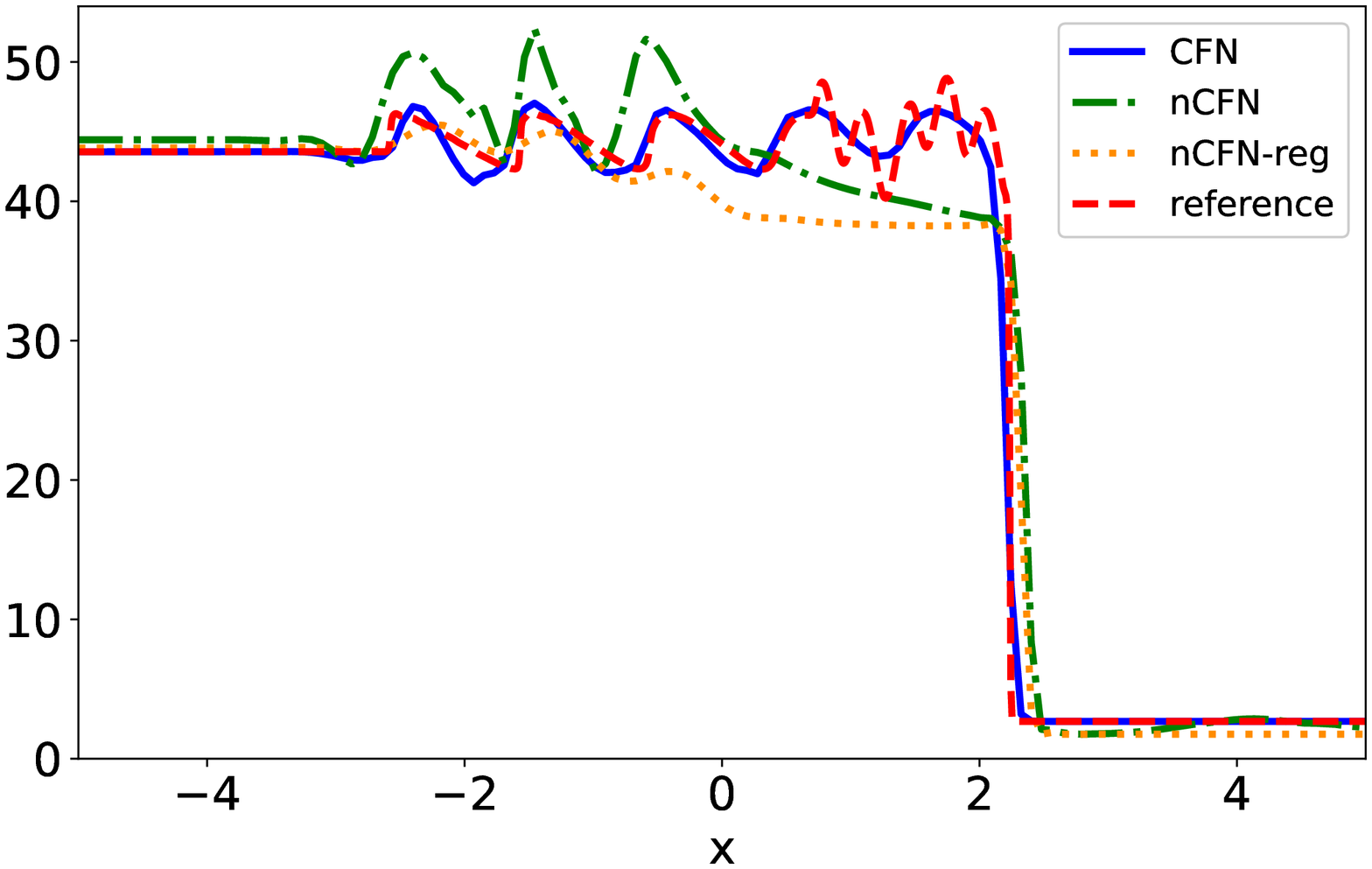}
	\caption{$E, N=128$}
	\label{fig:euler_sparse_E_N128}
\end{subfigure}%
~
\begin{subfigure}[b]{0.24\textwidth}
	\includegraphics[width=\textwidth]{%
		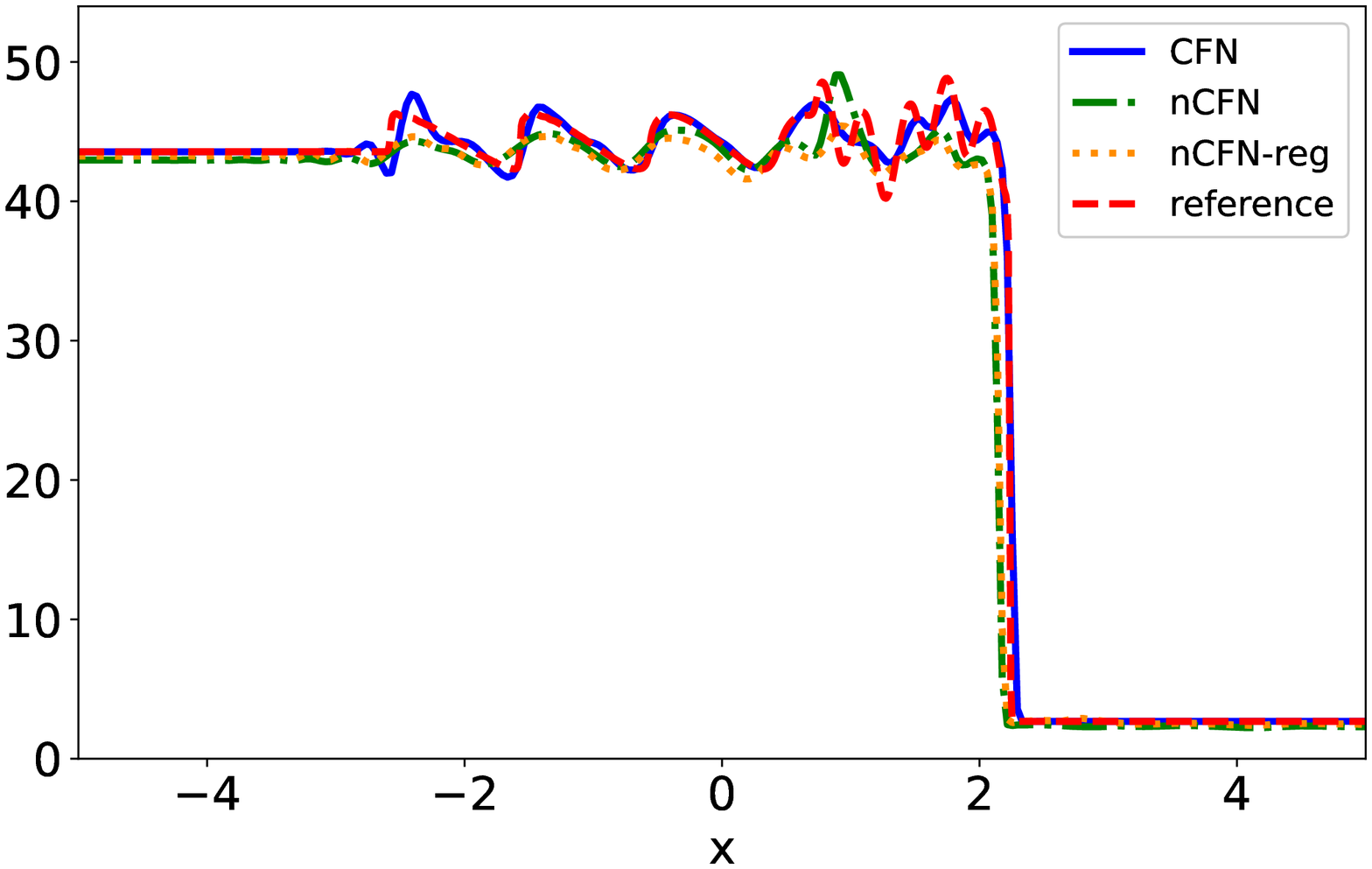}
	\caption{$E, N=256$}
	\label{fig:euler_sparse_E_N256}
\end{subfigure}%
	\caption{Comparison of the reference density $\rho$ and energy $E$ solutions to Example \ref{ex:euler} and the trained DNN model predictions at $t=1.6$ for sparse and noise-free observations.}
	\label{fig:euler_sparse}
\end{figure}

We now consider more sparsely observed data by choosing $N=128$ and $N=256$ {shown in \figref{fig:euler_sparse}}. Given the results in the idealized environment, it is not surprising that neither the nCFN or nCFN-reg is able to capture the dynamics of Example \ref{ex:euler} in the sparse observation case. While some solution details are lost, and there is noticeable error in the location of the shock front, it is evident that the CFN network still provides qualitative structure commensurate with the given resolution. 

\subsubsection*{Case III: Dense and noisy observations}
In this case the training data are
\begin{eqnarray}
	\label{eq:euler_noise_dense}
	\tilde{\boldsymbol{\rho}}^{(k)}(t_l) &=& \boldsymbol{\rho}^{(k)}(t_l) + \boldsymbol{\epsilon}_l^{(k)}, \nonumber\\
	\widetilde{(\boldsymbol{\rho}^{(k)}\odot\mathbf{u}^{(k)})}(t_l) &=& (\boldsymbol{\rho}^{(k)}\odot\mathbf{u}^{(k)})(t_l) + \boldsymbol{\eta}_l^{(k)}, \nonumber\\
	\tilde{\mathbf{E}}^{(k)}(t_l) &=& \mathbf{E}^{(k)}(t_l) + \boldsymbol{\delta}_l^{(k)},
\end{eqnarray}
for $l=1,\dots,L$  and $k = 1,\dots, N_{traj}$.
The added noise $\boldsymbol{\epsilon}_l^{(k)}$, $\boldsymbol{\eta}_{l}^{(k)}$, and $\boldsymbol{\delta}_l^{(k)}$  are  i.i.d. Gaussian with zero mean and variance determined by various noise levels \eqref{eq:noiselevel}.
We again consider the same noise levels, 100\%, 50\%, 20\%, and 10\%.

\begin{figure}[!h]
	\centering
	\begin{subfigure}[b]{0.24\textwidth}
		\includegraphics[width=\textwidth]{%
			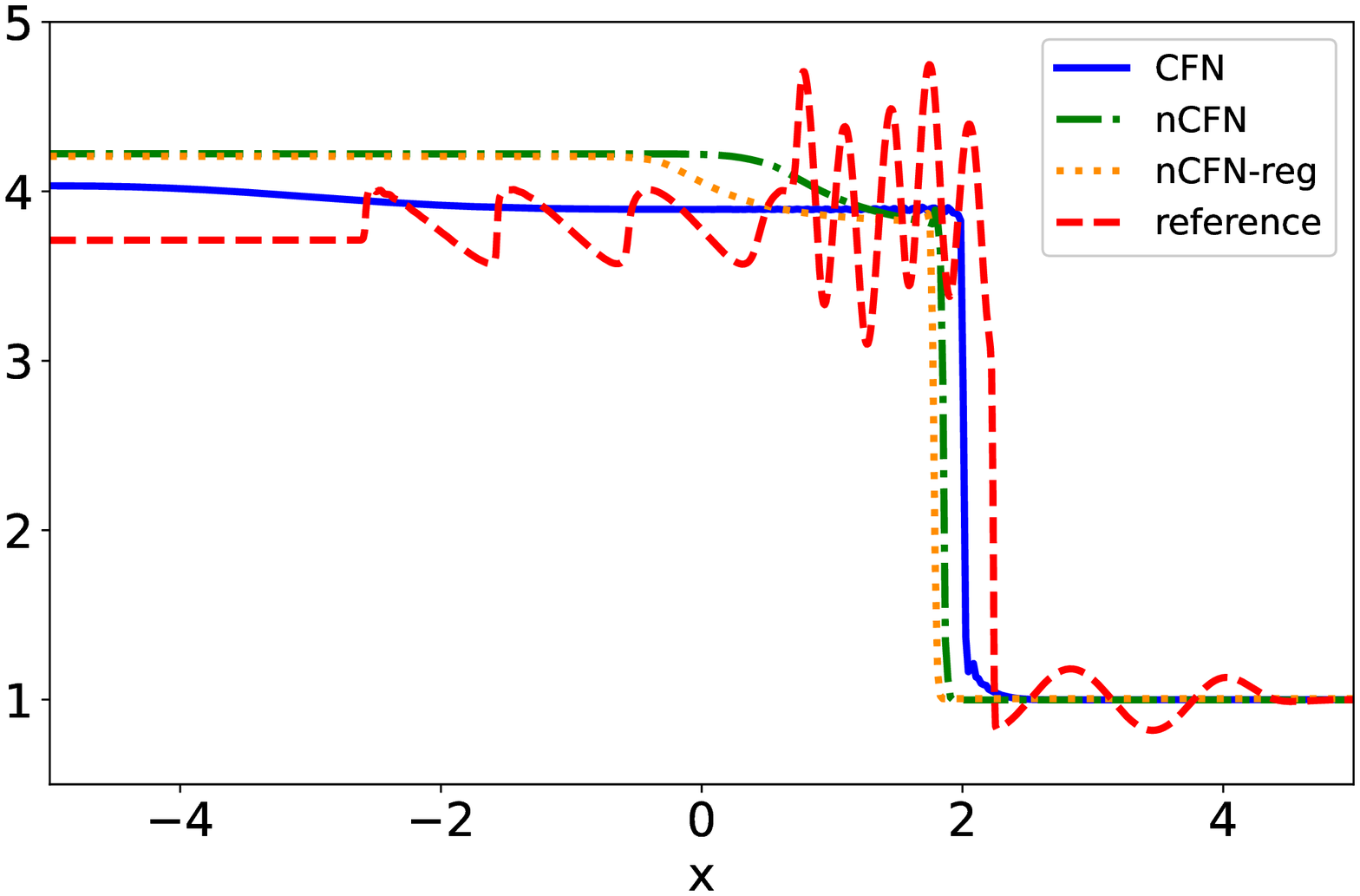}
		\caption{$\rho$, 100\% noise}
		\label{fig:euler_snr1_rho}
	\end{subfigure}%
	~
	\begin{subfigure}[b]{0.24\textwidth}
	\includegraphics[width=\textwidth]{%
		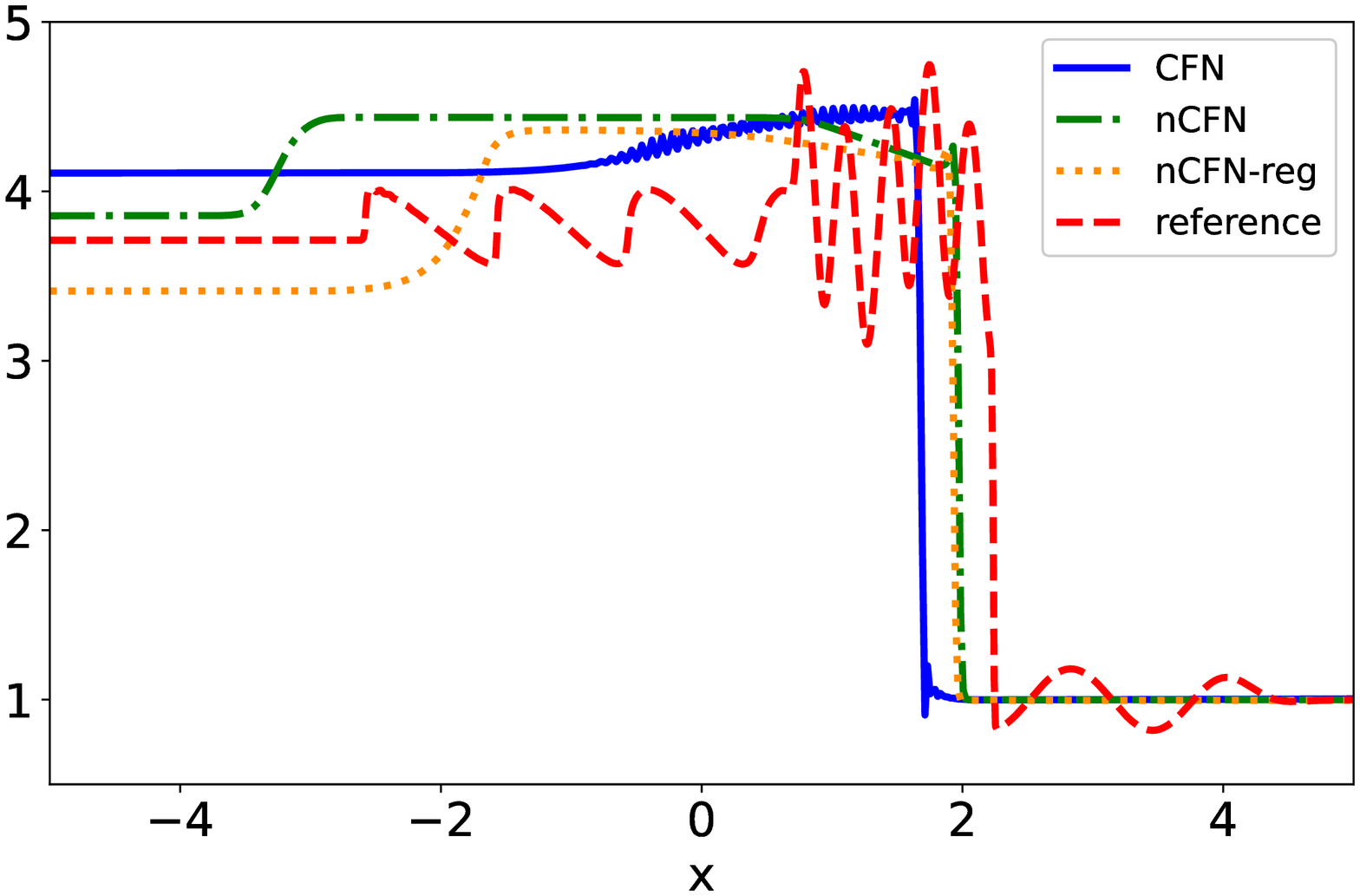}
	\caption{$\rho$, 50\% noise}
	\label{fig:euler_snr2_rho}
	\end{subfigure}%
		\begin{subfigure}[b]{0.24\textwidth}
		\includegraphics[width=\textwidth]{%
			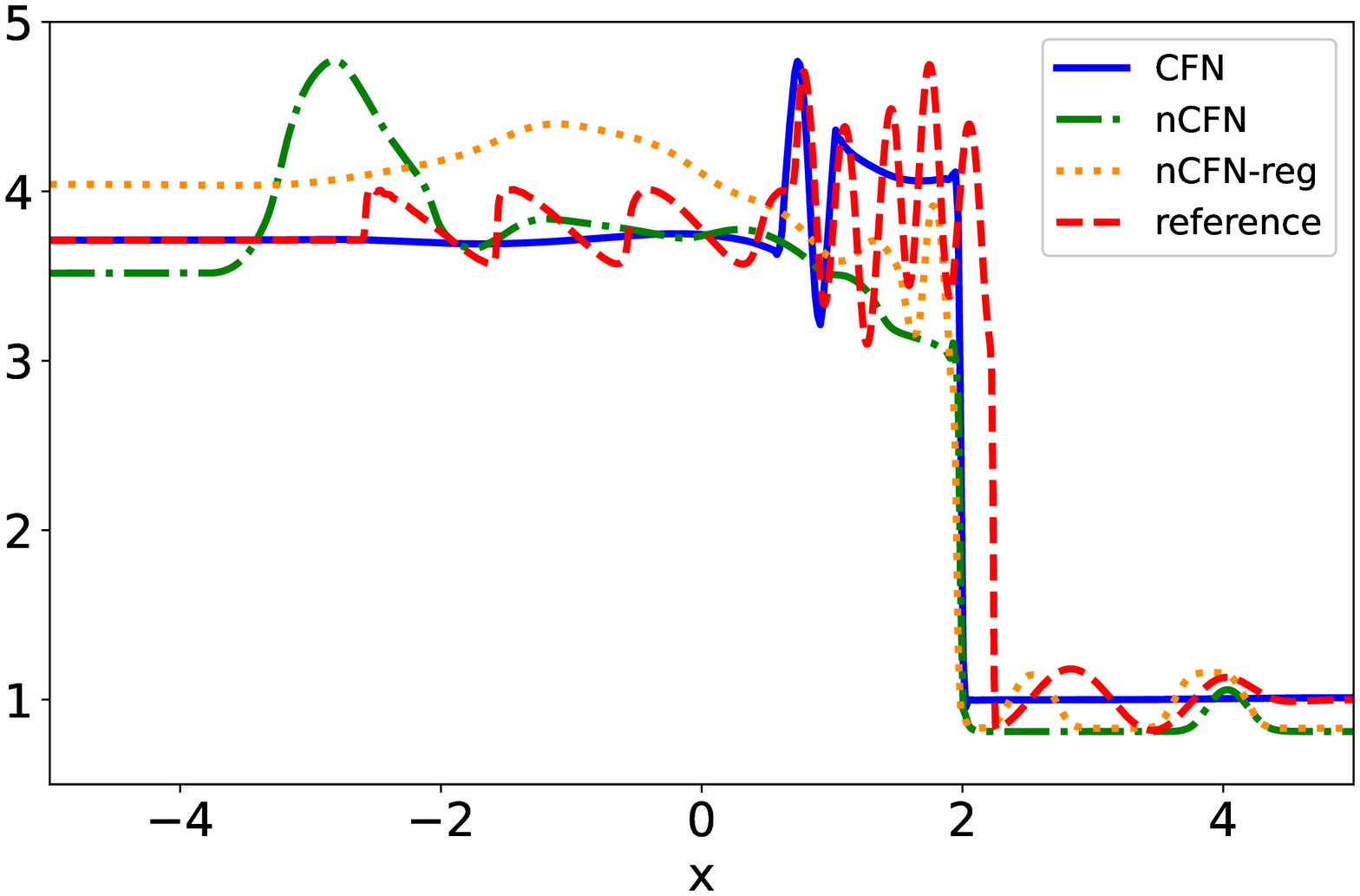}
		\caption{$\rho$, 20\% noise}
		\label{fig:euler_snr5_rho}
	\end{subfigure}%
	\begin{subfigure}[b]{0.24\textwidth}
	\includegraphics[width=\textwidth]{%
		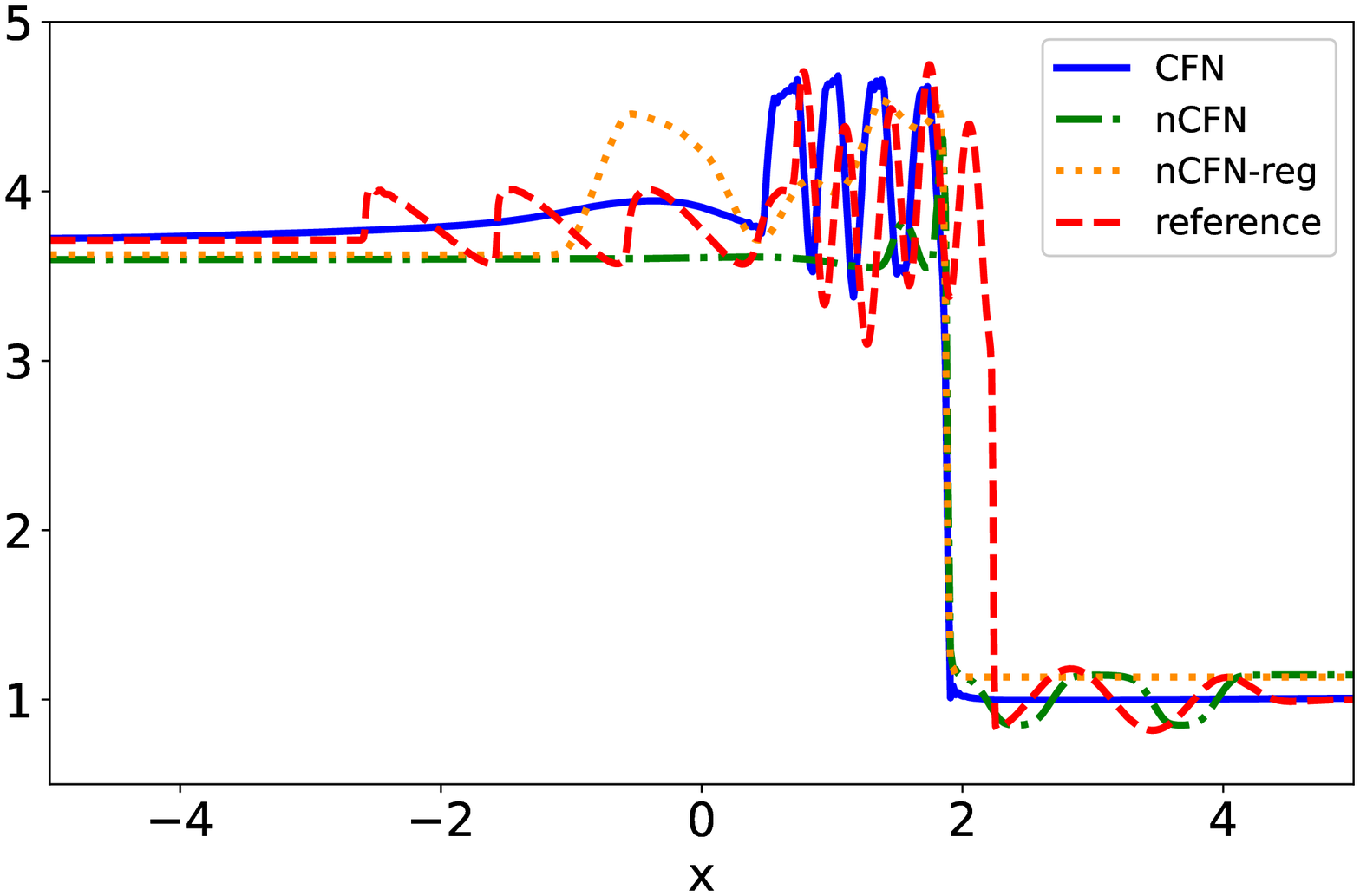}
	\caption{$\rho$, 10\% noise}
	\label{fig:euler_snr10_rho}
\end{subfigure}%
\\
	\begin{subfigure}[b]{0.24\textwidth}
		\includegraphics[width=\textwidth]{%
			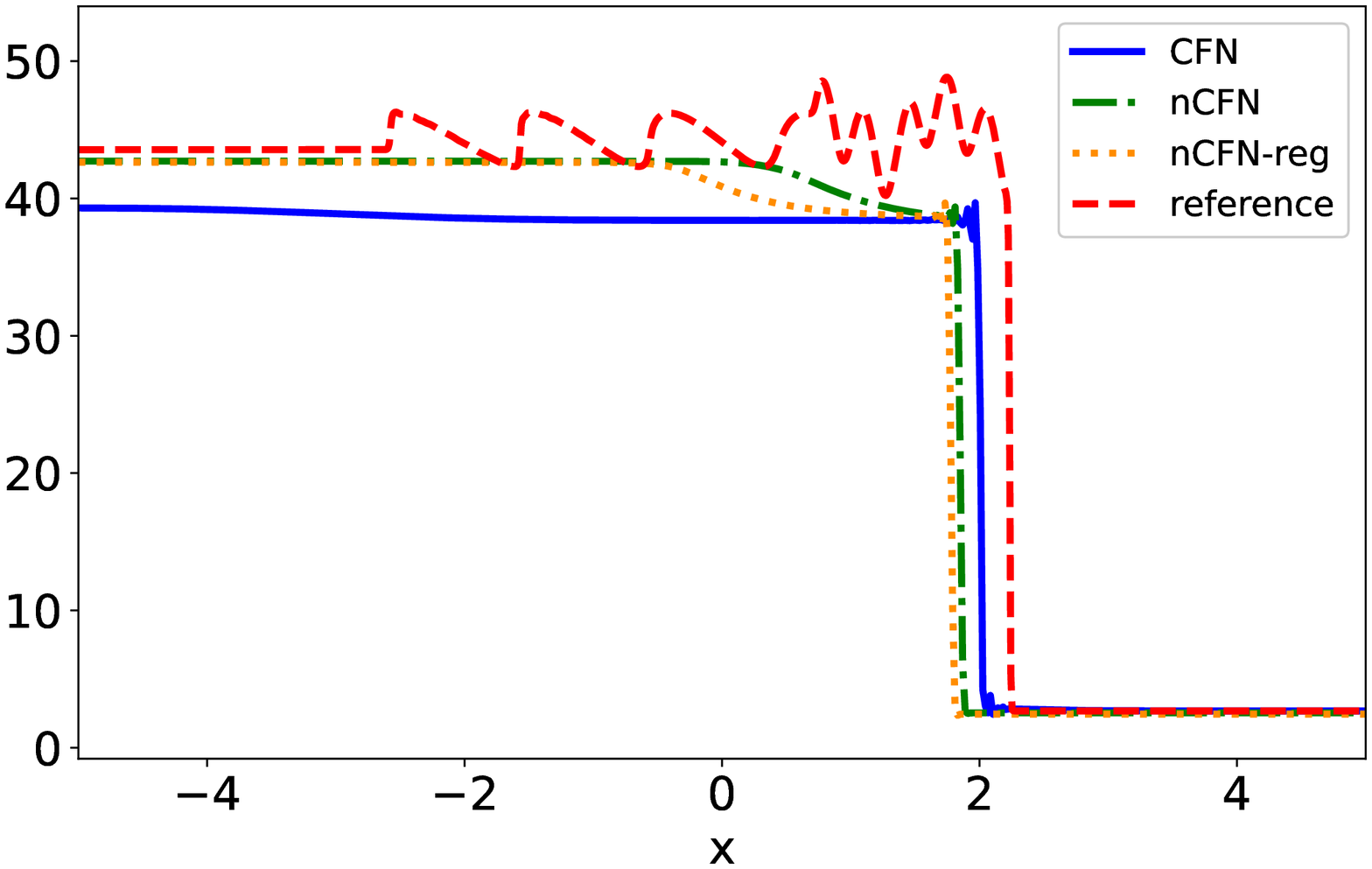}
		\caption{$E$, 100\% noise}
		\label{fig:euler_snr1_E}
	\end{subfigure}%
	~
	\begin{subfigure}[b]{0.24\textwidth}
		\includegraphics[width=\textwidth]{%
			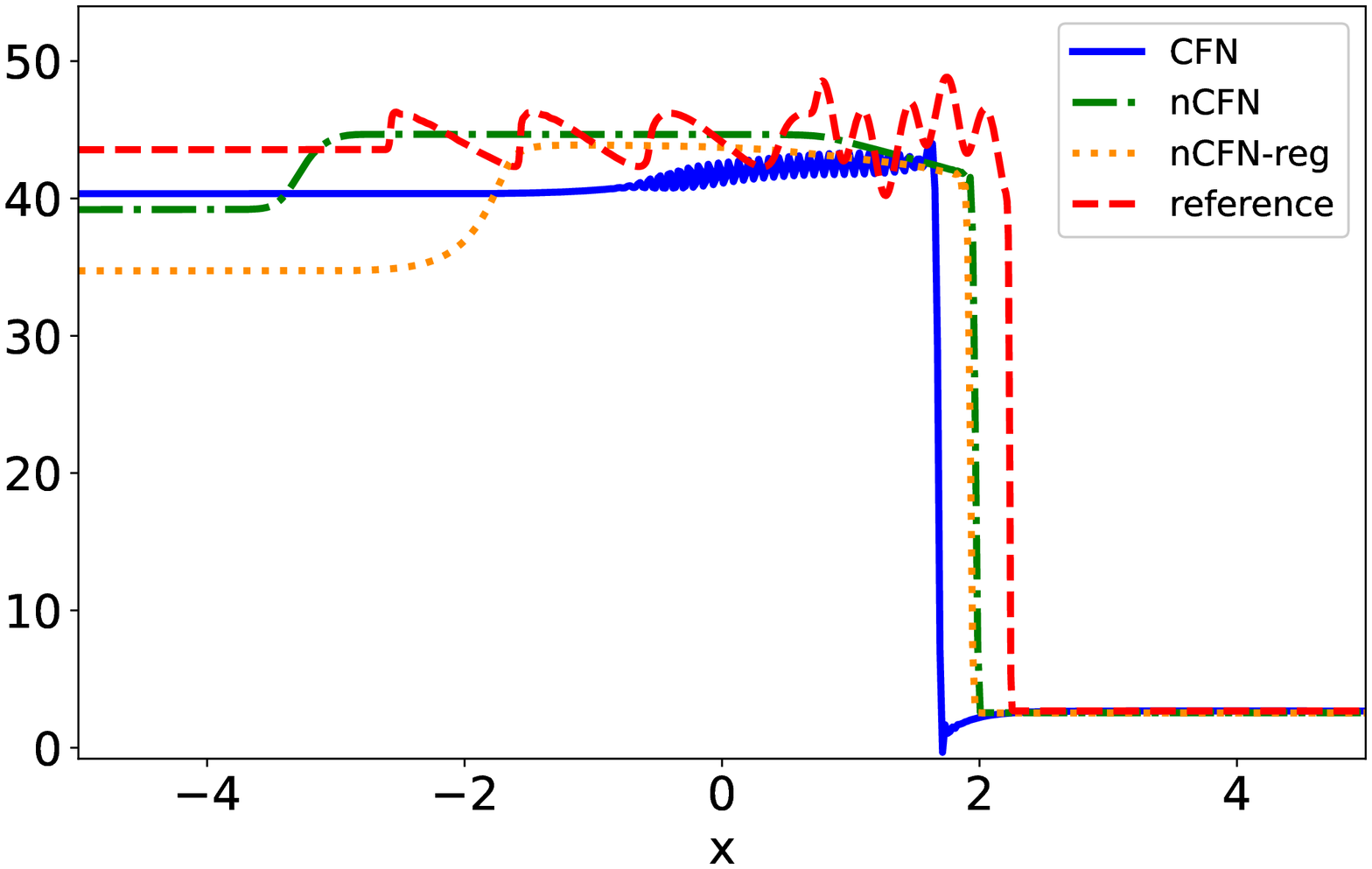}
		\caption{$E$, 50\% noise}
		\label{fig:euler_snr2_E}
	\end{subfigure}%
	\begin{subfigure}[b]{0.24\textwidth}
		\includegraphics[width=\textwidth]{%
			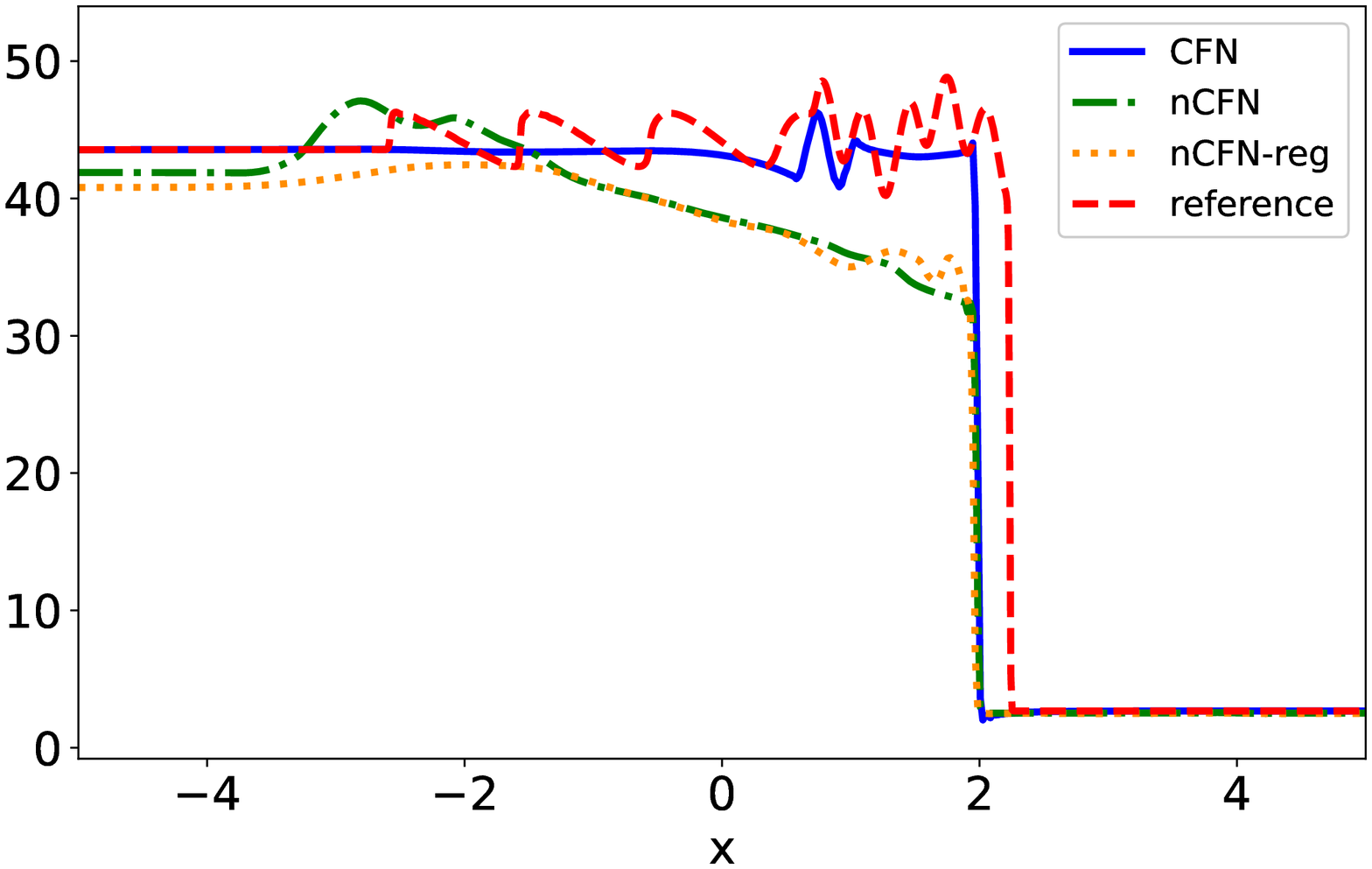}
		\caption{$E$, 20\% noise}
		\label{fig:euler_snr5_E}
	\end{subfigure}%
	\begin{subfigure}[b]{0.24\textwidth}
		\includegraphics[width=\textwidth]{%
			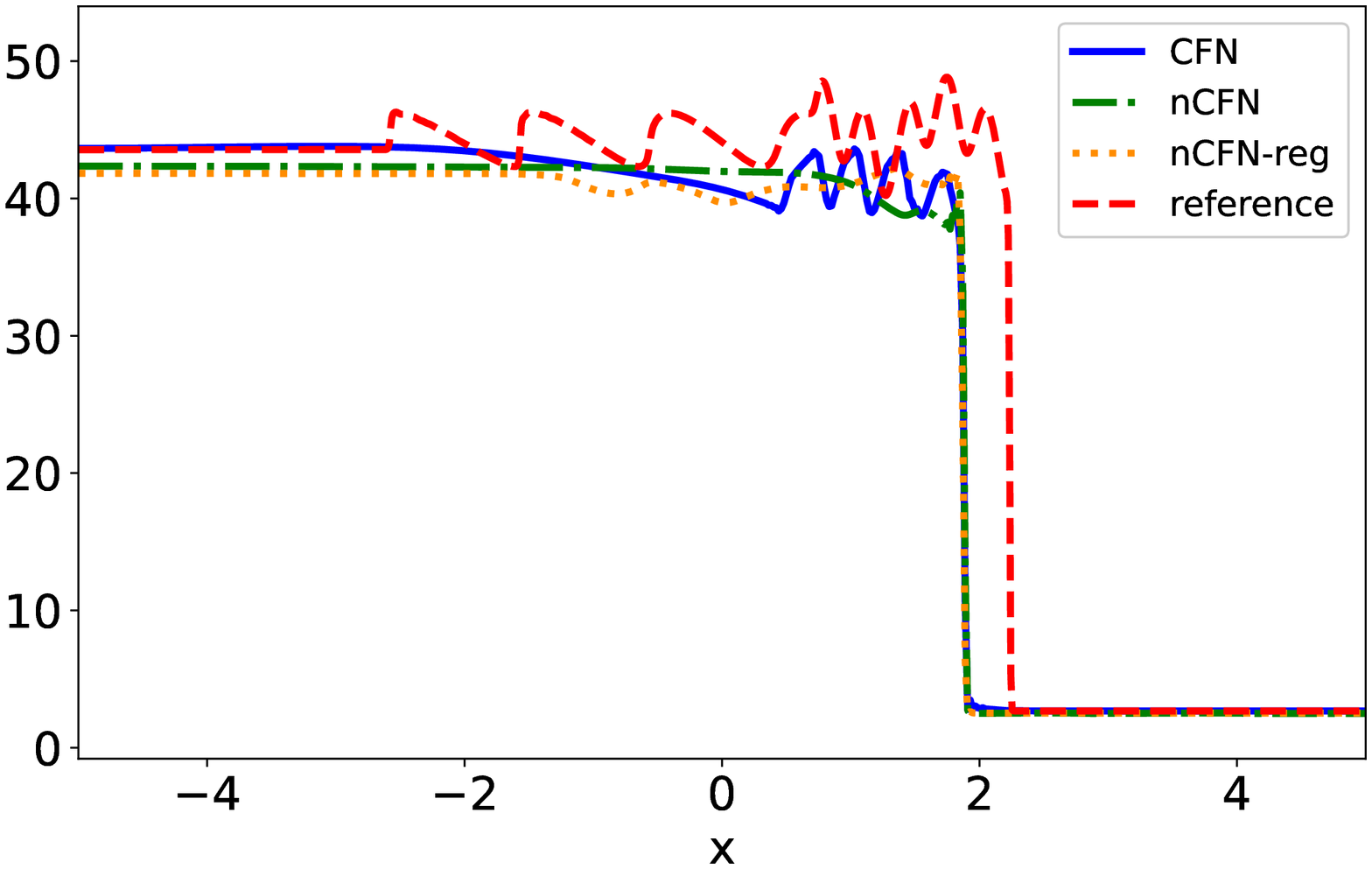}
		\caption{$E$, 10\% noise}
		\label{fig:euler_snr10_E}
	\end{subfigure}%
	\caption{Comparison of the  reference solution of density $\rho$ (top) and Energy $E$ (bottom) in Example \ref{ex:euler} to the trained DNN model predictions at time $t=1.6$ for dense ($N=512$) and noisy observations.}
	\label{fig:euler_rho_snr}
\end{figure}

The solutions for density $\rho$ and energy $E$ are presented in \figref{fig:euler_rho_snr}. We observe similar behavior as was seen for Case III in Example \ref{ex:SWE}.  Specifically, all three methods yield significant diffusion in high noise environments, 100\% and 50\%, and cannot predict the oscillatory structure to the left of the shock front.   Unlike what was observed in Example \ref{ex:SWE}, neither the nCFN nor the nCFN-reg appear to learn the noise-related dynamics, as even in the 10\% noise level case the solutions still appear diffusive.  This is likely because loss function still promotes a diffuse solution as opposed to one that contains noise-related dynamics.  As the noise decreases, the CFN appears to capture some of the oscillatory details in the solution.  In this regard, we again see that the CFN is a more robust network with respect to noise.

We again omit figures comparing the discrete conserved quantity remainder, \eqref{eq:conserve_u}, of each method for  Cases  II and III since the methods all generate the same general behavior pattern as what is shown for Case I in  \figref{fig:euler_conserv}.

\section{Conclusion}
\label{sec:conclusion}
In this investigation we proposed a conservative form network (CFN) to learn the dynamics of unknown hyperbolic systems of conservation laws from observation data. Inspired by classical finite volume methods for hyperbolic conservation laws, our new method employs a neural network to learn the flux function of the unknown system. The predictions using CFN yield the appropriate conserved quantities and also recover the correct physical structures, including the shock speed, even outside the training domain. We validated the effectiveness and robustness of our CFN approach through a series of numerical experiments for three classic examples of  one-dimensional conservation laws. Even in non-ideal environments, our results consistently demonstrate that the CFN outperforms the traditional non-conservative form network (nCFN) and its regularized version (nCFN-reg) in terms of  accuracy, efficiency, and robustness, in particular since it does not require fine-tuning of regularization parameters.  

The current study does not attempt to optimize  model performance for the realistic data cases,  and we will attempt to do this in future investigations.  For the sparse observation environment, the Mori–Zwanzig formalism \cite{mori1965transport,zwanzig1973nonlinear}, for which memory is included in the network, may potentially enhance the overall performance. In the noisy data environment one might consider using a denoising technique such as regularization \cite{bishop1995training, golub1999tikhonov}.
Future investigations will also consider two-dimensional examples with more complex boundary conditions.  Finally, we will also study mixed-form systems, where the CFN may be used for equations representing  conserved quantities within the system.  

 
 	\bibliographystyle{plain}
 	\bibliography{cfn}     
\end{document}